%% file: h3main.tex
\documentclass[12pt,leqno,twoside]{article} 
\usepackage{bezier,amsmath,amssymb,stmaryrd,lscape}
\usepackage[colorlinks, linkcolor=black, citecolor=black, urlcolor=black, hyperfootnotes=true]{hyperref}                

\input xy
\xyoption{all}
\input{rk.sty}
\newcommand{\pp}{\pfk\pfk}
\newcommand{\ulFl}{\ul{\Fl\!}\,}
\newcommand{\ulEl}{\ul{\El\!}\,}
\newcommand{\uulEl}{\ulk{\smash{\ulEl}\rule[-0.2ex]{0mm}{0mm}}}
\newcommand{\ulk}[1]{\ul{#1\!}\,}
\newcommand{\ulF}{\ul{F\!\!}\,\,}
\newcommand{\ulE}{\ul{E\!}\,}

\newcommand{\uulE}{\ulk{\smash{\ulE}\rule[-0.2ex]{0mm}{0mm}}}
\newcommand{\uul}[2]{\ulk{\smash{\ulk{#2}}\rule[-0.#1ex]{0mm}{0mm}}}

\newcommand{\tru}{\vartriangle}
\newcommand{\trd}{\triangledown}
\renewcommand{\ind}{\!\!\shortuparrow}
\renewcommand{\sind}{\shortuparrow}
\newcommand{\dDe}{\dot\De}
\newdir{ >}{{}*!/-5pt/@{>}}
\DeclareMathOperator{\III}{{\sf I}}
\DeclareMathOperator{\PPP}{{\sf P}}
\DeclareMathOperator{\TTT}{{\sf T}}
\DeclareMathOperator{\per}{periodic}
\DeclareMathOperator{\qcyc}{qcyc}
\newcommand{\ac}{{\operatorname{ac}}}
\newcommand{\spac}{{\operatorname{sp\, ac}}}
\newcommand{\trud}{{\tru\!\trd}}

\begin{document}

\title{Heller triangulated categories}
\author{Matthias K\"unzer}
\maketitle

\begin{small}
\begin{quote}
\begin{center}{\bf Abstract}\end{center}\vspace*{2mm}

Let $\El$ be a Frobenius category. Let $\uulEl$ denote its stable category. The shift functor on $\uulEl$ induces, by pointwise application, an {\it inner} shift functor
on the category of acyclic complexes with entries in $\uulEl$. Shifting a complex by $3$ positions yields an {\it outer} shift functor on this category.
Passing to the quotient modulo split acyclic complexes, {\sc Heller} remarked that inner and outer shift become isomorphic, via an isomorphism satisfying still a further compatibility. Moreover, 
{\sc Heller} remarked that a {\sf choice} of such an isomorphism determines a Verdier triangulation on $\uulEl$, except for the octahedral axiom.
We generalise the notion of acyclic complexes such that the accordingly enlarged version of Heller's construction includes octahedra.
\end{quote}
\end{small}

\renewcommand{\thefootnote}{\fnsymbol{footnote}}
\renewcommand{\thefootnote}{\arabic{footnote}}

\begin{footnotesize}
\renewcommand{\baselinestretch}{0.7}%
\parskip0.0ex%
\tableofcontents%
\parskip1.2ex%
\renewcommand{\baselinestretch}{1.0}%
\end{footnotesize}%

\input{h3intro}   %
\input{h3def}     %
\input{h3equ}     %
\input{h3verdier} %
\input{h3frobhot} %
\input{h3qcyc}    %
\input{h3general} %
\input{h3ref}     %
\end{document}

%% file: h3intro.tex
\setcounter{section}{-1}

\section{Introduction}

\subsection[Heller's idea]{Heller's idea (\footnote{{\sc Heller} formulated his idea using Freyd categories. We will rephrase it using complexes, for this is the language we will use below. 
Cf.\ \S\S\,\ref{SecHellerTriangulations},\,\ref{SecEnlarge}.})}

\subsubsection{Stable Frobenius categories and an isomorphism between outer and inner shift}
\label{SecIntrHCFC}

Let $\El$ be a Frobenius category, i.e.\ an exact category with enough bijective objects. For instance, the category of complexes with values in an additive category, equipped with pointwise
split exact sequences, is a Frobenius category.

Let $\uulEl$ denote the stable category of $\El$; cf.\ \S \ref{SecResBeg}. Assume that $\uulEl$ has split idempotents.

A complex with entries in $\uulEl$ is {\it acyclic} if any Hom functor turns it into an acyclic complex of abelian groups.
Let $\uulEl^+(\b\De_2^\#)$ denote the category of acyclic complexes with entries in $\uulEl\;$ (\footnote{The notation using the diagram $\b\De_2^\#$ is chosen to fit into a larger framework; see 
\S\ref{SecHellerTriangulations} for more details.}). Let $\ulk{\uulEl^+(\b\De_2^\#)}\ru{-2.3}$ denote the homotopy category of the category $\uulEl^+(\b\De_2^\#)$ of acyclic complexes; that is, the 
quotient category of acyclic complexes modulo split acyclic complexes.

There is a shift automorphism $\TTT$ on $\uulEl$. It induces a first, {\it inner shift} automorphism 
$\TTT^+(\b\De_2^\#)$ on $\uulEl^+(\b\De_2^\#)$ by pointwise application. 

There is also a shift automorphism $\TTT_{\! 2}$ on the diagram $\b\De_2^\#$. It induces a second, {\it outer shift} automorphism 
$\uulEl^+(\TTT_{\! 2})$ on $\uulEl^+(\b\De_2^\#)$, shifting a complex by three positions. 

Both outer and inner shift induce automorphisms
\[
\phantom{.}\hspace*{34mm}\ulk{\uulEl^+(\TTT_{\! 2})} \;\;\;\text{resp.}\;\;\; \ulk{\TTT^+(\b\De_2^\#)} \hspace*{10mm}\text{on}\;\;\;\; \ulk{\uulEl^+(\b\De_2^\#)}\;\; . 
\]
{\sc Heller} remarked that these functors are isomorphic. But there is no {\sf a priori given} isomorphism. So he {\sf chose} an \mb{isomorphism}
\[
\ulk{\uulEl^+(\TTT_{\! 2})} \;\;\; \lraisoa{\tht_2} \;\;\; \ulk{\TTT^+(\b\De_2^\#)}\; ,
\]
satisfying, for technical reasons, still a further compatibility.

Then he remarked that the choice of such an isomorphism $\tht_2$ determines a triangulation on $\uulEl$ in the sense of {\sc Puppe}
\bfcite{Pu62}{Sec.\ 2}; that is, it satisfies all the axioms of {\sc Verdier} \mb{\bfcite{Ve63}{Def.\ 1-1}} except possibly for the octahedral axiom. 
Namely, as distinguished triangles we take acyclic complexes on which outer and inner shift coincide (i.e., which are ``$3$-periodic up to shift'') and on which $\tht_2$ is the identity. 

Whether this observation now fathoms Puppe triangulations remains to be discussed. Whenever two objects are isomorphic but lack a nature-given isomorphism, 
it is at any rate not unusual to pick an isomorphism. Once a suitable isomorphism between our shift functors chosen, a Puppe triangulation ensues. In nontechnical terms, we may let the relation 
between the two shifts govern the Puppe triangulations. This is a possible point of view, which we shall adopt and put into a larger framework; cf.\ \S\ref{SecHellerTriangulations}.

{\sc Heller} used this construction to parametrise Puppe triangulations on $\uulEl$. The non-uniqueness of such a Puppe triangulation on $\uulEl$, and hence the impossibility of an intrinsic 
definition of distinguished triangles, thus can be regarded as rooted in the possible nontriviality of the automorphism group of the inner shift 
functor $\ul{\TTT^+(\b\De_2^\#)}\ru{-2.3}$, or, by choice, of the outer shift functor $\ulk{\uulEl^+(\TTT_{\! 2})}\ru{-2.3}$. This is to be seen in contrast to the intrinsic characterisation of short 
exact sequences in an abelian category. 

\subsubsection{The stable Frobenius case models a general definition of Puppe triangulations}
\label{SecIntrWeakAbel}

A {\it weak kernel} in an additive category is defined by the universal property of a kernel, except for the uniqueness of the induced morphism; dually a {\it weak cokernel}.

A {\it weakly abelian category} is an additive category in which each morphism has a weak kernel and a weak cokernel, and in which each morphism is a weak kernel and a weak cokernel.
For instance, the stable category $\uulEl$ appearing in \S\ref{SecIntrHCFC} is a weakly abelian category.

Let $\Cl$ be a weakly abelian category with split idempotents carrying a shift automorphism $\TTT$. Now Heller's construction yields an alternative, equivalent definition of a Puppe 
triangulation on $(\Cl,\TTT)$ as being an isomorphism
\[
\ulk{\Cl^+(\TTT_{\! 2})}\;\; \lraisoa{\tht_2} \;\; \ulk{\TTT^+(\b\De_2^\#)} 
\]
satisfying still a further compatibility. In other words, a Puppe triangulated category can be defined to be such a triple $(\Cl,\TTT,\tht_2)$.

\subsubsection{From Puppe to Verdier and beyond}

In a Puppe triangulated category, Verdier's octahedral axiom \bfcite{Ve63}{Def.\ 1-1} does not seem to hold in general (\footnote{The author lacks an example of a category that
is Puppe but not Verdier triangulated, but strongly suspects that such an example exists, i.e.\ that the octahedral axiom is not a consequence of Puppe's axioms; cf.\ Question \ref{Question1}. 
In any case, such a deduction is unknown.}).

In a Verdier triangulated category, in turn, it seems to be impossible to derive the existence of the two extra triangles in a particular octahedron described in \bfcite{BBD84}{1.1.13}, or to 
distinguish crosses as in \bfcite{Il03}{App.}.

Moreover, to define a K-theory simplicial set of a triangulated category, one is inclined to take objects as $1$-simplices, distinguished triangles as $2$-simplices, {\sf distinguished} octahedra as 
$3$-simplices, {\sf etc.}

So we enlarge the framework, generalising from $\Cl^+(\b\De_2^\#)$ to $\Cl^+(\b\De_n^\#)$, as described next in \S\ref{SecIntroDef}. 

\subsection{Definition of Heller triangulated categories} 
\label{SecIntroDef}

\subsubsection{A diagram shape}

Given $n\ge 0$, we let $\De_n := \{ i\in\Z\; :\; 0\le i\le n\}$, considered as a linearly ordered set. Let $\b\De_n$ be the {\it periodic prolongation} of $\De_n$, consisting of $\Z$ copies 
of $\De_n$ put in a row. This is a {\it periodic linearly ordered set}; that is, a linearly ordered set equipped with a shift automorphism $i\lramaps i^{+1}$. For instance, 
$\b\De_2 = \{\dots,2^{-1},0,1,2,0^{+1},\dots\}$, equipped with $i\lramaps i^{+1}$. Let $\b\Deltab$ be the category consisting of 
periodic linearly ordered sets of the form $\b\De_n$ as objects, and of monotone shiftcompatible maps as morphisms.

Let $\b\De_n(\De_1)$ denote the category of morphisms in $\b\De_n$, i.e.\ the category of $\b\De_n$-valued diagrams of shape $\De_1$. 
Given $\al,\,\be\,\in\,\b\De_n$ such that $\al\le\be$, the object $(\al\lra\be)$ in $\b\De_n(\De_1)$ is abbreviated by $\be/\al$.

Let $\b\De_n^\#$ be the full subcategory of $\b\De_n(\De_1)$ that consists of objects $\be/\al$ within a single period, i.e.\ 
such that $\be^{-1}\le\al\le\be\le\al^{+1}$. For instance,
\[
\begin{picture}(300,0)
\put(0,-200){$\b\De_2^\# \Icm = \Icm$}
\end{picture}
\xymatrix{
                     &                      &                      & 0^{+1}/0^{+1}\ar[r]  & \cdots\\
                     &                      & 2/2\ar[r]            & 0^{+1}/2\ar[u]\ar[r] & \cdots\\
                     & 1/1\ar[r]            & 2/1\ar[u]\ar[r]      & 0^{+1}/1\ar[u]\ar[r] & \cdots\\
0/0\ar[r]            & 1/0\ar[u]\ar[r]      & 2/0\ar[u]\ar[r]      & 0^{+1}/0\ar[u]       &       \\
\vdots\ar[u]         & \vdots\ar[u]         & \vdots\ar[u]         &                      &       \\
}
\]

\subsubsection{Heller triangulations}
\label{SecHellerTriangulations}

Let $\Cl$ be a weakly abelian category; cf.\ \S\ref{SecIntrWeakAbel}. A sequence $X\lraa{f} Y\lraa{g} Z$ in $\Cl$ is called {\it exact at $Y$} if $f$ is a weak kernel of $g$, or, equivalently,
if $g$ is a weak cokernel of $f$. A commutative quadrangle in $\Cl$ whose diagonal sequence is exact at the middle object is called a {\it weak square}.

Let $\Cl^+(\b\De_n^\#)$ be the category of $\Cl$-valued diagrams of shape $\b\De_n^\#$ with a zero at
$\al/\al$ and at $\al^{+1}/\al$ for each $\al\in\b\De_n\,$, and such that the quadrangle on $(\ga/\al,\de/\al,\ga/\be,\de/\be)$
is a weak square whenever $\de^{-1}\le\al\le\be\le\ga\le\de\le\al^{+1}$. Let $\ulk{\Cl^+(\b\De_n^\#)}$ be the quotient of 
$\Cl^+(\b\De_n^\#)\ru{-2.3}$ modulo the full subcategory of diagrams therein that consist entirely of split morphisms. 

For instance, $\Cl^+(\b\De_2^\#)$ is the category of $\Cl$-valued acyclic complexes; and $\ulk{\Cl^+(\b\De_2^\#)}$ is its quotient modulo split acyclic complexes, i.e.\ the homotopy category of
$\Cl$-valued acyclic complexes.

Furthermore, suppose given an automorphism $\TTT$ on $\Cl$. We obtain two shift functors on $\ulk{\Cl^+(\b\De_n^\#)}\ru{-2.3}$, the {\it inner shift} given by pointwise
application of $\TTT$, and the {\it outer shift} induced by a diagram shift $j/i\lramaps i^{+1}/j$. 

A {\it Heller triangulation} on $(\Cl,\TTT)$ is a tuple of isomorphisms $\tht = (\tht_n)_{n\ge 0}$, where $\tht_n$ is an isomorphism from the outer to the inner shift on
$\ulk{\Cl^+(\b\De_n^\#)}$. This tuple is required to be compatible with the functors induced by periodic monotone maps between $\b\De_n$ and $\b\De_m$,
where $m,\, n\,\ge 0$. Moreover, it is required to be compatible with an operation called {\it folding,} which emerges from the fact that a weak square 
\[
\xymatrix{
X\ar[r]^f                         & Y \\
X'\ar[r]_{f'}\ar[u]^x\ar@{}[ur]|+ & Y'\ar[u]_y \\
}
\]
entails a {\it folded} weak square 
\[
\xymatrix{
0\ar[r]                                     & Y \\
X'\ar[r]_(0.4){\smatez{x}{f'}}\ar[u]\ar@{}[ur]|+ & X\dk Y'\ar[u]_{\rsmatze{f}{-y}} \; . \!\!\\
}
\]
A {\it Heller triangulated category} is a triple $(\Cl,\TTT,\tht)$ as just described, often just denoted by $\Cl$.

An {\it $n$-triangle} in a Heller triangulated category $\Cl$ is an object $X$ of $\Cl^+(\b\De_n^\#)$ that is {\it periodic} in the sense that outer shift and inner shift coincide on $X$, and 
that satisfies $X\tht_n = 1$. The usual properties of $2$-triangles generalise to $n$-triangles. 

If $\Cl$ is a Heller triangulated category in which idempotents split, then, taking the $2$-triangles as the distinguished triangles, it is also triangulated in the sense of {\sc Verdier} 
\bfcite{Ve63}{Def.\ 1-1}; see Proposition \ref{PropV5}.

\subsubsection{Strictly exact functors}

An additive functor $\Cl\lraa{F} \Cl'$ between Heller triangulated categories $(\Cl,\TTT,\tht)$ and $(\Cl',\TTT',\tht')$ is called {\it strictly exact} if, firstly, it respects weak kernels, or, 
equivalently, weak cokernels; if, secondly, $F\TTT' = \TTT F$; and if, thirdly, the functor 
\[
\ulk{\Cl^+(\b\De_n^\#)}\;\mrafl{30}{\ulk{F^+(\b\De_n^\#)}}\;\ulk{\Cl'^+(\b\De_n^\#)}\;\; ,
\]
induced by pointwise application of $F$, satisfies $\;\ulk{F^+(\b\De_n^\#)}\st\tht'_n \= \tht_n\st\ulk{F^+(\b\De_n^\#)}\;$ for $n\ge 0$.

\subsubsection{Enlarge to simplify}
\label{SecEnlarge}

Let $\dDe_n := \{ i\in\Z\; :\; 1\le i\le n\}$. We have an embedding $\dDe_n\hra\b\De_n^\#\ru{-2}$ via $\al\lramaps\al/0$. Let $\Cl$ be a weakly abelian category. Let $\Cl(\dDe_n)$ denote the
category of $\Cl$-valued diagrams of shape $\dDe_n$. Let $\ulk{\Cl(\dDe_n)}\ru{-2.3}$ be the quotient of $\Cl(\dDe_n)\ru{-2.3}$ modulo the full subcategory of split diagrams. Restriction induces an 
equivalence 
$$
\ulk{\Cl^+(\b\De_n^\#)}\;\;\;\lraisofl{35}{(-)|_{\dDe_n}}\;\;\; \ulk{\Cl(\dDe_n)}\; ,\ru{-2.3}
\leqno (\ast)
$$
which is also a useful technical tool; cf.\ Proposition \ref{PropR5}. 

At first sight, one might be inclined to prefer $\ulk{\Cl(\dDe_n)}\ru{-2.3}$ over $\ulk{\Cl^+(\b\De_n^\#)}$. It contains smaller diagrams and has a less elaborate definition.
By transport of structure along $(\ast)$, one obtains an outer shift on $\ulk{\Cl(\dDe_n)}\ru{-2.3}$ as well. By pointwise 
application of the shift functor on $\Cl$, we also obtain an inner shift on $\ulk{\Cl(\dDe_n)}\ru{-2.3}$. These could be compared in order to write down a definition of Heller triangulated
categories.

So why then did we prefer to use $\ulk{\Cl^+(\b\De_n^\#)}\ru{-2.3}$ in our definition of Heller triangulated categories in \S\ref{SecHellerTriangulations}?
Working with $\ulk{\Cl(\dDe_n)}\ru{-2.3}$, the indirect definition of the outer shift would cause problems. In practice, one would have to pass the equivalence $(\ast)$ back and forth.
The ``blown-up variant'' $\ulk{\Cl^+(\b\De_n^\#)}\ru{-2.3}$ of $\ulk{\Cl(\dDe_n)}\ru{-2.3}$ carries a directly defined outer shift functor and is thus easier to work with.

There is a further equivalence $\ulk{\Cl(\dDe_n)}\lraiso\ulk{\hat\Cl(\dDe_{n-1})}$, where $\hat\Cl$ denotes the Freyd category of $\Cl$, i.e.\ the universal abelian category containing $\Cl$,
and where $\ulk{\hat\Cl(\dDe_{n-1})}\ru{-2.3}$ is the quotient of $\hat\Cl(\dDe_{n-1})$ modulo split diagrams with entries in $\Cl$;
cf.\ Proposition \ref{LemIT1}. Originally, {\sc Heller} worked with $\ulk{\hat\Cl(\dDe_{n-1})}$ for $n = 2$, i.e.\ with $\hat\Cl/\Cl$.

\subsection{A result to begin with}
\label{SecResBeg}

Let $\El$ be a Frobenius category. We define its {\it stable category} $\uulEl$ to be the quotient category of the category of purely acyclic complexes with values
in the bijective objects of $\El$, modulo the subcategory of split acyclic such complexes. 

Then $\uulEl$ is equivalent to the {\it classical stable category} $\ulEl$ of $\El$, defined as the quotient category of $\El$ modulo bijective objects. But $\uulEl$ carries a shift automorphism 
$\TTT$ (invertible), whereas $\ulEl$ carries, in general, only a shift autoequivalence (invertible up to isomorphism). In this sense, $\uulEl$ is a ``strictified version'' of $\ulEl$.

{\bf Theorem} (Corollary \ref{Cor1Epi5}, Corollary \ref{CorExFun}). {\it Given a Frobenius category $\El$, there exists a Heller triangulation $\tht$ on $(\uulEl,\TTT)$. An exact functor 
$\El\lrafl{22}{E}\El'$ between Frobenius categories that sends all bijective objects of $\El$ to bijective objects of $\El'$ induces a strictly exact functor $\uulEl\lrafl{27}{\uulE}\uulEl'$.}

The Verdier triangulated version of this theorem is due to {\sc Happel} \bfcite{Ha88}{Th.\ 2.6}.

\subsection{A quasicyclic category}

Let $\Cl$ be a Heller triangulated category.

A {\it quasicyclic category} is a contravariant functor from $\b\Deltab^{\!\circ}$ to the ($1$-)category of categories. Letting 
$\qcyc_n\Cl$ be the subcategory of isomorphisms in $\Cl^+(\b\De_n^\#)\ru{-1.8}$ for $n\ge 0$, we obtain a quasicyclic category $\qcyc_\bt \Cl$.
There is a quasicyclic subcategory $\qcyc^{\tht=1}_\bt \Cl\ru{-2}$ that consists only of $n$-triangles and their isomorphisms instead of all objects in $\Cl^+(\b\De_n^\#)$ and their isomorphisms
(\footnote{Here, ``$\;\tht\!=\!1\;$'' is a mere symbol that should evoke the definition of $n$-triangles via $\tht$.}).

Restricting $\qcyc^{\tht=1}_\bt\Cl$ along the functor $\Deltab^{\!\circ}\hra\b\Deltab^{\!\circ}$ of ``periodic prolongation'', this yields a simplicial category, hence a topological space; 
depending functorially on $\Cl$. This space is the author's tentative proposal for the definition of the K-theory of $\Cl$; cf.\ \bfcite{Ne05}{Rem.\ 63}. Of course, this definition still needs to be 
justified by results one expects of such a K-theory, which has not yet been attempted.

\subsection{Some remarks}

A comparison of our theory to the derivator approach and related constructions in \bfcit{De68}, \mb{\bfcite{Il71}{chap.\ V.1}}, \bfcit{He88}, \bfcit{Gr90}, \bfcit{Ke91}, \bfcit{Fr96} and \bfcit{Malt02} 
would be interesting. One might ask whether the base category of a triangulated derivator in the sense of \bfcit{Malt02} carries a Heller triangulation; and if so, whether morphisms of triangulated
derivators give rise to strictly exact functors.

Our approach differs from the derivator approach in that we consider a {\sf single} category $\Cl$ with shift and an ``exactness structure'', i.e.\ a Heller triangulation, on it.
The categories $\Cl^+(\b\De_n^\#)$ needed to define this ``exactness structure'' on $\Cl$ consist of veritable $\Cl$-valued diagrams; cf.\ \S\ref{SecIntroDef}.
In particular, a ``structure preserving map'' between two such categories $\Cl$ and $\Cl'$, i.e.\ a strictly exact functor, is a {\sf single} additive functor 
$\Cl\lraa{F}\Cl'$ compatible with the ``exactness structures'' imposed on $\Cl$ and on $\Cl'$. In contrast, a ``structure preserving map'' of triangulated derivators is a compatible {\sf family} of 
additive functors.

The generalised triangles in \bfcite{BBD84}{1.1.14} are, in our language, $n$-pretriangles for which the \mb{$2$-pretriangle} obtained by restriction along any periodic monotone map 
$\b\De_2\lra\b\De_n$ is a \mb{$2$-triangle}. An $n$-triangle is such a generalised triangle, but the converse does not hold in general, as pointed out to me by {\sc A.\ Neeman.} For an example, 
see \bfcite{Ku07}{\S 2}.

It is conceivable that the concept of Heller triangulated categories is essentially equivalent to a direct axiomatisation via $n$-triangles, as worked out independently 
by {\sc G.\ Maltsiniotis}~\bfcit{Malt05} and myself~\bfcit{Ku96}. To compare these approaches, on the one hand, the Heller triangulated category should be closed in the sense of 
\bfcite{Ku07}{\S A.2, Def.\ 13}; on the other hand, the $n$-triangles should be stable under folding as in Lemma~\ref{LemV3}.(2).

Concerning the motivation to consider triangulated categories at all, and in particular derived categories, conceived by {\sc Grothendieck}, we refer the 
reader to the introduction of the thesis of {\sc Verdier} \bfcit{Ve67}; cf.\ also \bfcit{Il90} and \bfcite{We97}{p.\ 26}. 

\subsection{Acknowledgements}

I thank {\sc A.\ Wiedemann} for an introduction to derived categories. I thank {\sc A.\ Beligiannis} for directing me to
Heller's parametrisation of Puppe triangulations, and for helpful remarks. I thank {\sc B.\ Keller} for the hint how to ``strictify'' the classical stable category of 
a Frobenius category using acyclic complexes. I thank {\sc A.\ Neeman} for the hint leading to the folding operation, and for corrections. 

More than once I returned to {\sc A.\ Heller}'s original construction \bfcite{He68}{p.\ 53--54}, the reference not only for the basic idea, but also for arguments perfectly extendable to the 
more general framework used here.  

\subsection{Notations and conventions}

\begin{footnotesize}
\begin{itemize}
\item[(i)] The disjoint union of sets $X$ and $Y$ is written $X\disj Y$.

\item[(ii)] Given $a,\, b,\, c\,\in\,\Z$, the assertion $a\con_c b$ is defined to hold if there exists a $z\in\Z$ such that $a - b = cz$.

\item[(iii)] For $a,\,b\,\in\,\Z$, we denote by $[a,b] := \{z\in\Z\;:\; a\le z\le b\}$ the integral interval. Similarly, we let \linebreak
$[a,b[ \; := \{z\in\Z\;:\; a\le z < b\}$, $]a,b] := \{z\in\Z\;:\; a < z\le b\}$, $\Z_{\ge 0} := \{z\in\Z\;:\; z\ge 0\}$ and \linebreak
$\Z_{\le 0} := \{z\in\Z\;:\; z\le 0\}$.

\item[(iv)] All categories are supposed to be small with respect to a sufficiently big universe.

\item[(v)] Given a category $\Cl$, and objects $X$, $Y$ in $\Cl$, we denote the set of morphisms from $X$ to $Y$ by $\liu{\Cl}{(X,Y)}$, or simply
by $(X,Y)$, if unambiguous.

\item[(vi)] Given a poset $P$, we frequently consider it as a category, letting $\liu{P}(x,y)$ contain one element $y/x$ if $x\le y$, and
letting it be empty if $x\not\le y$, where $x,y\in\Ob P = P$. 

\item[(vii)] Given $n\ge 0$, we denote by $\De_n := [0,n]$ the linearly ordered set with ordering induced by the standard ordering on $\Z$.
Let $\dDe_n := \De_n\ohne\{ 0\} = [1,n]$, considered as a linearly ordered set.

\item[(viii)] Maps act on the right. Composition of maps, and of more general morphisms, is written on the right, i.e.\ $\lraa{a}\lraa{b} = \lraa{ab}$. 

\item[(ix)] Functors act on the right. Composition of functors is written on the right, i.e.\ $\lraa{F}\lraa{G} = \lraa{FG}$. Accordingly, the entry
of a transformation $a$ between functors at an object $X$ will be written $Xa$. 

\bq
The reason for this convention is that we will mainly consider functors of type ``restriction to a subdiagram'' or ``shift'', and such operations are usually written on the right. 
\eq

\item[(x)] A functor is called {\it strictly dense} if its map on the objects is surjective. It is called {\it dense} if its induced map on the isoclasses is surjective.

\item[(xi)] Given transformations 
$
\xymatrix{
\Cl \ar@/^/[r]^F="F"\ar@/_/[r]_G="G"  & \Cl' \ar@/^/[r]^{F'}="F'"\ar@/_/[r]_{G'}="G'" & \Cl'' \; ,\\
\ar@2"F"+<0mm,-2.9mm>;"G"+<0mm,2.4mm>^a
\ar@2@<-1mm>"F'"+<0mm,-2.9mm>;"G'"+<0mm,2.4mm>^(0.4){a'}
} 
$\vspace{-7mm}
we write $a\st a'$ for the transformation from $FF'$ to $GG'$ given at $X\in\Ob\Cl$ by $X(a\st a') := (XFa')(XaG') = (XaF')(XGa')$. In this context, we also write the object $F$ for the identity
$1_F$ on this object, i.e.\ e.g.\ $X(F\st a') = X(1_F\st a') = (XF)a'$.

\item[(xii)] The inverse of an isomorphism $f$ is denoted by $f^-$. Note that if we denote an iterated shift automorphism $f\lramaps f^{+1}$ by $f\lramaps f^{+z}$ for $z\in\Z$, then
we have to distinguish $f^-$ (inverse isomorphism if $f$ is an isomorphism) and $f^{-1}$ (inverse of the shift functor applied to $f$).

\item[(xiii)] In an exact category, pure monomorphy is indicated by $X\lramono Y$, pure epimorphy by $X\lraepi Y$.

\item[(xiv)] A morphism in an additive category $\Al$ is {\it split} if it is isomorphic, in $\Al(\De_1)$, to a morphism of the form $X\ds Y\lrafl{30}{\smatzz{0}{0}{1}{0}} Y\ds Z\ru{6}\,$.
A morphism being split is indicated by $\xymatrix@C=6mm{X\;\;\ar~+{|<*\dir{>}}[r] & Y}$ (not to be confused with monomorphy). Accordingly, a morphism being a split monomorphism is
indicated by $\xymatrix@C=6mm{X\;\;\ar~+{|<*\dir{>}|(0.46)*\dir{*}}[r] & Y}$, a morphism being a split epimorphism by $\xymatrix@C=6mm{X\;\;\ar~+{|<*\dir{>}|(0.46)*\dir{|}}[r] & Y}$. 
Cf.\ \S\ref{SecDefEx}.

\item[(xv)] We say that {\it idempotents split} in an additive category $\Al$ if every endomorphism $e$ in $\Al$ that satisfies $e^2 = e$ is split. 

\item[(xvi)] The category of functors and transformations from a category $D$ to a category $\Cl$ is denoted by $\fbo D,\Cl\fbc$ or by $\Cl(D)$. To objects in $\Cl(D)$, we also refer to as 
{\it diagrams on $D$} with {\it values} or {\it entries} in $\Cl$.

\item[(xvii)] If $\Cl$ and $D$ are categories, and $X\in\Ob\Cl(D)$, we usually write $(d\lraa{a} e)X =: (X_d\lraa{X_a} X_e)$ for a morphism $d\lraa{a} e$ in $D$.
If the morphism $a$ is unambiguously given by the context, we also use small letters to
write $(X_d\lraa{x} X_e) := (X_d\lraa{X_a} X_e)$ (similarly $X'_d\lraa{x'} X'_e\,$, $\w Y_d\lraa{\w y} \w Y_e\,$, \dots)

\item[(xviii)] Let $\Add$ denote the $2$-category of additive categories.

\item[(xix)] Given an additive category $\Al$ and a full additive subcategory $\Bl\tm\Al$, we denote by $\Al/\Bl$ the quotient of $\Al$ by $\Bl$, having 
as objects the objects of $\Al$ and as morphisms equivalence classes of morphisms of $\Al$; where two morphisms $f$ and $f'$ are equivalent,
written $f\con_\Bl f'$, if their difference factors over an object of $\Bl$.

\item[(xx)] In an exact category, an object $P$ is called {\it projective} if $(P,-)$ turns pure epimorphisms into epimorphisms. An object $I$ is called 
{\it injective} if $(-,I)$ turns pure monomorphisms into epimorphisms. It is called {\it bijective} if it is injective and projective. See 
\S\ref{SecExCat} for details.

\item[(xxi)] In an additive category, a morphism $K\lraa{i} X$ is called a {\it weak kernel} of a morphism $X\lraa{f} Y$ if for every morphism $T\lraa{t} X$ with 
$tf = 0$ there exists a morphism $T\lraa{t'} K$ with $t'i = t$. A {\it weak cokernel} is defined dually. An additive category is called {\it weakly abelian} if every morphism has a weak kernel
and a weak cokernel, and is a weak kernel and a weak cokernel.

\item[(xxii)] The Freyd category of a weakly abelian category $\Cl$ is written $\h\Cl$. See \S\ref{SecFreyd} for details.

\item[(xxiii)] In an abelian category, a commutative quadrangle 
\[
\xymatrix{
X\ar[r]^f\ar[d]_x     & Y\ar[d]^y \\
X'\ar[r]_{f'}         & Y' \\
}
\]
is called a {\it square} if its {\it diagonal sequence} $X\lrafl{25}{\smatez{x}{f}} X'\ds Y\lrafl{40}{\rsmatze{f'}{-y}} Y'$ is short exact. Being a square is indicated by a box sign 
``$\,\Box\,$'' in the quadrangle.

The quadrangle $(X,Y,X',Y')$ is called a {\it weak square} if its diagonal sequence is exact in the middle; cf.\ Definition \ref{DefSquare}. Being a weak square is indicated by a ``$\,+\,$''-sign in the
quadrangle.

In an exact category, $(X,Y,X',Y')$ is a {\it pure square} if it has a pure short exact diagonal sequence. Being a pure square is indicated by a box sign 
``$\,\Box\,$'' in the quadrangle.

In a weakly abelian category, $(X,Y,X',Y')$ is a {\it weak square} if it is a weak square in the Freyd category of that weakly abelian category.

\item[(xxiv)] In an abelian category, given a morphism $X\lraa{f} Y$, we sometimes denote its kernel by $\KK_f$, and its cokernel by $\CC_f$.
\end{itemize}
\end{footnotesize}

%% file: h3def.tex
\newpage

\section{Definition of a Heller triangulated category} 
\label{SecDef}

\subsection{Periodic linearly ordered sets and their strips}

Without further comment, we consider a poset $D$ as a category, whose set of objects is given by $D$, and for which $\#\liu{D}{(\al,\be)} = 1$ 
if $\al\le\be$, and $\#\liu{D}{(\al,\be)} = 0$ otherwise. If existent, i.e.\ if $\al\le\be$, the morphism from $\al$ to $\be$ is denoted by $\be/\al$.
A full subposet of a category is a full subcategory that is a poset. In particular, a full subposet of a poset is just a full subcategory of that poset.

A {\it periodic poset} is a poset $P$ together with an automorphism $\TTT:P\lraiso P$, \mb{$\al\lramaps\al \TTT =: \al^{+1}$.}
Likewise, we denote $\al \TTT^m =: \al^{+m}$ resp.\ $\al \TTT^{-m} =: \al^{-m}$ for $m\in\Z_{\ge 0}$. By abuse of notation, we denote a periodic 
poset $(P,\TTT)$ simply by $P$.

A morphism of periodic posets $P\llaa{p} P'$ is a monotone map $p$ of the underlying posets such that $(\al'^{+1})p = ((\al')p)^{+1}$ for all 
$\al'\in P'$. The category of periodic posets shall be denoted by $\pp$.

A {\it periodic linearly ordered set} is a periodic poset the underlying poset of which is linearly ordered, i.e.\ such that $\#\big(\liu{D}{(\al,\be)} \cup \liu{D}{(\be,\al)}\big) = 1$ for all 
$\al,\,\be\,\in\, D$.

To any linearly ordered set $D$ we can attach a periodic linearly ordered set $\b D$ by letting $\b D := D\ti\Z$, and $(\al,z)\le (\be,w)$ if $z\le w$, 
or if ($z = w$ and $\al\le\be$ in $D$). We let $(\al,z)^{+1} := (\al,z+1)$. Sending $D\lra \b D$, $\al\lramaps (\al,0)$, and identifying $D$ with its image, we obtain
$(\al,z) = \al^{+z}$, and the latter is the notation we will usually use. The periodic linearly ordered set $\b D$ is called the {\it periodic 
repetition} of $D$. Likewise, the functor $D\lramaps\b D$ from the category of linearly ordered sets to the category of periodic linearly ordered sets is
called {\it periodic repetition.}

Let $\Deltab$ be the full subcategory of the category of linearly ordered sets defined by \mb{$\Ob\Deltab := \{ \De_n \;:\; n\in\Z_{\ge 0}\}$.}

Let $\b\Deltab$ be the full subcategory of the category of periodic linearly ordered sets defined by 
$\Ob\b\Deltab := \{ \b\De_n \;:\; n\in\Z_{\ge 0}\}$ (\footnote{The category $\b\Deltab$ is isomorphic to the category $L$ defined by {\sc Elmendorf} in \bfcit{El93}.}).

\bq
The reason for considering periodic linearly ordered sets is that the functor periodic repetition from $\Deltab$ to $\b\Deltab$ is dense and faithful but not full.
We will require a naturality of a certain construction with respect to $P\in\Ob\b\Deltab$, which is stronger than setting $P = \b D$ and requiring naturality 
with respect to $D\in\Ob\Deltab$.
\eq

Given $n\ge 0$, the underlying linearly ordered set of $\b\De_n$ is isomorphic to $\Z$ via 
\mb{$\al^{+z}\lramaps \al + (n+1)z$.} We use this isomorphism to define the operation
\[
\b\De_n\ti\Z\;\lra\;\b\De_n\; ,\;\; (\al^{+z},x)\;\lramaps\; \al^{+z} + x \; :=\; (\ol{\al + x})^{+(z+\ul{\al + x})} \; , 
\]
where we write $k = (n+1)\ul{k} + \ol{k}$ with $\ul{k}\in\Z$ and $\ol{k}\in [0,n]$ for $k\in\Z$. 
For instance, if $n = 3$, then $2^{+1} + 7 = 1^{+3}$. 

To a periodic linearly ordered set $P$, we attach the poset 
\[
P^\# := \{\be/\al\in P(\De_1)\; :\; \be^{-1}\le\al\le\be\le\al^{+1}\}
\]
as a full subposet of $P(\De_1)$, called the {\it strip} of $P$. A morphism therein from $\be/\al$ to $\de/\ga$ is written $\de/\ga\bby\be/\al$, 
which is unique if it exists, i.e.\ if $\al\le\ga$ and $\be\le\de$.

The strip $P^\#$ carries the automorphism $\be/\al \lramaps (\be/\al)^{+1} := \al^{+1}/\be$ in $\pp$, where \mb{$\be/\al\in P^\#$.}

If $P = \b\De_n$, we also write $\be/\al \lramapsa{\TTT_{\! n}} (\be/\al)^{+1}$.

This construction defines a functor
\[
\barcl
\b\Deltab   & \lraa{(-)^\#} & \pp \\
P           & \lramaps      & P^\# \\
\ea
\]
which sends a morphism $P\llaa{p} P'$ in $\b\Deltab$ to 
\[
\barcl
P^\#        & \llaa{p^\#} & P'^\# \\
\be'p/\al'p & \llamaps & \be'/\al' \\
\ea
\]
In fact, $p^\#$ is welldefined, since if $\be'^{-1}\le\al'\le\be'\le\al'^{+1}$, then $(\be'p)^{-1}\le\al'p\le\be'p\le (\al'p)^{+1}$. Moreover, $p^\#$ is monotone and compatible with shift. 

\bq
\begin{Example}
\label{ExStrip}\rm
The periodic poset $\b\De_2^\#$, i.e.\ the strip of the periodic repetition of $\De_2$, can be displayed as
\[
\xymatrix@R=6mm@C=6mm{
                    &                      &                      &                      &                      & 1^{+1}/1^{+1}\ar[r]       & \cdots  \\
                    &                      &                      &                      & 0^{+1}/0^{+1}\ar[r]  & 1^{+1}/0^{+1}\ar[u]\ar[r] & \cdots  \\
                    &                      &                      & 2/2\ar[r]            & 0^{+1}/2\ar[u]\ar[r] & 1^{+1}/2\ar[u]\ar[r]      & \cdots  \\
                    &                      & 1/1\ar[r]            & 2/1\ar[u]\ar[r]      & 0^{+1}/1\ar[u]\ar[r] & 1^{+1}/1\ar[u]            & \\
                    & 0/0\ar[r]            & 1/0\ar[u]\ar[r]      & 2/0\ar[u]\ar[r]      & 0^{+1}/0\ar[u]       &                           & \\
2^{-1}/2^{-1}\ar[r] & 0/2^{-1}\ar[u]\ar[r] & 1/2^{-1}\ar[u]\ar[r] & 2/2^{-1}\ar[u]       &                      &                           & \\
\vdots\ar[u]        & \vdots\ar[u]         & \vdots\ar[u]         &                      &                      &                           & \\
}
\]
\end{Example}
\eq

\subsection{Heller triangulated categories}

Suppose given a weakly abelian category $\Cl$; cf.\ Definition \ref{DefFreydic}. From \S\ref{SecNatInC} on, we assume it to be equipped with an automorphism
\[
\barcl
\Cl           & \lraisoa{\TTT} & \Cl \\
(X\lraa{u} Y) & \lramaps       & (X\TTT\lraa{u\TTT} Y\TTT) \; =:\; (X^{+1}\lraa{u^{+1}} Y^{+1})\; .
\ea
\]
Similarly, we denote $(X\TTT^m\lraa{u\TTT^m} Y\TTT^m) =: (X^{+m}\lraa{u^{+m}} Y^{+m})$ for $m\in\Z$. 

Recall that its Freyd category $\h\Cl$ is an abelian Frobenius category, and that the image of $\Cl$ in $\h\Cl$, identified with $\Cl$, is a sufficiently
big subcategory of bijectives; cf.\ \S\ref{SecFreyd}.

\subsubsection{The stable category of pretriangles\ \ \ul{$\Cl^+(P^\#)\!$}}

\paragraph{Definition of\ \ \ul{$\Cl^+(P^\#)\!$}}
\label{SecDefUl}
\indent

Concerning the Freyd category $\h\Cl$ of $\Cl$, cf.\ \S\ref{SecFreyd}.
Concerning the notion of a weak square in $\h\Cl\ru{4.5}$, see Definition \ref{DefSquare}. A weak square in $\Cl$ is a weak square in $\h\Cl$ that has all four objects in $\Ob\Cl$.
Applying Remark \ref{RemBal}, we obtain an elementary way to characterise weak squares as having a diagonal sequence with first morphism being a weak kernel of the second; or, equivalently,
with second morphism being a weak cokernel of the first.

Given a periodic linearly ordered set $P$, we let $\Cl^+(P^\#)$ be the full subcategory of $\Cl(P^\#)$ defined by
\[
\Ob\,\Cl^+(P^\#) \; :=\; \left\{ \;\; X\in \Ob\,\Cl(P^\#)\;\;\; :\; 
\mb{
\begin{tabular}{rl}
1) & $X_{\al/\al} = 0$ and $X_{\al^{+1}/\al} = 0$ for all $\al\in P$. \\
2) & For all $\de^{-1}\le\al\le\be\le\ga\le\de\le\al^{+1}$ in $P$, \\
   & the quadrangle \\
   & $\xymatrix{
      X_{\ga/\be}\ar[r]^x                       & X_{\de/\be} \\
      X_{\ga/\al}\ar[r]^x\ar[u]^x\ar@{}[ur]|{+} & X_{\de/\al}\ar[u]_x \\      
      }$ \\
   & is a weak square (as indicated by $+$). \\
\end{tabular}
}\right\} \;\; .
\]
Note that we do {\sf not} require that $(X_{\al^{+1}/\ga}\lraa{x} X_{\be^{+1}/\de}) = (X_{\ga/\al}\lraa{x} X_{\de/\be})^{+1}$ for  
\mb{$\ga/\al,\,\de/\be\,\in\, P^\#$} with $\ga/\al\le\de/\be$. 

An object of $\Cl^+(P^\#)$ is called a {\it $P$-pretriangle.} Given $n\ge 0$, an object of $\Cl^+(\b\De_n^\#)$, i.e.\ a $\b\De_n$-pretriangle, is also called an {\it $n$-pretriangle.} 

Roughly put, an $n$-pretriangle is a diagram on the strip $\b\De_n^\#$ of the periodic repetition $\b\De_n$ of $\De_n$ consisting of weak squares with zeroes on the boundaries.

\bq
 \begin{Example}
 \label{Ex012pre}\rm
 A $0$-pretriangle consists of zero objects. A $1$-pretriangle is just a sequence $\dots,X_{0/1^{-1}},X_{1/0},X_{0^{+1}/1},\dots$ of objects of $\Cl$, decorated with some zero objects. 
 A $2$-pretriang\-le is a complex in $\Cl$ which becomes acyclic in $\h\Cl$ -- for short, which {\it \,is\,} acyclic --, decorated with some zero objects.
 \end{Example}
\eq

A morphism in $\Cl$ is split in $\h\Cl$ if and only if it factors in $\h\Cl$ into a retraction followed by a coretraction. Equivalently, its image,
taken in $\h\Cl$, is bijective as an object of $\h\Cl$.

Let $\Cl^{+,\,\spl}(P^\#)$ be the full subcategory of $\Cl^+(P^\#)$ defined by 
\[
\Ob\,\Cl^{+,\,\spl}(P^\#) \; := \; \left\{ X\in\Ob\,\Cl^+(P^\#) \; :\; 
\mb{
\begin{tabular}{l} $X_{\ga/\al}\lraa{x} X_{\de/\be}$ is split in $\h\Cl$ \\
for all $\ga/\al,\,\de/\be\,\in\, P^\#$ with $\ga/\al\le\de/\be$\\
\end{tabular}
}\right\}
\]
We denote the quotient category by 
\[
\ulk{\Cl^+(P^\#)}\; :=\; \Cl^+(P^\#)/\Cl^{+,\,\spl}(P^\#)\; ,
\]
called the {\it stable category of $P$-pretriangles.}

\bq
\begin{Example}
\label{ExTriv}\rm
We have $\ulk{\Cl^+(\b\De_0^\#)} = \ulk{\Cl^+(\b\De_1^\#)} = 0$\vspace{1mm}. The category $\ulk{\Cl^+(\b\De_2^\#)}$ can be regarded as the homotopy category of the category of acyclic complexes 
with entries in $\Cl$.
\end{Example}
\eq

\paragraph{Naturality of\ \ \ul{$\Cl^+(P^\#)\!$}\ \ in\ \ $P$}
\label{SecNatInP}
\indent

Suppose given periodic linearly ordered sets $P$, $P'$, and a morphism $P^\#\llaa{q} P'^\#$ of periodic posets such that either ($P = P'$ and $q = \TTT$, the shift functor on $P^\#$) or 
$q = p^\#$ for some morphism $P\llaa{p} P'$ of periodic linearly ordered sets. 

Recall that if $P = \b\De_n$, then we write alternatively $\TTT_{\! n}$ for the shift functor $\TTT$ on $\b\De_n^\#$.

We obtain an induced functor
\[
\barcl
\Cl^+(P^\#) & \lrafl{25}{\Cl^+(q)} & \Cl^+(P'^\#) \\
X           & \lramapsfl{0}{}      & X\big(\Cl^+(q)\big) \; :=\; q\, X \; ,\\
\ea
\]
given by composition of $q$, followed by $X$. 

In particular, the shift $\TTT$ on $P^\#$ induces a functor
\[
\barcl
\Cl^+(P^\#) & \lrafl{25}{\Cl^+(\TTT)} & \Cl^+(P^\#) \\
X           & \lramapsfl{0}{}         & [X]^{+1} \; :=\; X\big(\Cl^+(\TTT)\big)\; , \\
\ea
\]
called the {\it outer shift.} Note that if $P = \b\De_n$, then $[X]^{+1}_{\be/\al} = X_{(\be/\al)^{+1}} = X_{\al^{+1}/\be}$ for $\be/\al\in\b\De_n^\#$. On the stable category, this 
functor induces a functor
\[
\barcl
\ulk{\Cl^+(P^\#)} & \lrafl{30}{\ulk{\Cl^+(\TTT)}} & \ulk{\Cl^+(P^\#)} \\
X                 & \lramapsfl{0}{}               & [X]^{+1}, \;  \\
\ea
\]
likewise called the {\it outer shift.} 

Given a morphism $P\llaa{p} P'$ in $\b\Deltab$, we obtain an induced morphism $P^\#\llaa{p^\#} P'^\#$, and hence an induced functor usually
abbreviated by
\[
\barcl
\Cl^+(P^\#) & \mrafl{30}{p^\# \; :=\; \Cl^+(p^\#)} & \Cl^+(P'^\#) \\
X           & \mramaps\;                           & X p^\# := X\big(\Cl^+(p^\#)\big) \; .\\
\ea
\]
Likewise on the stable categories; we abbreviate $\ul{p}^\# \; :=\; \ulk{\Cl^+(p^\#)}$.

\bq
So altogether, we have defined $X p^\# := p^\# X$ ($X p^\#$: operation induced by $p$, applied to $X\; $; $p^\# X$: composition of $p^\#$ and $X$), 
which is a bit unfortunate, but convenient in practice.
\eq

Given a morphism $P\llaa{p} P'$ in $\b\Deltab$ and $X\in\Ob\,\Cl^+(P^\#)$, we have
\[
[X]^{+1} p^\# \= [X p^\#]^{+1}\; ,
\]
natural in $X$. Likewise on the stable categories.

Given $P,\, P'\,\in\,\Ob\b\Deltab$, a functor $\Cl^+(P^\#)\lraa{F}\Cl^+(P'^\#)$ is called {\it strictly periodic} if 
\[
[XF]^{+1} \= [X]^{+1}F\; ,
\]
natural in $X$. Likewise on the stable categories. 

\paragraph{Naturality of\ \ \ul{$\Cl^+(P^\#)\!$}\ \ in\ \ $\Cl$}
\label{SecNatInC}
\indent

An additive functor $\Cl\lraa{F}\Cl'$ is called {\it subexact} if the induced additive functor $\h\Cl\lraa{\h F}\h\Cl'$ is an exact functor of abelian 
categories; cf.\ \S\ref{SecFreyd}. Alternatively, it is subexact if and only if it preserves weak kernels, or, equivalently, weak cokernels; cf.\ Remark \ref{RemBal}.

Suppose given a subexact functor $\Cl\lraa{F}\Cl'$ and $P\in\Ob\Deltab$. We obtain an induced functor
\[
\barcl
\Cl^+(P^\#) & \mra{F^+(P^\#)} & \Cl'^+(P^\#) \\
X           & \mramaps        & X\,F^+(P^\#)\; , \\ 
\ea
\]
where, writing $Y := X\,F^+(P^\#)$, we let
\[
(Y_{\be/\al} \lraa{y} Y_{\de/\ga}) \; := \; (X_{\be/\al}F \lraa{xF} X_{\de/\ga}F)  
\]
for $\be/\al,\,\de/\ga\,\in\, P^\#$ with $\be/\al\le\de/\ga$.

In particular, the automorphism $\Cl\lraa{\TTT}\Cl$ induces an automorphism
\[
\barcl
\Cl^+(P^\#) & \mrafl{25}{\TTT^+(P^\#)} & \Cl^+(P^\#) \\
X           & \mramaps                 & [X^{+1}] \; := \; X\big(\TTT^+(P^\#)\big)\; , \\ 
\ea
\]
called the {\it inner shift.} Note that if $P = \b\De_n$, then $[X^{+1}]_{\be/\al} = X_{\be/\al}^{+1}$ for $\be/\al\in\b\De_n^\#$.

On the stable category, this induces an automorphism
\[
\barcl
\ulk{\Cl^+(P^\#)} & \mrafl{30}{\ulk{\TTT^+(P^\#)}} & \ulk{\Cl^+(P^\#)} \\
X                 & \mramaps                       & [X^{+1}] \; , \\ 
\ea
\]
likewise called the {\it inner shift.} 

\subsubsection{Folding}
\label{SecFold}

\bq
The following construction arose from a hint of {\sc A.\ Neeman,} who showed me a multitude of \mb{$2$-triangles} in an $n$-triangle similar to the 
two $2$-triangles explained in \bfcite{BBD84}{1.1.13}; cf.\ Definition~\ref{Def4}.(ii) below.
\eq

\paragraph{Some notation}
\indent

Given $P = (P,\TTT)\in\Ob\b\Deltab$, we denote by $2P$ the periodic poset $(P,\TTT^2)$. 

Given a linearly ordered set $D$, we let $\rh\disj D$ be the linearly ordered set having as underlying set $\{\rh\}\disj D$; and as partial order 
$\rh\le_{\rh\disj D}\al$ for all $\al\in D$, and $\al\le_{\rh\disj D}\be$ if $\al,\,\be\,\in\, D$ and $\al\le_D\be$. 

\bq
Roughly put, $2P$ is $P$ with doubled period, and $\rh\disj D$ is $D$ with an added initial object $\rh$. 
\eq

Let $n\ge 0$. We have an isomorphism of periodic linearly ordered sets
\[
\barcl
2\b\De_n & \lraiso  & \b\De_{2n+1} \\
k^{+l}   & \lramaps & 
\left\{
\ba{ll}
k^{+l/2}           & \mb{for $l\con_2 0$} \\
(k+n+1)^{+(l-1)/2} & \mb{for $l\con_2 1$} \\
\ea
\right. \\
\ea
\]
and an isomorphism of linearly ordered sets
\[
\barcl
\rh\disj\De_n & \lraiso  & \De_{n+1} \\
k             & \lramaps & k+1 \hspace*{5mm}\mb{for $k\in [0,n]$} \\
\rh           & \lramaps & 0\; . \\
\ea
\]
In order to remain inside $\b\Deltab$ resp.\ inside $\Deltab$, we use these isomorphisms as identifications. 

Then $P\lramaps 2P$ is natural in $P$ and therefore defines an endofunctor of $\b\Deltab$, and
$D\lramaps\rh\disj D$ is natural in $D$ and therefore defines an endofunctor of $\Deltab$.

\bq
 Given a linearly ordered set $D$, we will need to consider the periodic posets $2\b D$ and $\ol{\rh\disj D}$, formed using periodic repetition.
\eq

\paragraph{The folding operation}
\indent

Let $n\ge 0$. Let the strictly periodic functor
\[
\barcl
\Cl^+((2\b\De_n)^\#) & \lraa{\ffk_n} & \Cl^+(\ol{\rh\disj\De_n}^{\,\#}) \\
X                    & \lramaps      & X\ffk_n \\
\ea
\]
be determined on objects by the following data. Writing $Y = X\ffk_n$, we let
\[
\ba{lcl}
(Y_{\al/\rh} \lraa{y} Y_{\be/\rh}) 
& := & \left(X_{\al^{+1}/\al} \;\lraa{x}\; X_{\be^{+1}/\be}\right) \vspace*{4mm}\\
(Y_{\be/\rh} \lraa{y} Y_{\be/\al}) 
& := & \left(X_{\be^{+1}/\be} \;\lraa{\smatez{x}{x}}\; X_{\be^{+1}/\al^{+1}}\ds X_{\al^{+2}/\be}\right) \vspace*{4mm}\\ 
(Y_{\be/\al} \lraa{y} Y_{\de/\ga}) 
& := &\left(X_{\be^{+1}/\al^{+1}}\ds X_{\al^{+2}/\be}\;\lrafl{40}{\smatzz{x}{0}{0}{x}}\; X_{\de^{+1}/\ga^{+1}}\ds X_{\ga^{+2}/\de}\right)\vspace*{4mm}\\ 
(Y_{\de/\ga} \lraa{y} Y_{\rh^{+1}/\ga}) 
& := & \left(X_{\de^{+1}/\ga^{+1}}\ds X_{\ga^{+2}/\de} \;\lrafl{40}{\rsmatze{x}{-x}}\; X_{\ga^{+2}/\ga^{+1}}\right) \\
\ea
\]
for $\al,\,\be,\,\ga,\,\de\,\in\,\De_n$ with $\al\le\be$, with $\ga\le\de$ and with $\be/\al \le \de/\ga$. 
The remaining morphisms are given by composition. 

\bq
Note that $X\in\Ob\,\Cl^+((2\b\De_n)^\#)$, so e.g.\ $X_{\be^{+1}/\be} \neq 0$ is possible, whereas $X_{\be^{+2}/\be} = 0$ for $\be\in\De_n$.
\eq

We claim that $X\ffk_n$ is an object of $\Cl^+(\ol{\rh\disj\De_n}^{\,\#})$. 

In fact, by Lemma \ref{Lem1}, applied in the abelian category $\h\Cl$, we are reduced to \mb{considering} the quadrangles of $Y$ on 
$(\ga/\rh,\,\de/\rh,\,\ga/\be,\,\de/\be)$ for $\be,\,\ga,\,\de\,\in\,\De_n$
with $\be\le\ga\le\de$; on $(\ga/\al,\,\de/\al,\,\ga/\be,\,\de/\be)$ for $\al,\,\be,\,\ga,\,\de\,\in\,\De_n$
with $\al\le\be\le\ga\le\de$; and on 
\mb{$(\ga/\al,\,\rh^{+1}/\al,\,\ga/\be,\,\rh^{+1}/\be)$} for $\al,\,\be,\,\ga\,\in\,\De_n$ with $\al\le\be\le\ga$.

The quadrangle of $Y$ on $(\ga/\al,\,\de/\al,\,\ga/\be,\,\de/\be)$ is a weak square as the direct sum of two weak squares.

For the remaining quadrangles to be treated, Lemma \ref{Lem2} reduces us to \mb{considering} the 
quadrangles of $Y$ on $(\al/\rh,\,\be/\rh,\,\al/\al,\,\be/\al)$, on $(\be/\rh,\,\rh^{+1}/\rh,\,\be/\al,\,\rh^{+1}/\al)$ and on 
\mb{$(\be/\al,\,\rh^{+1}/\al,\,\be/\be,\,\rh^{+1}/\be)$} for $\al,\,\be\,\in\,\De_n$ with $\al\le\be$. These are in fact weak squares, as 
ensues from Lemma \ref{Lem3} and its dual assertion. This proves our claim.

This construction of $Y = X\ffk_n$ is functorial in $X$.

To prove that the folding operation passes to the stable categories, we have to show that for an object $X$ of $\Cl^{+,\,\spl}((2\b\De_n)^\#)$, the
folded object $X\ffk_n$ is in $\Cl^{+,\,\spl}(\ol{\rh\disj\De_n}^\#)$. Denote $Y := X\ffk_n$. Since $Y_{\al/\rh}\lraa{y} Y_{\be/\rh}$ is
split in $\h\Cl$ for all $\al,\,\be\,\in\,\De_n$ with $\al \le \be$, it suffices to prove the following lemma.

\begin{Lemma}
\label{LemR4}
Suppose given $m\ge 0$ and $Z\in\Ob\,\Cl^+(\b\De_m^\#)$ such that $Z_{\al/0}\lraa{z} Z_{\be/0}$ is split in $\h\Cl$ for all $\al,\,\be\,\in\,\De_m$ with $0 < \al \le \be$.
Then $Z\in\Ob\,\Cl^{+,\,\spl}(\b\De_m^\#)$.
\end{Lemma}

{\it Proof.} Consider the morphism $Z_{\ga/\al}\lraa{z} Z_{\de/\be}$ for $\ga/\al\le\de/\be$ in $\b\De_m^\#$. We have to show that it is split in $\h\Cl$, 
i.e.\ that its image, taken in $\h\Cl$, is bijective there. 
Unless $\al\le\be\le\ga\le\de\le\al^{+1}$, this morphism is zero, hence split in $\h\Cl$. If this condition holds, it is the diagonal morphism of the 
weak square $(Z_{\ga/\al},\; Z_{\de/\al},\; Z_{\ga/\be},\; Z_{\de/\be})$.

So by Lemma \ref{Lem3_1}, applied in the abelian category $\h\Cl$, we see that it suffices to show that the (horizontal)
morphism $Z_{\be/\al}\lraa{z} Z_{\ga/\al}$ is split in $\h\Cl$ and that the (vertical) morphism $Z_{\ga/\al}\lraa{z} Z_{\ga/\be}$ is split in $\h\Cl$ for all $\al,\,\be,\,\ga$ 
in $\b\De_m$ with $\ga^{-1}\le\al\le\be\le\ga\le\al^{+1}$.

The long exact sequence
\[
\hspace*{-4mm}
\cdots\;\lra\; Z_{\al/\be^{-1}}\;\lra\; Z_{\al/\ga^{-1}}\;\lra\; Z_{\be/\ga^{-1}}\;\lra\; Z_{\be/\al}
\;\lra\; Z_{\ga/\al}\;\lra\; Z_{\ga/\be}\;\lra\; Z_{\al^{+1}/\be}\;\lra\;\cdots
\]
in $\h\Cl$ shows that it suffices to show that the morphism $Z_{\be/\al}\lraa{z} Z_{\ga/\al}$ is split in $\h\Cl$ for all 
\mb{$0\le\al\le\be\le\ga < 0^{+1}$.}
In fact, first of all we may assume that $0 \le \al < 0^{+1}$, so that $0\le\al\le\be\le\ga\le\al^{+1} < 0^{+2}$.
Hence either $0\le\al\le\be\le\ga < 0^{+1}$, or $0\le\ga^{-1}\le\al\le\be< 0^{+1}$, or
$0\le\be^{-1}\le\ga^{-1}\le\al < 0^{+1}$.

Now we may assume that $0 < \al$ and apply Lemma \ref{Lem3_1} to the weak square
\mb{$(Z_{\be/0},\; Z_{\ga/0},\; Z_{\be/\al},\; Z_{\ga/\al})$,} in which $Z_{\be/0}\lraa{z} Z_{\ga/0}$ is split in $\h\Cl$ by assumption, in which 
\mb{$Z_{\be/0}\lraa{z} Z_{\be/\al}$} is split in $\h\Cl$ since $Z_{\al/0}\lraa{z} Z_{\be/0}$
is split in $\h\Cl$ by assumption, and in which $Z_{\ga/0}\lraa{z} Z_{\ga/\al}$ is split in $\h\Cl$ since $Z_{\al/0}\lraa{z} Z_{\ga/0}$
is split in $\h\Cl$ by assumption.
\qed

So the folding operation passes to an operation
\[
\barcl
\ulk{\Cl^+((2\b\De_n)^\#)} & \lrafl{30}{\ul{\ffk}_n} & \ulk{\Cl^+(\ol{\rh\disj\De_n}^{\,\#})} \\
X                          & \lramapsfl{0}{}         & X\ul{\ffk}_n \\
\ea
\]
on the stable categories.

\pagebreak

\paragraph{An example: folding from $\b\De_5^\#$ to $\b\De_3^\#$}
\indent

\bq
Let $D = \De_2$. Note that \mb{$2\b\De_2\iso\b\De_5$}. Let $X\in\Ob\,\Cl^+((2\b\De_2)^\#)$, depicted as follows. 
\[
\hspace*{-10mm}
\xymatrix@C3mm{
& & & & & & & & & & \\
&  &  &  & 0\ar[r]  & X_{1^{+1}/0^{+1}}\ar[u]\ar[r]^{x}\ar@{}[ur]|{+} & X_{2^{+1}/0^{+1}}\ar[u]\ar[r]^{x}\ar@{}[ur]|{+} & 
  X_{0^{+2}/0^{+1}}\ar[u]\ar[r]^{x}\ar@{}[ur]|{+} & X_{1^{+2}/0^{+1}}\ar[u]\ar[r]^{x}\ar@{}[ur]|{+} & X_{2^{+2}/0^{+1}}\ar[u]\ar[r]\ar@{}[ur]|{+} & 
  0\ar[u] \\
&  &  & 0\ar[r] & X_{0^{+1}/2}\ar[u]\ar[r]^{x}\ar@{}[ur]|{+} & X_{1^{+1}/2}\ar[u]^{x}\ar[r]^{x}\ar@{}[ur]|{+} & 
  \fbox{$X_{2^{+1}/2}$}\ar[u]^{x}\ar[r]^{x}\ar@{}[ur]|{+} & X_{0^{+2}/2}\ar[u]^{x}\ar[r]^{x}\ar@{}[ur]|{+} & X_{1^{+2}/2}\ar[u]^{x}\ar[r]\ar@{}[ur]|{+} & 
  0\ar[u] \\
&  & 0\ar[r] & X_{2/1}\ar[u]\ar[r]^{x}\ar@{}[ur]|{+} & X_{0^{+1}/1}\ar[u]^{x}\ar[r]^{x}\ar@{}[ur]|{+} &
  \fbox{$X_{1^{+1}/1}$}\ar[u]^{x}\ar[r]^{x}\ar@{}[ur]|{+} & X_{2^{+1}/1}\ar[u]^{x}\ar[r]^{x}\ar@{}[ur]|{+} & X_{0^{+2}/1}\ar[u]^{x}\ar[r]\ar@{}[ur]|{+} & 
  0\ar[u] \\
& 0\ar[r] & X_{1/0}\ar[u]\ar[r]^{x}\ar@{}[ur]|{+} & X_{2/0}\ar[u]^{x}\ar[r]^{x}\ar@{}[ur]|{+} & 
  \fbox{$X_{0^{+1}/0}$}\ar[u]^{x}\ar[r]^{x}\ar@{}[ur]|{+} & X_{1^{+1}/0}\ar[u]^{x}\ar[r]^{x}\ar@{}[ur]|{+} & X_{2^{+1}/0}\ar[u]^{x}\ar[r]\ar@{}[ur]|{+} 
  & 0\ar[u] \\
0\ar[r] & X_{0/2^{-1}}\ar[u]\ar[r]^{x}\ar@{}[ur]|{+} & X_{1/2^{-1}}\ar[u]^{x}\ar[r]^{x}\ar@{}[ur]|{+} & X_{2/2^{-1}}\ar[u]^{x}\ar[r]^{x}\ar@{}[ur]|{+} & 
 X_{0^{+1}/2^{-1}}\ar[u]^{x}\ar[r]^{x}\ar@{}[ur]|{+} & X_{1^{+1}/2^{-1}}\ar[u]^{x}\ar[r]\ar@{}[ur]|{+} & 0\ar[u] \\
\ar[u]\ar@{}[ur]|{+}  & \ar[u]\ar@{}[ur]|{+} & \ar[u]\ar@{}[ur]|{+} & \ar[u]\ar@{}[ur]|{+} & \ar[u]\ar@{}[ur]|{+} & \ar[u] \\
}
\]
Note that $\rh\disj\De_2\iso\De_3$. Folding turns $X$ into $X\ffk_2\in\Ob\,\Cl^+(\ol{\rh\disj\De_2}^{\,\#})$, depicted as follows. 
\[
\hspace*{-10mm}
\xymatrix{
& & & & & 0\ar[r] & \\
& & & & 0\ar[r] &  X_{2^{+2}/2^{+1}}\ar[u]\ar[r]\ar@{}[ur]|{+} & \\
& & & 0\ar[r] &  X_{2^{+1}/1^{+1}}\dk X_{1^{+2}/2}\ar[r]^(0.6){\rsmatze{x}{-x}}\ar[u]\ar@{}[ur]|{+} & X_{1^{+2}/1^{+1}}\ar[u]^{x}\ar[r]\ar@{}[ur]|{+} & \\
& & 0\ar[r] & X_{1^{+1}/0^{+1}}\dk X_{0^{+2}/1}\ar[r]^{\smatzz{x}{0}{0}{x}}\ar[u]\ar@{}[ur]|{+}
        & X_{2^{+1}/0^{+1}}\dk X_{0^{+2}/2}\ar[r]^(0.6){\rsmatze{x}{-x}}\ar[u]^{\smatzz{x}{0}{0}{x}}\ar@{}[ur]|{+} 
        & X_{0^{+2}/0^{+1}}\ar[u]^{x}\ar[r]\ar@{}[ur]|{+} &  \\ 
& 0\ar[r] & \fbox{$X_{0^{+1}/0}$}\ar[r]^{x}\ar[u]\ar@{}[ur]|{+} & \fbox{$X_{1^{+1}/1}$}\ar[r]^{x}\ar[u]^{\smatez{x}{x}}\ar@{}[ur]|{+}
        & \fbox{$X_{2^{+1}/2}$}\ar[r]\ar[u]^{\smatez{x}{x}}\ar@{}[ur]|{+} & 0\ar[u] \\
0\ar[r] & X_{2/2^{-1}}\ar[r]^(0.35){\smatez{x}{x}}\ar[u]\ar@{}[ur]|{+} 
        & X_{2/0}\dk X_{0^{+1}/2^{-1}}\ar[r]^{\smatzz{x}{0}{0}{x}}\ar[u]^{\rsmatze{x}{-x}}\ar@{}[ur]|{+} 
        & X_{2/1}\dk X_{1^{+1}/2^{-1}}\ar[r]\ar[u]^{\rsmatze{x}{-x}}\ar@{}[ur]|{+} & 0\ar[u] \\
\ar[u]\ar@{}[ur]|{+}  & \ar[u]\ar@{}[ur]|{+} & \ar[u]\ar@{}[ur]|{+} & \ar[u] &         \\
}
\]
\eq

\pagebreak

\subsubsection{A definition of Heller triangulated categories and strictly exact functors}

Recall that $\Cl$ is a weakly abelian category, and that $\TTT = (-)^{+1}$ is an automorphism of $\Cl$.

Suppose given $n\ge 0$. We have introduced the automorphisms
$$
\ba{rclcl}
\ulk{\Cl^+(\b\De_n^\#)} & \mraisofl{30}{\ulk{\Cl^+(\TTT_{\! n})}}  & \ulk{\Cl^+(\b\De_n^\#)} & \Icm & \mb{(outer shift; \S\ref{SecNatInP})} \\
X                       & \mramaps                                 & {[X]^{+1}}              &      & \vspace*{5mm}\\
\ulk{\Cl^+(\b\De_n^\#)} & \mraisofl{30}{\ulk{\TTT^+(\b\De_n^\#)}}  & \ulk{\Cl^+(\b\De_n^\#)} &      & \mb{(inner shift; \S\ref{SecNatInC})}\\
X                       & \mramaps                                 & {[X^{+1}]}\; .          &      & \\
\ea
$$

\bq
The outer shift shifts the whole diagram $X\in\Ob\,\ulk{\Cl^+(\b\De_n^\#)}$ one step downwards -- the object $X_{\al^{+1}/\be}$ is the entry of $[X]^{+1}$ at position $\be/\al$.

The inner shift applies the given shift automorphism $(-)^{+1}$ of $\Cl$ entrywise to a diagram 
\mb{$X\in\Ob\,\ulk{\Cl^+(\b\De_n^\#)}$.}
\eq

Furthermore, we write $[X^{+a}]^{+b} := X \TTT^+(\b\De_n^\#)^a\,\Cl^+(\TTT_{\! n})^b = X \Cl^+(\TTT_{\! n})^b\, \TTT^+(\b\De_n^\#)^a$ for $a,\, b\,\in\,\Z_{\ge 0}$ and $X\in\Ob\,\Cl^+(\b\De_n^\#)$; 
similarly for $a,\, b\,\in\,\Z$. Likewise in the stable case.

\begin{Definition}
\label{Def4}\rm
\Absit
\begin{itemize}
\item[(i)] \ \vspace*{-6.5mm} \\ 
        \fbox{\begin{tabular}{p{15cm}}
        A {\it Heller triangulation} on $(\Cl,\TTT)$ is a tuple of isomorphisms of functors
        $$
        \tht \;\; =\;\; \Big(\;\;\ulk{\Cl^+(\TTT_{\! n})} \;\;\lraisoa{\tht_n}\;\; \ulk{\TTT^+(\b\De_n^\#)} \;\; \Big)_{\! n\ge 0} 
             \;\; =\;\; \Big(\;\;[-]^{+1} \;\;\lraisoa{\tht_n}\;\; [-^{+1}]\;\;\Big)_{\! n\ge 0}
        $$
        such that 
        $$
        \ulk{p}^\#\st \tht_m \;\; =\;\; \tht_n\st \ulk{p}^\#
        \leqno (\ast)
        $$
        for all $n,\, m\,\ge\, 0$ and all periodic monotone maps $\b\De_n\llaa{p}\b\De_m$ in $\b\Deltab$, and such that
        $$
        \ul{\ffk}_n\st\tht_{n+1} \;\; =\;\; \tht_{2n+1}\st \ul{\ffk}_n
        \leqno (\ast\ast)
        $$
        for all $n\ge 0$. \\
        \end{tabular}}

        \bq
        Note that given $n\ge 0$, the isomorphism $\tht_n$ consists of isomorphisms
        $$
        [X]^{+1} \;\;\lraisoa{X\tht_n} \;\; [X^{+1}]
        $$%
        in the stable category $\ulk{\Cl^+(\b\De_n^\#)}$ of $n$-pretriangles, where $X$ runs over the set 
        $\Ob\,\ulk{\Cl^+(\b\De_n^\#)}$ of $n$-pretriangles.\vspace*{2mm}

        Condition $(\ast)$ asserts that the following diagram commutes in $\Add$ for all $n,\, m\,\ge\, 0$ and all periodic monotone maps $\b\De_n\llaa{p}\b\De_m$ in $\b\Deltab$. 
        $$
        \xymatrix{
        \ulk{\Cl^+(\b\De_n^\#)} \ar[rrr]^{\ul{p}^\# \; = \;\ulk{\Cl^+(p^\#)}}
              \ar@/_4mm/[dd]_{[-]^{+1}}="phf" 
              \ar@/^4mm/[dd]^{[-^{+1}]}="phg" 
        & & & \;\ulk{\Cl^+(\b\De_m^\#)} 
              \ar@/_4mm/[dd]_{[-]^{+1}}="psf"
              \ar@/^4mm/[dd]^{[-^{+1}]}="psg" \\
        & & & & \\
        \ulk{\Cl^+(\b\De_n^\#)} \ar[rrr]^{\ul{p}^\# \; =\;\ulk{\Cl^+(p^\#)}} & & & \;\ulk{\Cl^+(\b\De_m^{\,\#})}
        \ar@2 "phf"+<6.0mm,0mm>;"phg"+<-5.5mm,0mm>^{\tht_n\ru{-0.8}}
        \ar@2 "psf"+<6.0mm,0mm>;"psg"+<-5.5mm,0mm>^{\tht_m\ru{-0.8}}
        }
        $$ 

        Condition $(\ast\ast)$ asserts that the following diagram commutes in $\Add$ for all $n\ge 0$. 
        $$
        \xymatrix{
        \ulk{\Cl^+(\b\De_{2n+1}^\#)} \ar[rrr]^{\ul{\ffk}_n}
              \ar@/_5mm/[dd]_{[-]^{+1}}="phf" 
              \ar@/^5mm/[dd]^{[-^{+1}]}="phg" 
        & & & \;\ulk{\Cl^+(\b\De_{n+1}^\#)} 
              \ar@/_5mm/[dd]_{[-]^{+1}}="psf"
              \ar@/^5mm/[dd]^{[-^{+1}]}="psg" \\
        & & & & \\
        \ulk{\Cl^+(\b\De_{2n+1}^\#)} \ar[rrr]^{\ul{\ffk}_n} & & & \;\ulk{\Cl^+(\b\De_{n+1}^{\,\#})}
        \ar@2 "phf"+<6mm,0mm>;"phg"+<-6mm,0mm>^{\tht_{2n+1}\ru{-0.8}}
        \ar@2 "psf"+<6mm,0mm>;"psg"+<-6mm,0mm>^{\tht_{n+1}\ru{-0.8}}
        }
        $$ 
        \eq

\item[(ii)] \ \vspace*{-6.5mm} \\ 
        \fbox{\begin{tabular}{p{15cm}}
        Given a Heller triangulation $\tht$ on $(\Cl,\TTT)$, we use the following terminology. 

        \begin{tabular}{rp{12cm}}
        (1) & The triple $(\Cl,\TTT,\tht)$ forms a {\it Heller triangulated category,} usually just denoted by $\Cl$.\\
        (2) & Given $n\ge 0$, an {\it $n$-triangle} is an object $X$ of $\Cl^+(\b\De_n^\#)$ for which $[X^{+1}] = [X]^{+1}$ in $\Ob\,\Cl^+(\b\De_n^\#)$, and for which
              $$
              X\tht_n \; =\; 1_{[X]^{+1}} \; =\; 1_{[X^{+1}]}\Icm\text{(equality in $\ulk{\Cl^+(\b\De_n^\#)}$)}\; .
              $$ 
              \vspace*{-4mm}\\
            &  A {\it morphism of $n$-triangles} is a morphism $X\lraa{u} Y$ in $\Cl^+(\b\De_n^\#)$ between $n$-triangles $X$ and $Y$ such that 
               $[u]^{+1} = [u^{+1}]$.\\
            & The category of $n$-triangles and morphisms of $n$-triangles is denoted by $\Cl^{+,\,\tht = 1}(\b\De_n^\#)\ru{4.5}$. 
        \end{tabular} 
        \end{tabular}}

        \bq
        In the notation $\Cl^{+,\,\tht = 1}(\b\De_n^\#)$, the index ``$\;\tht\!\!=\!\!1\;$'' is to be read as a symbol, not as an actual equation. 

        The subcategory of $n$-triangles $\Cl^{+,\,\tht = 1}(\b\De_n^\#)$ in the category of $n$-pretriangles $\Cl^+(\b\De_n^\#)$
        is not full in general. 
        \eq

\item[(iii)] \ \vspace*{-6.5mm} \\ 
        \fbox{\begin{tabular}{p{15cm}}
        An additive functor $\Cl\lraa{F}\Cl'$ between Heller triangulated categories $(\Cl,\TTT,\tht)$ and $(\Cl',\TTT',\tht')$ is called 
        {\it strictly exact} if the following conditions hold. 

        \begin{tabular}{rp{12cm}}
        (1) & $F\TTT' = \TTT F$.  \\
        (2) & $F$ is subexact; cf.\ \S\ref{SecNatInC}. \\
        (3) & We have      
        $$
        \tht_n\st \ulk{F^+(\b\De_n^\#)} \;\; =\;\; \ulk{F^+(\b\De_n^\#)}\st \tht'_n
        \leqno (\ast\!\ast\!\ast)
        $$ 
        for all $n\ge 0$.
        \end{tabular} 
        \end{tabular}}

        \bq
        Such a functor $F$ is called {\sf strictly} exact because of the {\sf equality} in (1).
 
        Condition $(\ast\!\ast\!\ast)$ asserts that the following diagram commutes in $\Add$ for all $n\ge 0$.
        $$
        \xymatrix{
        \ulk{\Cl^+(\b\De_n^\#)} \ar[rrr]^{\ulk{F^+(\b\De_n^\#)}}\ar@/_5mm/[dd]_{[-]^{+1}}="phf" \ar@/^5mm/[dd]^{[-^{+1}]}="phg" 
                        & & & \;\ulk{\Cl'^+(\b\De_n^\#)} \ar@/_5mm/[dd]_{[-]^{+1}}="psf"\ar@/^5mm/[dd]^{[-^{+1}]}="psg" \\
        & & & & \\
        \ul{\Cl^+(\b\De_n^\#)} \ar[rrr]^{\ul{F^+(\b\De_n^\#)}} & & & \;\ul{\Cl'^+(\b\De_n^\#)}
        \ar@2 "phf"+<6mm,0mm>;"phg"+<-6mm,0mm>^{\tht_n\ru{-0.8}}
        \ar@2 "psf"+<6mm,0mm>;"psg"+<-6mm,0mm>^{\tht'_n\ru{-0.8}}
        }
        $$ 
        \eq
\end{itemize} 
\end{Definition}

\bq
To summarise Definition \ref{Def4} roughly, a Heller triangulation is an isomorphism $\tht$ from the outer shift to the inner shift, varying with $\De_n$, and compatible with 
folding. An $n$-triangle is a periodic \mb{$n$-pretriangle} at which $\tht$ is an identity. A strictly exact functor respects the weakly abelian structure and is compatible with shift and $\tht$.
\eq

\bq
Note that if $\tht$ is a Heller triangulation on $(\Cl,\TTT)$, so is $-\tht$.
\eq

\bq
Definition \ref{Def4} would make sense for periodic, but not necessarily linearly ordered posets, generalising $\b\De_n$. But then it is unknown whether, and, it seems to the author, not
very probable that the stable category of a Frobenius category is triangulated in this generalised sense. More specifically, it seems to be impossible to generalise Proposition \ref{PropFull} below
accordingly, which is the technical core of our approach.
\eq

\bq
\begin{Question}
\label{Question1}\rm
Does there exist an additive functor $\Cl\lraa{F}\Cl'$ between Heller triangulated categories that, in Definition \ref{Def4}.(iii), satisfies (1) and (2), but (3) only for $n\le 2\,$? If $F$ is an 
identity, this amounts to asking for the existence of two Heller triangulations $\tht$ and $\tht'$ on $(\Cl,\TTT)$, $\Cl$ weakly abelian, $\TTT$ automorphism of $\Cl$, such that $\tht_n = \tht'_n$ 
only for $n\le 2$.
\end{Question}
\eq

%% file: h3equ.tex
\section{Some equivalences}
\label{SecEqu}

Suppose given $n\ge 0$. Suppose given a weakly abelian category $\Cl$, together with an automorphism $\TTT:\Cl\lra\Cl$, $X\lramaps X^{+1}$. 
Concerning the Freyd category $\h\Cl$ of $\Cl$, we refer to \S\ref{SecFreyd}.

\bq
We shall show in Proposition \ref{PropR5} that the functor $\ulk{\Cl(\b\De_n^\#)}\lra\ulk{\Cl(\dDe_n)}$, induced on the stable categories 
by restriction from $\b\De_n^\#$ to $\dDe_n := [1,n]$, is an equivalence. 
\eq

\subsection{Some notation}
\label{SecSomNot}

\subsubsection{Some posets}
\label{SecSomPos}

Let $\dDe_n := \De_n\ohne\{ 0\} = [1,n]$, considered as a linearly ordered set.
We have an injection $\dDe_n\lra\b\De_n^\#$, $i\lramaps i/0$, and identify $\dDe_n$ with its image in $\b\De_n^\#$.

We define two subposets of $\b\De_n^\#$ by
\[
\barcl
\b\De_n^\tru & := & \{\be/\al\in \b\De_n^\# \; :\; 0\le\al \} \\
\b\De_n^\trd & := & \{\be/\al\in\b\De_n^\# \; :\;  \al\le 0\}\; . \\
\ea
\]
Then $\dDe_n = (\b\De_n^\tru\cap \b\De_n^\trd)\ohne\{0/0,\,0^{+1}/0\}$.

\begin{center}
\begin{picture}(600,400)
\put(   0,   0){\line(1,1){400}}
\put( 200,   0){\line(1,1){400}}
\put( 200, 200){\line(1,0){200}}
\put( 300, 225){$\b\De_n^\tru$}
\put( 250, 150){$\b\De_n^\trd$}
\put(  55, 185){$\dDe_n$}
\put( 480, 200){\vector(-1,0){50}}
\put( 120, 200){\vector(1,0){50}}
\end{picture}
\end{center}

\subsubsection{Fixing parametrisations $\ka^\tru$, $\ka^\trd$}
\label{ParKappa}

There exists a bijective morphism $\b\De_n^\tru\lra\Z_{\ge 0}$ of posets (``refining the partial to a linear order''). We fix such a morphism and denote by 
$\Z_{\ge 0}\lraa{\ka^\tru} \b\De_n^\tru\ru{-2}\ru{5}$ its 
inverse (as a map of sets; in general, $\ka^\tru$ is not monotone). So whenever $\ka^\tru(\ell) \le \ka^\tru(\ell')$, then $\ell\le \ell'$. In particular, $\ka^\tru(0) = 0/0$.

There exists a bijective morphism $\b\De_n^\trd\lra\Z_{\le 0}$ of posets. We fix such a morphism and denote by $\Z_{\ge 0}\lraa{\ka^\trd} \b\De_n^\trd\ru{-2}\ru{5}$
its inverse (as a map of sets). So whenever $\ka^\trd(\ell) \le \ka^\trd(\ell')$, then $\ell\le\ell'$. In particular, $\ka^\trd(0) = 0^{+1}/0$.

\subsubsection{The categories $\h\Cl^{+,\ast}(\b\De_n^\tru)$, $\Cl^+(\b\De_n^\tru)$ etc.}

Let $\Al$ be an abelian category, and let $\Bl\tm\Al$ be a full subcategory. Let $E\tm\b\De_n^\#$ be a full subposet. 

\bq
For example, for $E$ we may take the subposets $\b\De_n^\tru$, $\b\De_n^\trd$ or $\b\De_n^\tru\cap \b\De_n^{\trd,+1}$ of $\b\De_n^\#$.

Moreover, for example, we may take $\Al = \h\Cl$ and for $\Bl$ either $\h\Cl$ or $\Cl$.
\eq

Let $\Bl^{+,\ast}(E)$ be the full subcategory of $\Bl(E)$ defined by
\[
\Ob\,\Bl^{+,\ast}(E) \; :=\; \left\{ \;\; X\in \Ob\,\Bl(E)\;\;\; :\; 
\mb{
\begin{tabular}{l}
For all $\de^{-1}\le\al\le\be\le\ga\le\de\le\al^{+1}$ in $\b\De_n$ \\
such that $\ga/\al$, $\ga/\be$, $\de/\al$ and $\de/\be$ are in $E$, \\
the quadrangle \\
      $\xymatrix{
      X_{\ga/\be}\ar[r]^x                       & X_{\de/\be} \\
      X_{\ga/\al}\ar[r]^x\ar[u]^x\ar@{}[ur]|{+} & X_{\de/\al}\ar[u]_x \\      
      }$ \\
is a weak square (as indicated by $+$). \\
\end{tabular}
}\right\} \;\; .
\]
The symbol $\ast$ should remind us of the fact that we still allow $X_{\al/\al}$ resp.\ $X_{\al^{+1}/\al}$ to be arbitrary for $\al\in \b\De_n$ such that
$\al/\al\in E$ resp.\ $\al^{+1}/\al\in E$.

In turn, let $\Bl^+(E)$ be the full subcategory of $\Bl^{+,\ast}(E)$ defined by
\[
\Ob\,\Bl^+(E) \; :=\; \left\{ \;\; X\in\Ob\,\Bl^{+,\ast}(E)\;\;\; :\; 
\mb{
\begin{tabular}{l}
$X_{\al/\al} = 0$ for $\al\in \b\De_n$ such that $\al/\al\in E\;$, and\\
$X_{\al^{+1}/\al} = 0$ for $\al\in \b\De_n$ such that $\al^{+1}/\al\in E$. \\
\end{tabular}
}\right\} \;\; .
\]

\subsubsection{Reindexing}
\label{ParReindex}

Given a subposet $E\tm \b\De_n^\#$, we have a {\it reindexing equivalence}
\[
\ba{lcl}
\Cl(E)   & \lraiso  & \Cl(E^{+1}) \\
X        & \lramaps & X^{(-1)} \\
X^{(+1)} & \llamaps & X \\
\ea
\]
defined by 
\[
(X^{(-1)})_{\be/\al} \; := \; X_{(\be/\al)^{-1}} \= X_{\al/\be^{-1}} \; ,
\] 
where $\be/\al\in E^{+1}$; and inversely by
\[
(X^{(+1)})_{\be/\al} \; :=\; X_{(\be/\al)^{+1}} \= X_{\al^{+1}/\be} \; , 
\]
where $\be/\al\in E$. This equivalence restricts to an equivalence between $\Cl^+(E)$ and $\Cl^+(E^{+1})$. 

For instance, if $E = \b\De_n^\#$, then $X^{(+1)} = [X]^{+1}$. The outer shift and reindexing will play different roles, and so we distinguish in notation.

\subsection{Density of the restriction functor from $\b\De_n^\#$ to $\dDe_n$} 

\subsubsection{Upwards and downwards spread}
\label{SubsubUpDown}

Let the {\it upwards spread} $S^\tru$ be defined by
\[
\barcl
\h\Cl(\dDe_n) & \lraa{S^\tru} & \h\Cl^+(\b\De_n^\tru) \\
X             & \lramaps      & XS^\tru\; , \\
\ea
\]
where $XS^\tru$ is given by
\[
\ba{lcll}
(XS^\tru)_{0/0}     & := & 0                            & \\
(XS^\tru)_{\be/0}   & := & X_\be                        & \mb{for $\be\in \dDe_n$} \\
(XS^\tru)_{\be/\al} & := & \Cokern(X_\al\lraa{x} X_\be) & \mb{for $\al,\,\be\,\in\, \dDe_n$ with $\al\le\be$} \\
(XS^\tru)_{\be/\al} & := & 0                            & \mb{for $\al,\,\be\,\in\, \b\De_n$ with $0^{+1}\le\be\le\al^{+1}\le\be^{+1}$} \; ,\\
\ea
\]
the diagram being completed with the induced morphisms between the cokernels and zero morphisms elsewhere.

This construction is functorial in $X$. The functor $S^\tru$ is left adjoint to the restriction functor from $\h\Cl^+(\b\De_n^\tru)$ to $\h\Cl(\dDe_n)$, with unit being the identity, 
i.e.\ $X = XS^\tru|_{\dDe_n}$.

Dually, let the {\it downwards spread} $S^\trd$ be 
\[
\barcl
\h\Cl(\dDe_n) & \lraa{S^\trd} & \h\Cl^+(\b\De_n^\trd) \\
X             & \lramaps      & XS^\trd\; , \\
\ea
\]
where $XS^\trd$ is given by
\[
\ba{lcll}
(XS^\trd)_{0^{+1}/0}     & := & 0                            & \\
(XS^\trd)_{\al/0}        & := & X_\al                        & \mb{for $\al\in \dDe_n$} \\
(XS^\trd)_{\al/\be^{-1}} & := & \Kern(X_\al\lraa{x} X_\be)   & \mb{for $\al,\,\be\,\in\, \dDe_n$ with $\al\le\be$} \\
(XS^\trd)_{\al/\be^{-1}} & := & 0                            & \mb{for $\al,\,\be\,\in\,\b\De_n$ with $\al^{-1}\le\be^{-1}\le\al\le 0$}\; ,\\
\ea
\]
the diagram being completed with the induced morphisms between the kernels and zero morphisms elsewhere.

This construction is functorial in $X$. The functor $S^\trd$ is right adjoint to the restriction functor from $\h\Cl^+(\b\De_n^\trd)$ to $\h\Cl(\dDe_n)$, with counit being the identity,
i.e.\ $XS^\trd|_{\dDe_n} = X$.

\subsubsection{Resolutions}

\paragraph{A stability under pointwise pushouts and pullbacks}
\indent

Let $E\tm \b\De_n^\#$ be a full subposet. Moreover, assume that $E$ is a {\it convex} subposet, i.e.\ that whenever given $\xi,\,\ze\in E$ and $\eta\in\b\De_n^\#$ such that $\xi\le\eta\le\ze$, then
$\eta\in E$.

An element $\de/\be\in E$ is {\it on the left boundary of $E$} if we may conclude from $\ga/\be\in E$ and $\ga\le\de$ that $\ga = \de$.
It is {\it on the lower boundary of $E$} if we may conclude from $\de/\al\in E$ and $\al\le\be$ that $\al = \be$.
 
An element $\ga/\al\in E$ is {\it on the right boundary of $E$} if we may conclude from $\de/\al\in E$ and $\ga\le\de$ that $\ga = \de$.
It is {\it on the upper boundary of $E$} if we may conclude from $\ga/\be\in E$ and $\al\le\be$ that $\al = \be$.

Let $\Al$ be an abelian category. Concerning pointwise pullbacks and pointwise pushouts, we refer to \S\ref{SubsecPPP}.

\begin{Lemma}
\label{LemR1}
Suppose given $\eps\in E$ and an object $X$ of $\Al^{+,\ast}(E)$.
\begin{itemize}
\item[{\rm (1)}] Given a monomorphism $X_\eps\lraa{x'} X'$ in $\Al$, the pointwise pushout $X\ind^{x'}$ of $X$ along $x'$ is an object of $\Al^{+,\ast}(E)$ again.
\item[{\rm (2)}] Given an epimorphism $X_\eps\llaa{x'} X'$ in $\Al$, the pointwise pullback $X\ind_{x'}$ of $X$ along $x'$ is an object of $\Ob\Al^{+,\ast}(E)$ again.
\item[{\rm (3)}] Suppose that $\eps$ is on the left boundary or on the lower boundary of $E$. Given a morphism $X_\eps\lraa{x'} X'$ in $\Al$, the pointwise 
pushout $X\ind^{x'}$ of $X$ along $x'$ is an object of $\Al^{+,\ast}(E)$ again.
\item[{\rm (4)}] Suppose that $\eps$ is on the right boundary or on the upper boundary of $E$. Given a morphism $X_\eps\llaa{x'} X'$ in $\Al$, the pointwise pullback 
$X\ind_{x'}$ of $X$ along $x'$ is an object of $\Ob\Al^{+,\ast}(E)$ again.
\end{itemize}
\end{Lemma}

{\it Proof.} Ad (1). First we remark that by Lemma \ref{Lem1_5}, the quadrangle 
\mb{$(X_{\be/\al}\, ,\;X_{\de/\ga}\, ,\;(X\ind^{x'})_{\be/\al}\, ,\; (X\ind^{x'})_{\de/\ga})$}
is a pushout for $\eps\le\be/\al\le\de/\ga$ in $E$. 

We have to show that the quadrangle of $\mb{$X\ind^{x'}$}$ on $(\ga/\al,\de/\al,\ga/\be,\de/\be)$, where 
\mb{$\de^{-1}\le\al\le\be\le\ga\le\de\le\al^{+1}$} in $\b\De_n$,
is a weak square, provided its corners have indices in $E$. Using Lemmata \ref{Lem1}, \ref{Lem1_7} and convexity of $E$, we are reduced to the case $\eps\le\ga/\al$. In this case, the 
assertion follows by Lemma \ref{Lem1_5}. 

Ad (3). Here we need only Lemma \ref{Lem1} and convexity of $E$ to reduce to the case $\eps\le\ga/\al$, the rest of the argument is as in (1). Hence the morphism $x'$ may be arbitrary.
\qed

\paragraph{Upwards and downwards resolution}
\indent

\begin{Remark}
\label{RemR2}
\indent\rm
\begin{itemize}
\item[(1)] Given a direct system $X_0\lra X_1\lra X_2\lra\cdots$ in $\h\Cl^{+,\ast}(\b\De_n^\tru)$ such that its restriction to any finite full
subposet $E\tm \b\De_n^\tru$ eventually becomes constant, then the direct limit $\limd_i X_i$ exists in $\h\Cl^{+,\ast}(\b\De_n^\tru)$.
\item[(2)] Given an inverse system $X_0\lla X_1\lla X_2\lla\cdots$ in $\h\Cl^{+,\ast}(\b\De_n^\trd)$ such that its restriction to any finite
subposet $E\tm \b\De_n^\trd$ eventually becomes constant, then the inverse limit $\limp_i X_i$ exists in $\h\Cl^{+,\ast}(\b\De_n^\trd)$.
\end{itemize}
\end{Remark}

For $k\ge 0$, we let 
\[
\barcl
\Ob\h\Cl^{+,\ast}(\b\De_n^\tru) & \lraa{R^\tru_k} & \Ob\h\Cl^{+,\ast}(\b\De_n^\tru) \\
X                               & \lramaps        & X\ind^{x'(k)}\; , \\
\ea
\]
where
\[
x'(k) \; := \; \left\{
\ba{ll}
(X_{\ka^\tru(k)} \lra 0)                            & \mb{if $\ka^\tru(k) \in \{\al/\al,\,\al^{+1}/\al\}$ for some $\al\in \b\De_n$ with $0\le\al$} \\
(X_{\ka^\tru(k)} \lraa{\io} X_{\ka^\tru(k)}{\sf I}) & \mb{if $\ka^\tru(k) = \be/\al$ for some $\al,\,\be\,\in\, \b\De_n$ with $0\le\al < \be < \al^{+1}$} \\
\ea
\right.
\]

Define the {\it upwards resolution} map by
\[
\barcl
\Ob\h\Cl^{+,\ast}(\b\De_n^\tru) & \lraa{R^\tru} & \Ob\,\Cl^+(\b\De_n^\tru) \\
X                               & \lramaps      & XR^\tru \; :=\; \dis{\limd_m}\, X R^\tru_0 \cdots R^\tru_m\; , \\
\ea
\]
the direct system being given by the transition morphisms 
\[
X R^\tru_0 \cdots R^\tru_m \;\lraa{i}\; (X R^\tru_0 \cdots R^\tru_m) R^\tru_{m+1}\; .
\]
We have $X R^\tru = X$ for $X\in\Ob\,\Cl^+(\b\De_n^\tru)$.

\bq
 Note that we apparently cannot turn the upwards resolution into a functor unless we are in a particular case in which the map $\III$ on objects can be turned into a functor.
\eq

Dually, for $k\ge 0$, we let 
\[
\barcl
\Ob\h\Cl^{+,\ast}(\b\De_n^\trd) & \lraa{R^\trd_k} & \Ob\h\Cl^{+,\ast}(\b\De_n^\trd) \\
X                               & \lramaps        & X\ind_{x''(k)}\; , \\
\ea
\]
where
\[
x''(k) \; :=\; \left\{
\ba{ll}
(X_{\ka^\trd(k)} \lla 0)                            & \mb{if $\ka^\trd(k) \in\{\al/\al,\,\al^{+1}/\al\}$ for some $\al\in \b\De_n$ with $\al\le 0$} \\
(X_{\ka^\trd(k)} \llaa{\pi} X_{\ka^\trd(k)}{\sf P}) & \mb{if $\ka^\trd(k) = \be/\al$ for some $\al,\,\be\,\in\, \b\De_n$ with $\al^{-1} < \be^{-1} < \al \le 0$} \\
\ea
\right.
\]

Define the {\it downwards resolution} map by
\[
\barcl
\Ob\h\Cl^{+,\ast}(\b\De_n^\trd) & \lraa{R^\trd} & \Ob\,\Cl^+(\b\De_n^\trd) \\
X                               & \lramaps      & XR^\trd \; :=\; \dis{\limp_m}\, X R^\trd_0 \cdots R^\trd_m\; , \\
\ea
\]
the inverse system being given by the transition morphisms 
\[
X R^\trd_0 \cdots R^\trd_m \;\llaa{p}\; (X R^\trd_0 \cdots R^\trd_m) R^\trd_{m+1}\; .
\]
We have $X R^\trd = X$ for $X\in\Ob\,\Cl^+(\b\De_n^\trd)$.

\pagebreak[3]

\begin{Lemma}
\label{LemR2_5}\indent
\begin{itemize}
\item[{\rm (1)}] Given a morphism $Y\lraa{g} X$ in $\h\Cl^+(\b\De_n^\tru)$ with $X\in\Ob\,\Cl^+(\b\De_n^\tru)$, there exists a factorisation
\[
(Y\lraa{g} X) \= (Y\lra Y R^\tru \lra X)\; .
\]
\item[{\rm (2)}] Given a morphism $Y\llaa{g} X$ in $\h\Cl^+(\b\De_n^\trd)$ with $X\in\Ob\,\Cl^+(\b\De_n^\trd)$, there exists a factorisation
\[
(Y\llaa{g} X) \= (Y\lla Y R^\trd \lla X)\; .
\]
\end{itemize}
\end{Lemma}

{\it Proof.} Ad (1). Since the entries of $X$ are injective in $\h\Cl$ and since $X_{\al/\al} = 0$ and $X_{\al^{+1}/\al} = 0$ for $\al\ge 0$, we obtain, using the universal property of the pointwise 
pushout, a factorisation 
\[
(Y\lraa{g} X) \= (Y\lra Y R_0^\tru\cdots R_m^\tru \lra X)
\]
for every $m\ge 0$, compatible with the transition morphisms, resulting in a factorisation over $Y R^\tru = \limd_m Y R_0^\tru\cdots R_m^\tru$.
\qed

\paragraph{Both-sided resolutions}
\label{ParRes}
\indent

Let the {\it resolution} map
\[
\barcl
\Ob\,\Cl(\dDe_n) & \lraa{R} & \Ob\,\Cl^+(\b\De_n^\#) \\
X                & \lramaps & XR \\
\ea
\]
be defined by gluing an upper and a lower part along $\dDe_n$ as follows.
\[
\barcl
XR|_{\b\De_n^\tru} & := & XS^\tru R^\tru \\
XR|_{\b\De_n^\trd} & := & XS^\trd R^\trd \\
\ea
\] 
This is welldefined, since $XS^\tru R^\tru |_{\dDe_n} = X = XS^\trd R^\trd |_{\dDe_n}$. In particular, we obtain 
\[
XR|_{\dDe_n} \= X\; .
\]

We summarise.

\begin{Proposition}
\label{PropR3}
The restriction functor
\[
\barcl
\Cl^+(\b\De_n^\#) & \lrafl{35}{(-)|_{\dDe_n}}     & \Cl(\dDe_n) \\
Y                 & \lramapsfl{0}{} & Y|_{\dDe_n} \\
\ea
\]
is strictly dense, i.e.\ it is surjective on objects.
\end{Proposition}

\subsection{Fullness of the restriction functor from $\b\De_n^\#$ to $\dDe_n$}

\begin{Proposition}
\label{PropFull}
The restriction functors
\[
\barcl
\Cl^+(\b\De_n^\#) & \lrafl{35}{(-)|_{\dDe_n}} & \Cl(\dDe_n) \\
Y                 & \lramapsfl{0}{}           & Y|_{\dDe_n} \\
\ea
\]
and
\[
\barcl
\Cl^+(\b\De_n^\tru) & \lrafl{35}{(-)|_{\dDe_n}} & \Cl(\dDe_n) \\
Y                   & \lramapsfl{0}{}           & Y|_{\dDe_n} \\
\ea
\]
are full. 
\end{Proposition}

{\it Proof.} By duality and gluing along $\dDe_n$, it suffices to consider the restriction from $\b\De_n^\tru$ to $\dDe_n$. So suppose given $X,\, Y\,\in\,\Ob\,\Cl^+(\b\De_n^\tru)$ and a morphism 
$X|_{\dDe_n}\lraa{f} Y|_{\dDe_n}$. We have to find a morphism $X\lraa{f^\tru} Y$ such that $f^\tru|_{\dDe_n} = f$.

We construct $f^\tru_{\ka^\tru(\ell)}$ for $\ell\ge 0$ by induction on $\ell$.

At $\ka^\tru(0)$, we let $f^\tru_{\ka^\tru(0)} := 1_0$. Suppose given $\ell\ge 1$. If $\ka^\tru(\ell)\in \dDe_n$, we let $f^\tru_{\ka^\tru(\ell)} := f_{\ka^\tru(\ell)}$. If 
$\ka^\tru(\ell)\in\{ \al/\al,\,\al^{+1}/\al\}$ for some $\al\ge 0$, we let $f^\tru_{\ka^\tru(\ell)} := 1_0$. If $\ka^\tru(\ell) =: \be/\al$ with \mb{$0 < \al < \be < \al^{+1}$}, then we let 
$\al' := \al - 1$ be the predecessor of $\al$ in $\b\De_n$, and we let $\be' := \be - 1$ be the predecessor of $\be$ in $\b\De_n$, using that $\b\De_n$ is {\sf linearly} ordered. We may complete 
the diagram
\[
\xymatrix{
 & Y_{\be'/\al}\ar[rrr]^y & & & Y_{\be/\al} \\
X_{\be'/\al}\ar[ur]^*+<-1mm,-1mm>{\scm f^\tru_{\be'/\al}}\ar[rrr]^x & & & X_{\be/\al} \\ 
 & & & &  \\
 & Y_{\be'/\al'}\ar[rrr]_y|(0.61)\hole\ar[uuu]|(0.63)\hole^y\ar@{}[uuurrr]|+ & & & Y_{\be/\al'}\ar[uuu]_y \\
X_{\be'/\al'}\ar[ur]^*+<-1mm,-1mm>{\scm f^\tru_{\be'/\al'}}\ar[rrr]_x\ar[uuu]^x\ar@{}[uuurrr]|+ & & & X_{\be/\al'}\ar[ur]^*+<-1mm,-1mm>{\scm f^\tru_{\be/\al'}}\ar[uuu]_x \\ 
}
\]
to a commutative cuboid, inserting a morphism $X_{\be/\al}\lrafl{30}{f^\tru_{\be/\al}} Y_{\be/\al}$.\qed

\bq
Since we need the restriction functor $\Cl^+(\b\De_n^\#)\lrafl{30}{(-)|_{\dDe_n}} \Cl(\dDe_n)$ to be full, we are not able to generalise from linearly ordered periodic 
posets to arbitrary periodic posets.
\eq

\subsection{The equivalence between\ \ \ulk{$\Cl^+(\b\De_n^\#)$}\ \ and\ \ \ulk{$\Cl(\dDe_n)$}}

Let $\Cl^{\spl}(\dDe_n)$ be the full subcategory of $\Cl(\dDe_n)$ defined by
\[
\Ob\,\Cl^{\spl}(\dDe_n) \; :=\; \{ X\in\Ob\,\Cl(\dDe_n)\; :\; \mb{$X_\al\lra X_\be$ is split in $\h\Cl$ for all $\al,\,\be\,\in\, \dDe_n$ with $\al\le\be$}\}\; .
\]
We denote the quotient category by 
\[
\ulk{\Cl(\dDe_n)} \; :=\; \Cl(\dDe_n)/\Cl^{\spl}(\dDe_n)\; .
\]

\begin{Proposition}
\label{PropR5}
The functor
\[
\barcl
\ulk{\Cl^+(\b\De_n^\#)} & \lrafl{35}{(-)|_{\dDe_n}} & \ulk{\Cl(\dDe_n)} \\
X                       & \lramapsfl{0}{}           & X|_{\dDe_n}\; , \\
\ea
\]
induced by restriction from $\b\De_n^\#$ to $\dDe_n$, is an equivalence.
\end{Proposition}

{\it Proof.} By Propositions \ref{PropR3} and \ref{PropFull}, we may invoke Lemma \ref{Lem5}. Moreover, Lemma \ref{LemR4} gives the inverse
image of $\Ob\,\Cl^{\spl}(\dDe_n)$ under $\Cl^+(\b\De_n^\#) \lrafl{35}{(-)|_{\dDe_n}} \Cl(\dDe_n)\ru{6}$ as $\Ob\,\Cl^{+,\,\spl}(\b\De_n^\#)$. 

Consider a morphism $X\lraa{f}X'$ in $\Cl^+(\b\De_n^\#)$ such that $(X\lraa{f}X')|_{\dDe_n}$ is zero in $\Cl(\dDe_n)$. We have to prove that it factors over an object of
$\Cl^{+,\,\spl}(\b\De_n^\#)$.

Let $Y^\tru$ be the cokernel in $\h\Cl(\b\De_n^\tru)$ of the counit $X|_{\dDe_n}S^\tru\lra X|_{\b\De_n^\tru}$ at $X|_{\b\De_n^\tru}$. Note that $Y^\tru|_{\dDe_n} = 0$. By Lemma \ref{Lem3_2}, 
we have $Y^\tru\in\Ob\,\h\Cl^+(\b\De_n^\tru)$. Consider the following diagram in $\h\Cl^+(\b\De_n^\tru)$.
\[
\xymatrix{
X|_{\dDe_n}S^\tru\ar[r]\ar[d]_0 & X|_{\b\De_n^\tru}\ar[r]\ar[d]_{f|_{\b\De_n^\tru}} & Y^\tru\ar[r]\ar@{.>}[dl] & Y^\tru R^\tru\ar@{.>}[dll] \\
X'|_{\dDe_n}S^\tru\ar[r]        & X'|_{\b\De_n^\tru}                                &                          &               \\
}
\]
The indicated factorisation 
\[
(X|_{\b\De_n^\tru}\lrafl{35}{f|_{\b\De_n^\tru}}X'|_{\b\De_n^\tru}) \= (X|_{\b\De_n^\tru}\lra Y^\tru\lra X'|_{\b\De_n^\tru})
\]
ensues from the universal property of the cokernel $Y^\tru$. By Lemma \ref{LemR2_5}.(1), we can factorise further to obtain
$$
(X|_{\b\De_n^\tru}\lrafl{35}{f|_{\b\De_n^\tru}}X'|_{\b\De_n^\tru}) \= (X|_{\b\De_n^\tru}\lra Y^\tru R^\tru\lra X'|_{\b\De_n^\tru})\; .
\leqno (\ast^\tru)
$$

Dually, we obtain a factorisation
$$
(X|_{\b\De_n^\trd}\lrafl{35}{f|_{\b\De_n^\trd}}X'|_{\b\De_n^\trd}) \= (X|_{\b\De_n^\trd}\lra Y^\trd R^\trd\lra X'|_{\b\De_n^\trd})
\leqno (\ast^\trd)
$$
for some $Y^\trd\in\Ob\,\h\Cl^+(\b\De_n^\trd)$ such that $Y^\trd|_{\dDe_n} = 0$.

Since $Y^\tru R^\tru|_{\dDe_n} = 0 = Y^\trd R^\trd|_{\dDe_n}\ru{-2}$, there is a unique $N\in\Ob\,\Cl^+(\b\De_n^\#)$ such that $N|_{\b\De_n^\tru} = Y^\tru R^\tru$
and $N|_{\b\De_n^\trd} = Y^\trd R^\trd$. By Lemma \ref{LemR4}, we have $N\in\Ob\,\Cl^{+,\,\spl}(\b\De_n^\#)$.

Moreover, since both factorisations $(\ast^\tru)$ and $(\ast^\trd)$ restrict to the factorisation
\[
(X|_{\dDe_n}\lraa{0} X'|_{\dDe_n}) \= (X|_{\dDe_n}\lra 0\lra X'|_{\dDe_n})
\]
in $\Cl(\dDe_n)$, we may glue to a factorisation
$$
(X\lraa{f} X') \= (X\lra N \lra X')
\leqno (\ast)
$$
that restricts to $(\ast^\tru)$ in $\Cl^+(\b\De_n^\tru)$ and to $(\ast^\trd)$ in $\Cl^+(\b\De_n^\trd)$.\qed

\subsection{Auxiliary equivalences}
\label{SecAuxEq}

\bq
We shall extend the equivalence $\ulk{\Cl^+(\b\De_n^\#)}\lraiso\ulk{\Cl(\dDe_n)}$ to a diagram of equivalences
\[
\ulk{\Cl^+(\b\De_n^\#)}\;\lraiso\;\ulk{\Cl^+(\De_n^{\trud})}\;\lraiso\;\ulk{\Cl(\dDe_n)}\;\lraiso\;\ulk{\h\Cl(\dDe_{n-1})}\; .
\]
\eq

\subsubsection{Factorisation into two equivalences}

Abbreviate $\b\De_n^{\tru,+1} := (\b\De_n^\tru)^{+1}\tm \b\De_n^\#$, $\b\De_n^{\trd,+1} := (\b\De_n^\trd)^{+1}\tm \b\De_n^\#$ and 
$\dDe_n^{+1} := (\dDe_n)^{+1} \tm \b\De_n^\#$. 

Abbreviate $\b\De_n^{\trud} := \b\De_n^\tru\cap \b\De_n^{\trd,+1} = \{ \be/\al \; :\; \al,\,\be\,\in\, \b\De_n,\; 0\le\al\le\be\le0^{+1}\} \tm \b\De_n^\#$.

\begin{center}
\begin{picture}(600,480)
\put(   0,   0){\line(1,1){420}}
\put( 200,   0){\line(1,1){420}}
\put( 200, 200){\line(1,0){200}}
\put( 400, 200){\line(0,1){200}}
\put( 300, 245){$\b\De_n^{\trud}$}
\put(  50, 188){$\dDe_n$}
\put( 250, 150){$\b\De_n^\trd$}
\put( 410, 310){$\b\De_n^{\tru,+1}$}
\put( 480, 200){\vector(-1,0){50}}
\put( 120, 200){\vector(1,0){50}}
\put( 400, 480){\vector(0,-1){50}}
\put( 400, 120){\vector(0,1){50}}
\put( 380,  55){$\dDe_n^{+1}$}
\end{picture}
\end{center}

Let $\Cl^{+,\,\spl}(\b\De_n^{\trud})$ be the full subcategory of $\Cl^+(\b\De_n^{\trud})$ defined by
\[
\Ob\,\Cl^{+,\,\spl}(\b\De_n^{\trud}) 
\; := \; \left\{ X\in\Ob\,\Cl^+(\b\De_n^{\trud}) \; :\; 
\mb{
\begin{tabular}{l} $X_{\ga/\al}\lraa{x} X_{\de/\be}$ is split in $\h\Cl$ \\
for all $\ga/\al,\,\de/\be\,\in\, \b\De_n^{\trud}$\\ 
with $\ga/\al\le\de/\be$\\
\end{tabular}
}\right\}
\]
We denote the quotient category by 
\[
\ulk{\Cl^+(\b\De_n^{\trud})} \; :=\; \Cl^+(\b\De_n^{\trud})/\Cl^{+,\,\spl}(\b\De_n^{\trud})
\]

\begin{Lemma}
\label{LemR6}\Absit
\begin{itemize}
\item[{\rm (1)}] The restriction functors
\[
\ba{rcccl}
\Cl^+(\b\De_n^\#) & \mrafl{35}{(-)|_{\b\De_n^{\trud}}} & \Cl^+(\b\De_n^{\trud}) & \lrafl{35}{(-)|_{\dDe_n}} & \Cl^+(\dDe_n) \\
X                 & \mramaps                           & X|_{\b\De_n^{\trud}}   &                           &               \\
                  &                                    &                     Y  & \lramapsfl{0}{}           & Y|_{\dDe_n}   \\
\ea
\]
are full and strictly dense.

\item[{\rm (2)}] The functors
\[
\ba{rcccl}
\ulk{\Cl^+(\b\De_n^\#)} & \mrafl{35}{(-)|_{\b\De_n^{\trud}}} & \ulk{\Cl^+(\b\De_n^{\trud})} & \lrafl{35}{(-)|_{\dDe_n}} & \ulk{\Cl^+(\dDe_n)} \\
X                       & \mramaps                           & X|_{\b\De_n^{\trud}}         &                           &                     \\
                        &                                    & Y                            & \lramapsfl{0}{}           & Y|_{\dDe_n}\; ,     \\
\ea
\]
induced by restriction, are equivalences.
\end{itemize}
\end{Lemma}

{\it Proof.} Ad (1). The composition
\[
\ba{rcccl}
\Cl^+(\b\De_n^\#) & \lra     & \Cl^+(\b\De_n^{\trud}) & \lra     & \Cl(\dDe_n)     \\
X                 & \lramaps & X|_{\b\De_n^{\trud}}   & \lramaps & X|_{\dDe_n}     \\
\ea
\]
is strictly dense by Proposition \ref{PropR3} and full by Proposition \ref{PropFull}. Therefore, the restriction functor from 
$\Cl^+(\b\De_n^{\trud})$ to $\Cl^+(\dDe_n)$ is full and strictly dense. 

We claim that the restriction functor from $\Cl^+(\b\De_n^\#)$ to $\Cl^+(\b\De_n^{\trud})$ is strictly dense. Let
\[
\barcl
\Cl(\dDe_n^{+1}) & \lraa{S^{\tru,+1}} & \h\Cl^+(\b\De_n^{\tru,+1}) \\
X                & \lramaps           & (X^{(+1)} S^\tru)^{(-1)}\; , \\
\ea
\]
cf.\ \S\ref{ParReindex}. Similarly,
\[
\barcl
\Ob\h\Cl^+(\b\De_n^{\tru,+1}) & \lraa{R^{\tru,+1}} & \Ob\,\Cl^+(\b\De_n^{\tru,+1}) \\
X                             & \lramaps           & (X^{(+1)} R^\tru)^{(-1)} \; .\\
\ea
\]
Given $X\in\Ob\,\Cl^+(\b\De_n^{\trud})$, we may define $X'\in\Ob\,\Cl^+(\b\De_n^\#)$ letting
\[
\ba{lcl}
X'|_{\b\De_n^{\tru,+1}}         & := & X|_{\dDe_n^{+1}} S^{\tru,+1} R^{\tru,+1} \\
X'|_{\b\De_n^{\trud}}           & := & X  \\
X'|_{\b\De_n^\trd}              & := & X|_{\dDe_n} S^\trd R^\trd \; .\\
\ea
\]

We claim that the restriction functor from $\b\De_n^\#$ to $\b\De_n^{\trud}$ is full. Suppose given $X,\, Y\,\in\,\Ob\,\Cl^+(\De_n^\#)$ and a morphism 
$X|_{\b\De_n^{\trud}}\lraa{f} Y|_{\b\De_n^{\trud}}$. By Proposition \ref{PropR5} and a shift,
there exists a morphism $X|_{\De_n^{\tru,+1}}\lraa{f^\tru} Y|_{\De_n^{\tru,+1}}$ such that $f^\tru|_{\dDe_n^{+1}} = f|_{\dDe_n^{+1}}$. By Proposition \ref{PropR5} and by duality, there exists a 
morphism $X|_{\De_n^\trd}\lraa{f^\trd} Y|_{\De_n^\trd}$ such that $f^\trd|_{\dDe_n} = f|_{\dDe_n}$. We may define a morphism $X\lraa{f'} Y$ letting
\[
\ba{lcl}
f'|_{\b\De_n^{\tru,+1}}         & := & f^\tru \\
f'|_{\b\De_n^{\trud}}           & := & f  \\
f'|_{\b\De_n^\trd}              & := & f^\trd \; .\\
\ea
\]

Ad (2). The composition
\[
\ba{rclcl}
\ulk{\Cl^+(\b\De_n^\#)} & \lra     & \ulk{\Cl^+(\b\De_n^{\trud})} & \lra     & \ulk{\Cl(\dDe_n)}     \\
X                       & \lramaps & X|_{\b\De_n^{\trud}}         & \lramaps & X|_{\dDe_n}           \\
\ea
\]
is an equivalence by Proposition \ref{PropR5}. Therefore, the functor induced by restriction from $\b\De_n^\#$ to $\b\De_n^{\trud}$ 
is faithful. By (1), it is full and dense, and so it is an equivalence. Therefore, 
also the functor induced by restriction from $\b\De_n^{\trud}$ to $\dDe_n$ is an equivalence.\qed

\subsubsection{Cutting off the last object}

\bq
Putting $n = 2$, the equivalence given in Lemma \ref{LemIT1}, composed with the equivalence in Proposition \ref{PropR5}, can be used to retrieve 
{\sc Heller}'s original isomorphism, called $\delr(\Delta)$ in \mb{\bfcite{He68}{p.\ 53}}.
\eq

In this section, we suppose that $n\ge 2$. Consider the functor
\[
\barcl
\Cl(\dDe_n) & \lraa{K}     & \h\Cl(\dDe_{n-1}) \\
X           & \lramaps     & XK \; := (X\ind_0)|_{\dDe_{n-1}} \; ,\\
\ea
\]
where $0$ denotes the morphism $0\lraa{0} X_n\,$; cf.\ \S\ref{SubsecPPP}.

Explicitly, we have $(XK)_i := \Kern(X_i\lraa{x} X_n)$ for $i\in [1,n-1]$, taken in $\h\Cl$, equipped with the induced morphisms 
$(XK)_i\lra (XK)_j$ for $i,\, j\,\in\, [1,n-1]$ with $i \le j$, fitting into a pullback $((XK)_i,(XK)_j,X_i,X_j)$.

Let $\Cl^{\spl}(\dDe_n)\tm\Cl(\dDe_n)$ be the full subcategory defined by 
\[
\Ob\,\Cl^{\spl}(\dDe_n) \; :=\; \{ X\in\Ob\,\Cl(\dDe_n)\; :\; \mb{$(X_i\lra X_j)$ is split in $\h\Cl$ for all $i,\, j\,\in\, [1,n]$ with $i\le j$} \} \; ,
\]
and let $\ulk{\Cl(\dDe_n)} := \Cl(\dDe_n)/\Cl^{\spl}(\dDe_n)$ and $\ulk{\h\Cl(\dDe_{n-1})} := \h\Cl(\dDe_{n-1})/\Cl^{\spl}(\dDe_{n-1})$.

For $Y\in\Ob\h\Cl(\dDe_n)$ and $i\in [1,n]$, we let $YR_i := Y\ind^{Y_i\io}\in\Ob\h\Cl(\dDe_n)$.

We have a {\it resolution map}
\[
\barcl
\Ob\h\Cl(\dDe_n) & \lraa{R'_n} & \Ob\,\Cl(\dDe_n) \\
Y                & \lramaps    & YR'_n \; :=\; YR_0 \cdots R_n \; . \\
\ea
\]
If $Y\in\Ob\h\Cl(\dDe_n)$ consists of monomorphisms, then so does $YR'_n$, whence
$YR'_n\in\Ob\,\Cl^{\spl}(\dDe_n)$.

Given a morphism $Y\lra Y'$ in $\h\Cl(\dDe_n)$ with $Y'$ having bijective entries, this morphism factors over $Y\lramono Y R'_n$ by injectivity of the entries of $Y'$ and by the universal
property of the pointwise pushout.

\begin{Lemma}
\label{LemIT0_5} 
The functor $K$ is dense.
\end{Lemma}

{\it Proof.} Suppose given $X\in\Ob\h\Cl(\dDe_{n-1})$. Let $X'\in\Ob\h\Cl(\dDe_n)$ be defined by $X'|_{\dDe_{n-1}} := X$ and $X'_n := 0$. Then $X'R'_n\in\Ob\,\Cl(\dDe_n)$ has $(X'R'_n) K \iso X$.
\qed

\begin{Lemma}
\label{LemIT0_7} 
The functor $K$ is full.
\end{Lemma}

{\it Proof.} Suppose given $X,\, Y\,\in\,\Ob\,\Cl(\dDe_n)$ and a morphism $XK\lraa{f} YK$. We claim that \mb{there} exists a morphism $X\lraa{\w f} Y$ such that $\w fK = f$. 
We construct its components $\w f_\ell$ by induction on $\ell$. For $\ell = 1$, we obtain a morphism $X_1\lraa{\w f_1} Y_1\ru{6}$ such that $((XK)_1,(YK)_1,X_1,Y_1)$ commutes, by injectivity of $Y_1$ in 
$\h\Cl$. For $\ell \ge 2$, we obtain a morphism $X_\ell\lraa{\w f_\ell} Y_\ell\ru{6}$ such that  
$((XK)_\ell,(YK)_\ell,X_\ell,Y_\ell)$ and $(X_{\ell-1},Y_{\ell-1},X_\ell,Y_\ell)$ commute, by the fact that 
$((XK)_{\ell-1},(XK)_\ell,X_{\ell-1},X_\ell)$ is a weak square and by injectivity of $Y_\ell$.\qed

\begin{Proposition}
\label{LemIT1}
The functor $K$ induces an equivalence
\[
\barcl
\ulk{\Cl(\dDe_n)} & \lraisoa{\ulk{K}} & \ulk{\h\Cl(\dDe_{n-1})} \\
X                 & \lramaps          & XK \; .\\
\ea
\]
\end{Proposition}

{\it Proof.} Let $\w\Cl\tm\h\Cl$ denote the full subcategory of bijective objects in $\h\Cl$. Every object in $\w\Cl$ is a direct summand of an 
object in $\Cl$. Let $\w\Cl^{\spl}(\dDe_{n-1})\tm\w\Cl(\dDe_{n-1})$ be the full subcategory defined by 
\[
\Ob\w\Cl^{\spl}(\dDe_{n-1}) \; :=\; \{ X\in\w\Cl(\dDe_{n-1})\; :\; \mb{$(X_i\lra X_j)$ is split for all $i,\, j\,\in\, [1,n-1]$ with $i\le j$} \} \; ,
\]
Let $Y$ be an object of $\w\Cl^{\spl}(\dDe_{n-1})$. Then $Y R'_{n-1}$ is an object of 
$\Cl^{\spl}(\dDe_{n-1})$ that has $Y$ as a direct summand since the identity on $Y$ factors over $Y\lramono YR_{n-1}$.

Therefore, any morphism that factors over an object of $\w\Cl^{\spl}(\dDe_{n-1})$ already factors over an object of $\Cl^{\spl}(\dDe_{n-1})$.
We infer that
\[
\ulk{\h\Cl(\dDe_{n-1})} \= \h\Cl(\dDe_{n-1})/\w\Cl^{\spl}(\dDe_{n-1}).
\]

Suppose given $X\in\Ob\,\Cl(\dDe_n)$. Denote $X' := XK\in\Ob\h\Cl(\dDe_{n-1})$.

We {\it claim} that if $X\in\Ob\,\Cl^{\spl}(\dDe_n)$, then $X'\in\Ob\w\Cl^{\spl}(\dDe_{n-1})$. First of all, $X'_i$ is bijective for $i\in [1,n-1]$, since the 
image of $X_i\lra X_n$ is bijective, and since $X'_i$ is the kernel of this morphism. Now suppose given $i,\, j\,\in [1,n-1]$ with $i < j$.
In $\h\Cl$, we let $B$ be the image of $X_i\lra X_j$ and form a pullback $(B',X'_j,B,X_j)$. Then there is an induced morphism $X'_i\lra B'$ turning 
$(X'_i,B',X_i,B)$ into a commutative quadrangle, which is a pullback by composition to a pullback $(X'_i,X'_j,X_i,X_j)$.
We insert the common kernel $Z$ of $X_i\lra X_j$ and $X'_i\lra X'_j$.
\[
\xymatrix{
Z \ar@{ >->}~+{|*\dir{*}}[r]    & X_i \ar@{ >->}~+{|*\dir{|}}[r]  & B \ar@{ >->}~+{|*\dir{*}}[r]      & X_j \\
Z \ar@{=}[u]\ar~+{|*\dir{*}}[r] & X'_i \ar~+{|*\dir{|}}[r]\ar~+{|*\dir{*}}[u]\ar@{}[ur]|(0.3){\hookpb} 
                                                                  & B'\ar~+{|*\dir{*}}[r]\ar~+{|*\dir{*}}[u]\ar@{}[ur]|(0.3){\hookpb}  
                                                                                                      & X'_j\ar~+{|*\dir{*}}[u] \\
}
\]
Hence $Z\lramono X'_i$ is split monomorphic, and therefore $X'\lraepi B'$ is split epimorphic. Thus $B'$ is bijective, and so finally 
$B'\lramono X'_j$ is split monomorphic. This proves the {\it claim.}

We {\it claim} that if $X'\in\Ob\w\Cl^{\spl}(\dDe_{n-1})$, then $X\in\Ob\,\Cl^{\spl}(\dDe_n)$. Suppose given $i,\, j\,\in [1,n-1]$ with $i < j$. We
have to show that $X_i\lra X_j$ is split in $\h\Cl$. In $\h\Cl$, we insert the image $B$ of $X_i\lra X_j$ and form a pullback $(B',X'_j,B,X_j)$. Since
$(X'_i,B',X_i,B)$ is a square, and since $X'_i$ is bijective, its diagonal sequence is split short exact. Hence $B$ is bijective
as a direct summand of $X_i\ds B'$, which proves the {\it claim.}

Invoking Lemma \ref{Lem5} to prove the equivalence, it remains to show that given $X\lraa{f} Y$ in $\Cl(\dDe_n)$ such that $fK = 0$, there exists an object $V$ in $\Cl^{\spl}(\dDe_n)$
such that there exists a factorisation
\[
(X\lraa{f} Y) \= (X\lra V\lra Y)\; .
\]
Denote by $XK'\in\Ob\,\Cl(\dDe_n)$ the object that restricts to $XK$ on $\dDe_{n-1}$ and that has $(XK')_n := 0$.
Let $U$ be the cokernel of $XK'\lramono X$ and consider the following diagram.
\[
\xymatrix{
XK' \ar~+{|*\dir{*}}[r]\ar[d]^0 & X \ar~+{|*\dir{|}}[r]\ar[d]^f & U \ar~+{|<>(0.4)*\dir{*}}[r]\ar@{.>}[dl] & UR_n\ar@{.>}[dll] \\
YK' \ar~+{|*\dir{*}}[r]         & Y                             &                                   &                   \\
}
\]
The morphism $U\lra Y$ is induced by the universal property of the cokernel. Its factorisation $(U\lra Y) = (U\lramono UR_n\lra Y)$ exists since
$Y$ consists of bijective objects.

Since the morphism $XK'\lramono X$ consists of pullbacks, its cokernel $U$ consists of monomorphisms. Hence so does $V := UR_n$, which is therefore
in $\Ob\,\Cl^{\spl}(\dDe_n)$, as required.\qed

\subsubsection{Not quite an equivalence}
\label{SecNotQuite}

Let $\Cl^{+,\,\per}(\b\De_n^\#)$ be the subcategory of $\Cl^+(\b\De_n^\#)$ that consists of morphisms $X\lraa{f} Y$ for which
\[
([X^{+1}]\lraa{[f^{+1}]} [Y^{+1}]) \= ([X]^{+1}\lraa{[f]^{+1}} [Y]^{+1})\; ;
\]
which is in general not a full subcategory. The objects $\Cl^{+,\,\per}(\b\De_n^\#)$ are called {\it periodic} $n$-pretriangles, the morphisms are called {\it periodic} morphisms
of periodic $n$-pretriangles. Let $\Cl^{+,\,\spl,\,\per}(\b\De_n^\#) := \Cl^{+,\,\per}(\b\De_n^\#)\cap\Cl^{+,\,\spl}(\b\De_n^\#)$. 

\bq
 For instance, if $(\Cl,\TTT,\tht)$ is a Heller triangulated category, then $\Cl^{+,\,\tht = 1}(\b\De_n^\#)\tm \Cl^{+,\,\per}(\b\De_n^\#)$ is a full subcategory.
\eq

For $Y\in\Ob\,\Cl(\dDe_n)$, we define $YS\in\Ob\,\Cl^+(\b\De_n^\#)$ by 
\[
\barcl
YS|_{\b\De_n^\tru} & := & YS^\tru      \\
YS|_{\b\De_n^\trd} & := & YS^\trd \; , \\
\ea
\]
and similarly on morphisms. If $Y\in\Ob\,\Cl^{\spl}(\dDe_n)$, then $YS\in\Ob\,\Cl^{+,\,\spl}(\b\De_n^\#)$ by Lemma \ref{LemR4}, since $YS|_{\dDe_n} = Y\in\Ob\,\Cl^{\spl}(\dDe_n)$.

To any $X\in\Ob\,\Cl^+(\b\De_n^\#)$ for which $X_{\be/\al}$ is zero for all but finitely many $\be/\al\in \b\De_n^\#\,$, we can assign its {\it periodification}
\[
\b X \; :=\; \Ds_{i\in\Z} [X^{+i}]^{-i}\;\in\; \Ob\,\Cl^{+,\,\per}(\b\De_n^\#)\; ,
\]
and similarly for morphisms between such objects. 

If $X\in\Ob\,\Cl^{+,\,\spl}(\b\De_n^\#)$, then $\b X\in\Ob\,\Cl^{+,\,\spl,\,\per}(\b\De_n^\#)$.

We have the restriction functor
\[
\barcl
\Cl^{+,\,\spl,\,\per}(\b\De_n^\#) & \lrafl{35}{(-)|_{\dDe_n}} & \Cl^{\spl}(\dDe_n) \\
X                                 & \lramapsfl{0}{}           & X|_{\dDe_n} \\
\ea
\] 
which is not faithful in general, as the case $n = 2$ shows. In the inverse direction, we dispose of the functor
\[
\barcl
\Cl^{+,\,\spl,\,\per}(\b\De_n^\#) & \llaa{\b S} & \Cl^{\spl}(\dDe_n) \\
\ol{YS} \; =:\; Y\b S             & \llamaps    & Y\; . \\
\ea
\] 

\begin{Lemma}
\label{LemNQE}
For $X\in\Ob\,\Cl^{+,\,\spl,\,\per}(\b\De_n^\#)$, we have $X\iso X|_{\dDe_n}\b S$. 
\end{Lemma}

\bq
Note that we do {\sf not} claim that $1\iso (-)|_{\dDe_n}\b S$ as endofunctors of $\Cl^{+,\,\spl,\,\per}(\b\De_n^\#)$.
\eq

{\it Proof.} We have a short exact sequence
\[
X|_{\dDe_n}S|_{\b\De_n^{\trud}}\;\lra\; X|_{\b\De_n^{\trud}} \;\lra\; [X|_{\dDe_n}S^{+1}]^{-1}|_{\b\De_n^{\trud}}
\]
in $\h\Cl(\b\De_n^{\trud})$, and it suffices to show that it splits. Write $C := X|_{\dDe_n}S|_{\b\De_n^{\trud}}$.

It suffices to show that there exists a retraction to $C \lra X|_{\b\De_n^{\trud}}$, which we will construct by induction. 
Suppose given $0 <\al < \be \le 0^{+1}$. We may assume that after restriction of $C \lra X|_{\b\De_n^{\trud}}$ to $\{\de/\ga\in \b\De_n^{\trud}\; :\; \de/\ga < \be/\al\}\ru{4}$, there 
exists a retraction. Let $\al' := \al - 1$ be the predecessor of $\al$, and let $\be' := \be - 1$ be the predecessor of $\be$, using that $\b\De_n$ is linearly ordered. 
It suffices to show that the morphism from the quadrangle $(C_{\be'/\al'},C_{\be'/\al},C_{\be/\al'},C_{\be/\al})$ to the quadrangle
$(X_{\be'/\al'},X_{\be'/\al},X_{\be/\al'},X_{\be/\al})$ has a retraction. 

Let $(X_{\be'/\al'},X_{\be'/\al},X_{\be/\al'},T)$ be the pushout in $\h\Cl$. The quadrangle $(C_{\be'/\al'},C_{\be'/\al},C_{\be/\al'},C_{\be/\al})$ is a pushout. The induced morphism from 
$(C_{\be'/\al'},C_{\be'/\al},C_{\be/\al'},C_{\be/\al})$ to $(X_{\be'/\al'},X_{\be'/\al},X_{\be/\al'},T)$ has a retraction by functoriality of the pushout. The morphism $T\lra X_{\be/\al}$ 
induced by pushout is a monomorphism in $\h\Cl$, since $(X_{\be'/\al'},X_{\be'/\al},X_{\be/\al'},X_{\be/\al})$ is a weak square. Note that $(C_{\be/\al},T,X_{\be/\al})$ is a commutative
triangle.

The morphism $T\lra C_{\be/\al}$ that is part of the retraction of quadrangles, factors as
\[
(T\lra C_{\be/\al}) \= (T\lramono X_{\be/\al}\lra C_{\be/\al})\; ,
\]
since $C_{\be/\al}$ is injective in $\h\Cl$ as a summand of $C_{\be/0} = X_{\be/0}$. Now $X_{\be/\al}\lra C_{\be/\al}$ completes the three morphisms on the other vertices to a retraction of 
quadrangles from $(X_{\be'/\al'},X_{\be'/\al},X_{\be/\al'},X_{\be/\al})$ to $(C_{\be'/\al'},C_{\be'/\al},C_{\be/\al'},C_{\be/\al})$ as sought.\qed

%% file: h3verdier.tex
\pagebreak

\section{Verification of Verdier's axioms}
\label{SecVerdier}

Let $(\Cl,\TTT,\tht)$ be a Heller triangulated category. 

\subsection{Restriction from $\Cl^{+,\,\tht = 1}(\b\De_n^\#)$ to $\Cl(\dDe_n)$ is dense and full}

Let $n\ge 1$. Let
\[
\barcl
\Cl(\dDe_n^{+1}) & \lraa{S'^\trd} & \h\Cl^+(\b\De_n^{\trd,+1}) \\
U                & \lramaps       & US'^\trd \; :=\; (U^{(+1)} S^\trd)^{(-1)} \\
\ea
\]
be the conjugate by reindexing, i.e.\ a ``shifted version'' of $S^\trd$; cf.\ \S\ref{ParReindex}. Note that
\[
(U S'^\trd)_{\be/\al} \= \Kern(U_{0^{+1}/\al} \lraa{u} U_{0^{+1}/\be})
\]
for $\al,\, \be\,\in\, \b\De_n$ with $0\le\al\le\be\le 0^{+1}$.

\begin{Lemma}
\label{LemV1}
Suppose that idempotents split in $\Cl$. Given $X\in\Ob\,\Cl(\dDe_n)$, there exists an \mb{$n$-triangle} $\w X\in\Ob\,\Cl^{+,\,\tht=1}(\b\De_n^\#)$ that restricts to 
\[
\w X|_{\dot\De_n} \= X\; .
\]
In other words, the restriction functor $\;\Cl^{+,\,\tht = 1}(\b\De_n^\#)\;\lrafl{35}{(-)|_{\dot\De_n}}\;\Cl(\dDe_n)\;$ is strictly dense.
\end{Lemma}

{\it Proof.} Let $Y := XR\in\Ob\,\Cl^+(\b\De_n^\#)$; cf.\ \S\ref{ParRes}. We have an isomorphism \mb{$[Y]^{+1}\;\lraisoa{Y\tht_n}\;[Y^{+1}]$} in $\ulk{\Cl^+(\b\De_n^\#)}\ru{-2.5}$. Let 
$[Y^{+1}]\lraa{\theta}[Y]^{+1}$ be a representative in $\Cl^+(\b\De_n^\#)$ of the inverse isomorphism $(Y\tht_n)^-$ in $\ulk{\Cl^+(\b\De_n^\#)}$.
Consider the morphism (\footnote{We recall the convention that the inverse
of the outer shift applied to a morphism $f$ is written $[f]^{-1}$, whereas $f^-$ denotes the inverse morphism, if existent.})
\[
[Y^{+1}]^{-1}|_{\dDe_n^{+1}} \;\;\mrafl{35}{[\theta]^{-1}|_{\dDe_n^{+1}}}\;\; [Y]^{+1-1}|_{\dDe_n^{+1}} \;\; =\;\; Y|_{\dDe_n^{+1}} \; .
\]
We have an induced pointwise epimorphism
\[
Y|_{\b\De_n^{\trud}}\;\lraepi\;  Y|_{\dDe_n^{+1}} S'^\trd|_{\b\De_n^{\trud}}\; ,
\]
which we may use to form the pullback 
\[
\xymatrix{
Z \ar[d]_f\ar~+{|(0.25)*\dir{|}}[r]\ar@{}[rd]_(0.2){\hookpbo} 
              & [Y^{+1}]^{-1}|_{\dDe_n^{+1}}S'^\trd|_{\b\De_n^{\trud}} \ar[d]^{[\theta]^{-1}|_{\dDe_n^{+1}}S'^\trd|_{\b\De_n^{\trud}}} \\
Y|_{\b\De_n^{\trud}} \ar~+{|(0.38)*\dir{|}}[r]        
              & Y|_{\dDe_n^{+1}}S'^\trd|_{\b\De_n^{\trud}}\\
}
\]
in the abelian category $\h\Cl(\b\De_n^{\trud})$, i.e.\ pointwise. An application of Lemma \ref{Lem3_2}.(2) to the diagonal sequence of this pullback 
shows that $Z\in\Ob\h\Cl^{+,\ast}(\b\De_n^{\trud})\ru{4.5}$. We obtain $Z_{\al/\al} = 0$ for all $0\le\al\le 0^{+1}$; and we obtain $Z_{0^{+1}/0} = 0$. Hence 
we have $Z\in\Ob\h\Cl^+(\b\De_n^{\trud})\ru{4.5}$.

Suppose given $\be/\al\in\b\De_n^{\trud}$. We {\it claim} that $Z_{\be/\al}\lrafl{28}{f_{\be/\al}} Y_{\be/\al}$ represents an isomorphism in $\h\Cl/\Cl$. By Lemma \ref{Lem3_4}, it suffices to show that 
\[
([Y^{+1}]^{-1}|_{\dDe_n^{+1}}S'^\trd)_{\be/\al} 
\;\;\mvlrafl{35}{([\theta]^{-1}|_{\dDe_n^{+1}}S'^\trd)_{\be/\al}}\;\;
(Y|_{\dDe_n^{+1}}S'^\trd)_{\be/\al}
\]
represents an isomorphism in $\h\Cl/\Cl$. Since evaluation at $\be/\al$ induces a functor from $\h\Cl^+(\b\De_n^{\trud})/\Cl^{+,\,\spl}(\b\De_n^{\trud})$ to $\h\Cl/\Cl$, where 
$\Cl^{+,\,\spl}(\b\De_n^{\trud})$ denotes the full subcategory of $\Cl^+(\b\De_n^{\trud})$ consisting of diagrams all of whose morphisms split (in $\h\Cl$ or, equivalently, in $\Cl$), 
it suffices to show that 
\[
[Y^{+1}]^{-1}|_{\dDe_n^{+1}}S'^\trd|_{\b\De_n^{\trud}}
\;\;\mvlrafl{35}{[\theta]^{-1}|_{\dDe_n^{+1}}S'^\trd|_{\b\De_n^{\trud}}}\;\;
Y|_{\dDe_n^{+1}}S'^\trd|_{\b\De_n^{\trud}}
\]
represents an isomorphism in $\h\Cl^+(\b\De_n^{\trud})/\Cl^{+,\,\spl}(\b\De_n^{\trud})$. 

Now $(-)^{(-1)} S'^\trd |_{\b\De_n^{\trud}}$ induces a functor from $\ulk{\Cl(\dDe_n)}$ to $\h\Cl^+(\b\De_n^{\trud})/\Cl^{+,\,\spl}(\b\De_n^{\trud})$, since it maps $\Cl^{\spl}(\dDe_n)$
to $\Cl^{+,\,\spl}(\b\De_n^{\trud})$ by Lemma \ref{LemR4}, using that {\sf idempotents are assumed to split in $\Cl$.} Therefore, it suffices to show that
\[
([Y^{+1}]^{-1}|_{\dDe_n^{+1}})^{(+1)}
\;\;\mrafl{35}{([\theta]^{-1}|_{\dDe_n^{+1}})^{(+1)}}\;\;
(Y|_{\dDe_n^{+1}})^{(+1)}
\]
represents an isomorphism in $\ulk{\Cl(\dDe_n)}$. Since $([-]^{-1}|_{\dDe_n^{+1}})^{(+1)} = (-)|_{\dDe_n}$, this means that it suffices to show that
\[
[Y^{+1}]|_{\dDe_n}
\;\;\lrafl{32}{\theta|_{\dDe_n}}\;\;
[Y]^{+1}|_{\dDe_n}
\]
represents an isomorphism in $\ulk{\Cl(\dDe_n)}$. Since $(-)|_{\dDe_n}$ induces a functor from $\ulk{\Cl^+(\De_n^\#)}$ to $\ulk{\Cl(\dDe_n)}$, it suffices to show that
\[
[Y^{+1}]\;\;\lraa{\theta}\;\;[Y]^{+1}
\]
represents an isomorphism in $\ulk{\Cl^+(\b\De_n^\#)}$. This, however, follows by choice of $\theta$. This proves the {\it claim.}

{\sf Since idempotents are assumed to split in $\Cl$,} we can conclude from the claim that for all $\be/\al$ in $\b\De_n^{\trud}$, the entry $Z_{\be/\al}$ is isomorphic in $\h\Cl$ 
to an object in $\Cl$, hence without loss of generality, the entry $Z_{\be/\al}$ is an object of $\Cl$. So $Z\in\Ob\,\Cl^+(\b\De_n^{\trud})$.

We remark that 
\[
\ba{lclcl}
Z|_{\dDe_n}      & = & Y|_{\dDe_n}                  & = & X \\
Z|_{\dDe_n^{+1}} & = & [Y^{+1}]^{-1}|_{\dDe_n^{+1}} & = & (X^{+1})^{(-1)}\; , \\
\ea
\leqno (\ast)
\]
where $X^{+1}$ arises from $X$ by pointwise application of $(-)^{+1}$. Concerning morphisms, we remark that 
\[
\ba{lcl}
(Z\lraa{f} Y|_{\b\De_n^{\trud}})|_{\dDe_n}      & = & (X\lraa{1_X} X) \\
(Z\lraa{f} Y|_{\b\De_n^{\trud}})|_{\dDe_n^{+1}} & = & ([Y^{+1}]^{-1} \lraa{[\theta]^{-1}} Y)|_{\dDe_n^{+1}} \; .\\
\ea
\leqno (\ast\ast)
\]
In fact, on $\dDe_n$, the right hand side column of our pullback vanishes; and on $\dDe_n^{+1}$, the lower row of our pullback is an identity.

Now, $(\ast)$ allows to define the periodic prolongation $\b Z\in\Ob\,\Cl^+(\b\De_n^\#)$ of $Z\in\Ob\,\Cl^+(\b\De_n^{\trud})$ by $\b Z|_{\b\De_n^{\trud}} := Z$ and by the requirement that 
$[\b Z]^{+1} = [\b Z^{+1}]$. 

We {\it claim} that $\b Z\tht_n = 1_{[\b Z]^{+1}}$ in $\ulk{\Cl^+(\b\De_n^\#)}$. Let $\b Z\lraa{\h f} Y$ be an inverse image of $Z\lraa{f} Y|_{\b\De_n^{\trud}}$ under
$\Cl^+(\b\De_n^\#)\;\lrafl{35}{(-)|_{\b\De_n^{\trud}}}\; \Cl^+(\b\De_n^{\trud})$; cf.\ Lemma \ref{LemR6}.(1). By $(\ast\ast)$, we get
\[
\ba{lcl}
(\b Z\lraa{\h f} Y)|_{\dDe_n}      & = & (X\lraa{1_X} X) \\
(\b Z\lraa{\h f} Y)|_{\dDe_n^{+1}} & = & ([Y^{+1}]^{-1} \lraa{[\theta]^{-1}} Y)|_{\dDe_n^{+1}} \; .\\
\ea
\leqno (\ast\ast')
\]
We consider the commutative quadrangle
\[
\xymatrix{
[\b Z]^{+1}\ar[r]^{\b Z\tht_n}\ar[d]_{[\h f]^{+1}} & [\b Z^{+1}]\ar[d]^{[\h f^{+1}]} \\
[Y]^{+1}\ar[r]^{Y\tht_n}                         & [Y^{+1}] \\
}
\]
in $\ulk{\Cl^+(\b\De_n^\#)}$. We restrict it to $\dDe_n$ to obtain the commutative quadrangle
\[
\xymatrix{
[\b Z]^{+1}|_{\dDe_n}\ar[rr]^{\b Z\tht_n|_{\dDe_n}}\ar[d]_{[\h f]^{+1}|_{\dDe_n}} & & [\b Z^{+1}]|_{\dDe_n}\ar[d]^{[\h f^{+1}]|_{\dDe_n}} \\
[Y]^{+1}|_{\dDe_n}\ar[rr]^{Y\tht_n|_{\dDe_n}}                                     & & [Y^{+1}]|_{\dDe_n} \\
}
\]
in $\ulk{\Cl(\dDe_n)}$, which, using $(\ast\ast')$, can be rewritten as
\[
\xymatrix{
X^{+1}\ar[rr]^{\b Z\tht_n|_{\dDe_n}}\ar[d]_{\theta|_{\dDe_n}}   & & X^{+1}\ar[d]^{1_{X^{+1}}} \\
[Y]^{+1}|_{\dDe_n}\ar[rr]^{Y\tht_n|_{\dDe_n}}                   & & X^{+1}\; , \!\!\!\!\\
}
\]
where we did not distinguish in notation between $\theta|_{\dDe_n}$ and its residue class in $\ulk{\Cl(\dDe_n)}$, etc. 

Since $\theta (Y\tht_n) = 1_{[Y^{+1}]}$ in $\ulk{\Cl^+(\b\De_n^\#)}$, we have $\theta|_{\dDe_n} (Y\tht_n|_{\dDe_n}) = 1_{X^{+1}}$ in $\ulk{\Cl(\dDe_n)}$. Thus the last quadrangle shows that
$\b Z\tht_n|_{\dDe_n} = 1_{X^{+1}} = 1_{[\b Z]^{+1}}|_{\dDe_n}$ in $\ulk{\Cl(\dDe_n)}$ as well. Since $\ulk{\Cl^+(\b\De_n^\#)}\;\lraisofl{35}{(-)|_{\b\De_n^{\trud}}}\; \ulk{\Cl^+(\b\De_n^{\trud})}$ 
is an equivalence, we conclude that $\b Z\tht_n = 1_{[\b Z]^{+1}}$ in $\ulk{\Cl^+(\b\De_n^\#)}$; cf.\ Proposition \ref{PropR5}. This proves the {\it claim;} i.e.\ we have shown that $\b Z$ is an 
$n$-triangle.

Since $\b Z|_{\dDe_n} = X$ by $(\ast)$, this proves the lemma.\qed

\bq
In the proof of Lemma \ref{LemV1}, we needed the assumption that idempotents split in $\Cl$ in the equivalent form that the residue class functor 
$\h\Cl\lra\h\Cl/\Cl$ maps precisely the objects isomorphic to objects of $\Cl$ to zero -- just as {\sc Heller} did at that point.
\eq

\begin{Lemma}
\label{LemV2}
Given $n$-triangles $X$ and $Y$ and a morphism 
\[
X|_{\dot\De_n}\;\lraa{f}\; Y|_{\dot\De_n}
\]
in $\Cl(\dot\De_n)$, there exists a morphism $X\lraa{\w f} Y$ of $n$-triangles such that $\w f|_{\dot\De_n} = f$. In other words,
the restriction functor $\;\Cl^{+,\,\tht = 1}(\b\De_n^\#)\;\lrafl{35}{(-)|_{\dot\De_n}}\;\Cl(\dDe_n)\;\ru{6}$ is full.
\end{Lemma}

{\it Proof.} Since the restriction functor $\Cl^+(\b\De_n^\#)\lrafl{35}{(-)|_{\dDe_n}}\Cl(\dDe_n)$ is full by Proposition \ref{PropFull}, 
we find a morphism $X\lraa{g} Y$ in $\Cl^+(\b\De_n^\#)$ such that $g|_{\dDe_n} = f$. 

Let $\ulk{g}$ denote the residue class of $g$ in $\ulk{\Cl^+(\b\De_n^\#)}$.
Since $\tht_n$ is a transformation, we have $[\ulk{g}]^{+1} (Y\tht_n) = (X\tht_n)[\ulk{g}^{+1}]$. Since $X$ and $Y$ are $n$-triangles, both $X\tht_n$ and $Y\tht_n$ 
are identities, and this equality amounts to $[\ulk{g}]^{+1} = [\ulk{g}^{+1}]$, i.e.\ the difference
$[g]^{+1} - [g^{+1}]$ factors over an object of $\Cl^{+,\,\spl}(\b\De_n^\#)$. Restricting to $\dDe_n$, the difference
\[
([g]^{+1} - [g^{+1}])|_{\dDe_n} \= (g|_{\dDe_n^{+1}})^{(+1)} - f^{+1} 
\]
factors over an object of $\Cl^{\spl}(\dDe_n)$. Therefore, $g|_{\dDe_n^{+1}} - (f^{+1})^{(-1)}$ factors over an object $Z$ of $\Cl^{\spl}(\dDe_n^{+1})$,
say, as 
\[
\left(X|_{\dDe_n^{+1}}\mvlrafl{40}{g|_{\dDe_n^{+1}} \; -\; (f^{+1})^{(-1)}} Y|_{\dDe_n^{+1}}\right)
\;\;\=\;\; \left(X|_{\dDe_n^{+1}}\lraa{a} Z \lraa{b} Y|_{\dDe_n^{+1}}\right)\; .
\]

By periodic continuation, it suffices to find a morphism $X|_{\b\De_n^{\trud}}\lraa{\breve f} Y|_{\b\De_n^{\trud}}$ in 
$\Cl^+(\b\De_n^{\trud})$ such that \mb{$\breve f|_{\dDe_n} = f$} and such that $\breve f|_{\dDe_n^{+1}} = (f^{+1})^{(-1)}$.
I.e.\ we have to find a morphism $X|_{\b\De_n^{\trud}}\lraa{h} Y|_{\b\De_n^{\trud}}$ such that $h|_{\dDe_n} = 0$ and such that \mb{$h|_{\dDe_n^{+1}} = ab$},
for then we may take $\breve f := g|_{\b\De_n^\trud} - h$. 

Note that $(ZS'^\trd|_{\b\De_n^{\trud}})|_{\dDe_n} = 0$.
Note that $ZS'^\trd|_{\b\De_n^{\trud}}$ is in $\Cl^+(\b\De_n^{\trud})$, hence in $\Cl^{+,\,\spl}(\b\De_n^{\trud})$ by Lemma \ref{LemR4}.

Since $S'^\trd|_{\b\De_n^{\trud}}$ is right adjoint to restriction to $\dDe_n^{+1}$, we have a morphism  
\mb{$X|_{\b\De_n^{\trud}} \lraa{a'} ZS'^\trd|_{\b\De_n^{\trud}}$} such that $a'|_{\dDe_n^{+1}} = a$. 

Since $\Cl^+(\b\De_n^{\trud}) \;\lrafl{40}{(-)|_{\dDe_n^{+1}}}\; \Cl(\dDe_n^{+1})$ is full by the dual and shifted assertion of Lemma \ref{LemR6}.(1), there is a 
morphism $ZS'^\trd|_{\b\De_n^{\trud}} \lraa{b'} Y|_{\b\De_n^{\trud}}\ru{5}$ such that $b'|_{\dDe_n^{+1}} = b$. 

We may take $h := a'b'$.\qed

\bq
In Lemmata \ref{LemV1} and \ref{LemV2}, we do not claim the existence of a coretraction from $\Cl(\dDe_n)$ to $\Cl^{+,\,\tht = 1}(\b\De_n^\#)$ to restriction.
The construction made in the proof of Lemma \ref{LemV2} involves e.g.\ a choice of a lift $b'$ of $b$. Cf.\ \bfcite{Ve67}{II.1.2.13}.

The fullness used in the proof of Lemma \ref{LemV2} to lift $b$, can also be used to lift $a$. We have used the direct argument and thus seen that the lift $a'$ of $a$ does not
involve a choice.
\eq

\begin{Remark}
\label{RemV1+2}\rm
Suppose that idempotents split in $\Cl$. By Lemmata \ref{LemV1} and \ref{LemV2}, the restriction functor
\[
\Cl^{+,\,\tht = 1}(\b\De_n^\#)\;\lrafl{35}{(-)|_{\dot\De_n}}\; \Cl(\dot\De_n)
\]
is full and strictly dense. By Proposition \ref{PropR5}, the restriction functor
\[
\ulk{\Cl^+(\b\De_n^\#)}\;\lraisofl{35}{(-)|_{\dot\De_n}}\; \ulk{\Cl(\dot\De_n)}
\]
is an equivalence. Denoting by $\ulk{\Cl^{+,\,\tht = 1}(\b\De_n^\#)}$ the image of $\Cl^{+,\,\tht = 1}(\b\De_n^\#)$ in $\ulk{\Cl^+(\b\De_n^\#)}$, we 
obtain a full and strictly dense functor $\ulk{\Cl^{+,\,\tht = 1}(\b\De_n^\#)}\;\lrafl{35}{(-)|_{\dot\De_n}}\; \ulk{\Cl(\dot\De_n)}\ru{6}$.
Since it factors as a faithful embedding $\ulk{\Cl^{+,\,\tht = 1}(\b\De_n^\#)}\hra \ulk{\Cl^+(\b\De_n^\#)}\ru{4.5}\ru{-2.5}$ followed by an equivalence, 
it is also faithful. We end up with equivalences
\[
\ulk{\Cl^{+,\,\tht = 1}(\b\De_n^\#)}\;\lraisofl{35}{(-)|_{\dot\De_n}}\; \ulk{\Cl(\dot\De_n)}\;\; , \Icm \ulk{\Cl^{+,\,\tht = 1}(\b\De_n^\#)}\;\lraiso\; \ulk{\Cl^+(\b\De_n^\#)}\; .
\]
\end{Remark}

\subsection{An omnibus lemma}

Suppose given $n,\, m\,\ge\, 1$. Concerning the category $\;\Cl^{+,\,\per}(\b\De_n^\#)\;$ of periodic $n$-pretriangles and its full subcategory $\Cl^{+,\,\spl,\,\per}(\b\De_n^\#)$, 
cf.\ \S\ref{SecNotQuite}; concerning the category $\;\Cl^{+,\,\tht = 1}(\b\De_n^\#)\;$ of $n$-triangles, cf.\ Definition \ref{Def4}.(ii). Note that 
\[
\Cl^{+,\,\tht = 1}(\b\De_n^\#)\;\tm\;\Cl^{+,\,\per}(\b\De_n^\#)\;\tm\;\Cl^+(\b\De_n^\#)\; ,
\]
and that the first inclusion is full.

\begin{Lemma}
\label{LemV3}\Absit
\begin{itemize}
\item[{\rm (1)}] Let $X$ be an $n$-triangle, and let $\b\De_n\llaa{p}\b\De_m$ be a morphism of periodic linearly 
ordered sets. Then $X p^\#$, obtained by ``restriction along $p$'', is an $m$-triangle.
\item[{\rm (2)}] Let $X$ be a $(2n+1)$-triangle. Then $X\ffk_n$, obtained by folding, is an $(n+1)$-triangle. 
\item[{\rm (3)}] The category $\;\Cl^{+,\,\tht = 1}(\b\De_n^\#)\;$ of $n$-triangles is a full additive subcategory of the category $\;\Cl^{+,\,\per}(\b\De_n^\#)\;$
of periodic $n$-pretriangles, closed under direct summands.
\item[{\rm (4)}] Suppose given an isomorphism $X\lraa{f} Y$ in $\Cl^{+,\,\per}(\b\De_n^\#)$. 
If $X$ is an $n$-triangle, then $Y$ is an $n$-triangle.
\item[{\rm (5)}] Let $X\lraa{f} Y$ be a morphism in $\Cl^{+,\,\per}(\b\De_n^\#)$ such that $f|_{\dDe_n}$ is an isomorphism. 
Then $f$ is an isomorphism.
\item[{\rm (6)}] Let $X$ and $Y$ be $n$-triangles. Suppose given an isomorphism 
$X|_{\dDe_n}\lraisoa{u} Y|_{\dDe_n}$ in $\Cl(\dDe_n)$. 
Then there exists an isomorphism $X\lraisoa{\w u} Y$ in $\;\Cl^{+,\,\tht = 1}(\b\De_n^\#)\;$ such that $\,\w u|_{\Cl(\dDe_n)} = u\,$.
\item[{\rm (7)}] If $X\in\Ob\,\Cl^{+,\,\spl,\,\per}(\b\De_n^\#)$, then $X$ is an $n$-triangle.
\end{itemize}
\end{Lemma}

\bq
Note that Lemma \ref{LemV3}.(5) applies in particular to $n$-triangles and a morphism of $n$-triangles.
\eq

{\it Proof.} Ad (1). In $\ulk{\Cl^+(\b\De_m^\#)}$, we have
\[
(X \ul{p}^\#)\tht_m \= (X\tht_n)\ul{p}^\# \= (1_{[X]^{+1}})\ul{p}^\# \= 1_{[X \ul{p}^\#]^{+1}}\; .
\]

Ad (2). In $\ulk{\Cl^+(\b\De_{n+1}^\#)}$, we have
\[
(X\ul{\ffk}_n)\tht_{n+1} \= (X\tht_{2n+1})\ul{\ffk}_n \= (1_{[X]^{+1}})\ul{\ffk}_n \= 1_{[X\ul{\ffk}_n]^{+1}}\; .
\]

Ad (3). We have to show that 
\[
X,\, Y\,\in\,\Ob\,\Cl^{+,\,\tht = 1}(\b\De_n^\#) \Icm\equ\Icm X\ds Y\,\in\,\Ob\,\Cl^{+,\,\tht = 1}(\b\De_n^\#) \;\; .
\]
But since $\tht_n$ is a morphism between additive functors, we have $(X\ds Y)\tht_n = 1_{[X\ds Y]^{+1}}$ if and only if $X\tht_n = 1_{[X]^{+1}}$ 
and $Y\tht_n = 1_{[Y]^{+1}}$. In fact, $(X\ds Y)\tht_n$ identifies with $\smatzz{X\tht_n}{0}{0}{Y\tht_n}$.

Ad (4). Since $f|_{\dDe_n}$ is an isomorphism in $\Cl(\dDe_n)$, so is its image in $\ulk{\Cl(\dDe_n)}$. Hence the image of $f$ in $\ulk{\Cl^+(\b\De_n^\#)}\ru{4}$ 
is an isomorphism by Proposition \ref{PropR5}. Consider the commutative quadrangle
\[
\xymatrix{
[X]^{+1}\ar[r]^{[f]^{+1}}_\sim\ar[d]_{X\tht_n}^\wr & [Y]^{+1}\ar[d]^{Y\tht_n}_\wr \\
[X^{+1}]\ar[r]^{[f^{+1}]}_\sim                     & [Y^{+1}] \\
}
\]
in $\ulk{\Cl^+(\b\De_n^\#)}$. Since $[f]^{+1} = [f^{+1}]$ by assumption, we conclude from $X\tht_n = 1_{[X]^{+1}}$ that $Y\tht_n = 1_{[Y]^{+1}}\ru{4}$.

Ad (5). It suffices to show that given $0\le i\le j\le n$, the morphism $f_{j/i}$ is an isomorphism in $\Cl$. In fact, we have a morphism
of exact sequences 
\[
(\; f_{i/0},\; f_{j/0},\; f_{j/i},\; f_{i/0}^{+1},\; f_{j/0}^{+1} \;)
\] 
in $\h\Cl$, whose entries except possibly $f_{j/i}$ are isomorphisms; hence also $f_{j/i}$ is isomorphic.

Ad (6). This follows by Lemma \ref{LemV2} using (5).

Ad (7). We have $[X]^{+1}\iso 0$ in $\ulk{\Cl^+(\b\De_n^\#)}$, whence $X\tht_n = 1_{[X]^{+1}}\in\liu{\ulk{\Cl^+(\b\De_n^\#)}}([X]^{+1},[X]^{+1})$.
\qed

\subsection{Turning $n$-triangles}

Let $n\ge 2$.

\begin{Lemma}
\label{LemV4}
Given an $n$-triangle $X\in\Ob\,\Cl^{+,\,\tht = 1}(\b\De_n^\#)$, we define $Y\in\Ob\,\Cl^{+,\,\per}(\b\De_n^\#)$ by letting
\[
(Y_{j/i} \lraa{y} Y_{j'/i'}) \; :=\;  (X_{i^{+1}/j}\lraa{x} X_{i'^{+1}/j'})
\]
for $0\le i\le j\le n$ and $0\le i'\le j'\le n$ such that $i\le i'$ and $j\le j'$, and by letting
\[
(Y_{n/i} \lraa{y} Y_{0^{+1}/i}) \; :=\; (X_{i^{+1}/n} \lraa{-x} X_{i^{+1}/0^{+1}})
\]
for $0\le i\le n$. Then $[X]^{+1}_- := Y$ is an $n$-triangle.
\end{Lemma}

{\it Proof.} Let 
\[
\barcl
2\b\De_{n-1} & \lraa{h_n} & \b\De_n \\
i^{+j}       & \lramaps   & 
\left\{
\ba{ll}
(i+1)^{+j/2} & \mb{if $j\con_2 0$} \\
0^{+(j+1)/2} & \mb{if $j\con_2 1$} \; , \\
\ea
\right.
\ea
\]
where $i\in [0,n-1]$ and $j\in\Z$. The map $h_n$ is a morphism of periodic posets. We claim that 
\[
Y \= X h_n^\# \ffk_{n-1}\; .
\]
Once this claim is shown, we are done by Lemma \ref{LemV3}.(1,\,2).

Note that $(X h_n^\#)_{l/k} = X_{lh_n/kh_n}$ for $k,\, l\,\in\, 2\b\De_{n-1}$ with $k\le l$. For $0\le i\le n$ and $1\le j\le 0^{+1}$, we obtain
\[
\barcl
(X h_n^\# \ffk_{n-1})_{j/i}
& = & 
\left\{
\ba{ll}
(X h_n^\#)_{(j-1)^{+1}/(j-1)}                                       & \mb{for $i = 0$ and $1\le j\le n$} \\
(X h_n^\#)_{(i-1)^{+2}/(i-1)^{+1}}                                  & \mb{for $1\le i\le n$ and $j = 0^{+1}$} \\
(X h_n^\#)_{(j-1)^{+1}/(i-1)^{+1}}\ds (X h_n^\#)_{(i-1)^{+2}/(j-1)} & \mb{for $1\le i\le j\le n$} \\
0                                                                   & \mb{for $i = 0$ and $j = 0^{+1}$} \\
\ea
\right.\vs\\
& = & 
\left\{
\ba{ll}
X_{0^{+1}/j}           & \mb{for $i = 0$ and $1\le j\le n$} \\
X_{i^{+1}/0^{+1}}      & \mb{for $1\le i\le n$ and $j = 0^{+1}$} \\
X_{i^{+1}/j}           & \mb{for $1\le i\le j\le n$} \\
0                      & \mb{for $i = 0$ and $j = 0^{+1}$}\; , \\
\ea
\right. \\
\ea
\]
and also the morphisms result as claimed.
\qed

\subsection{Application to the axioms of Verdier}

Recall that $(\Cl,\TTT,\tht)$ is a Heller triangulated category.

\begin{Proposition}
\label{PropV5}
Suppose that idempotents split in $\Cl$.
The tuple $(\Cl,\TTT)$, equipped with the set of \mb{$2$-triangles} as the set of distinguished triangles, is a triangulated category in the sense of Verdier \bfcite{Ve63}{Def.\ 1-1}.
\end{Proposition}

{\it Proof.} We number the axioms of Verdier as in loc.\ cit.

Ad (TR 1). Stability under isomorphism of the set of distinguished triangles follows from Lemma \ref{LemV3}.(4).

The possible extension of a morphism to a distinguished triangle follows by Lemma \ref{LemV1}.

The distinguished triangle $(X,X,0)$ on the identity of an object $X$ in $\Cl$ follows by Lemma \ref{LemV3}.(7). 
Alternatively, one can use that each morphism is contained in a distinguished triangle 
and the fact that a distinguished triangle is a long exact sequence in $\h\Cl$. 

Ad (TR 2). Suppose given a distinguished triangle 
\[
X\;\lraa{u}\; Y\;\lraa{v}\; Z\;\lraa{w}\; X^{+1}\; . 
\]
By Lemma~\ref{LemV4}, we obtain the distinguished triangle
\[
X^{+1}\;\lraa{u^{+1}}\; Y^{+1}\;\lraa{v^{+1}}\; Z^{+1}\;\lraa{-w^{+1}}\; X^{+2}\; . 
\]
By Lemma~\ref{LemV3}.(1), applied to the morphism $\b\De_2\lla\b\De_2$ that sends $0$ to $2^{-1}$, $1$ to $0$ and $2$ to $1$, we obtain the distinguished 
triangle
\[
Y\;\lraa{v}\;Z\;\lraa{-w}\; X^{+1}\;\lraa{u^{+1}}\; Y^{+1}\; . 
\]
By Lemma~\ref{LemV3}.(4), we obtain the distinguished triangle
\[
Y\;\lraa{v}\;Z\;\lraa{w}\; X^{+1}\;\lraa{-u^{+1}}\; Y^{+1}\; . 
\]

Ad (TR 3). The possible completion of a morphism in $\Cl(\dDe_2)$ to a morphism of distinguished triangles follows from Lemma \ref{LemV2}.

Ad (TR 4). The octahedral axiom, i.e.\ the compatibility of forming cones with composition of morphisms, follows from Lemma \ref{LemV1}, applied to the case $n = 3$, for 
by Lemma \ref{LemV3}.(6), we may arbitrarily choose completions to distinguished triangles.
\qed

Note that $3$-triangles are particular octahedra, in the language of \bfcite{BBD84}{1.1.6}. Using \mb{$3$-triangles,} we will now verify the axiom proposed in 
\bfcite{BBD84}{1.1.13}.

\begin{Lemma}
\label{LemV6}
Suppose given a $3$-triangle $T$ in $\Ob\,\Cl^{+,\,\tht = 1}((2\b\De_1)^\#)$, depicted as follows.
\[
\xymatrix{
        &                             &                                 &                                      & 0 \\
        &                             &                                 & 0 \ar[r]                             & Z^{+1}\ar[u] \\
        &                             & 0 \ar[r]                        & Z''\ar[r]^{w''}\ar[u]\ar@{}[ur]|+    & Y^{+1}\ar[u]_{v^{+1}} \\
        & 0\ar[r]                     & Y'\ar[r]^{v'}\ar[u]\ar@{}[ur]|+ & Z'\ar[r]^{w'}\ar[u]_{z'}\ar@{}[ur]|+ & X^{+1}\ar[u]_{u^{+1}} \\
0\ar[r] & X\ar[r]^u\ar[u]\ar@{}[ur]|+ & Y\ar[r]^v\ar[u]_y\ar@{}[ur]|+   & Z\ar[r]\ar[u]_z\ar@{}[ur]|+          & 0\ar[u] \\
}
\]
Then
\[
T\ffk_1 \= ( \;Y\;\lrafl{30}{vz}\; Z'\;\lrafl{30}{\smatez{z'}{w'}}\; Z''\ds X^{+1}\;\lrafl{40}{\smatze{w''}{-u^{+1}}}\; Y^{+1}\; )
\]
and
\[
Ts^\#\ffk_1 \= (Z'^{-1}\;\lrafl{30}{w'^{-1}u} \; Y \;\lrafl{30}{\smatez{y}{v}}\; Y'\ds Z\;\lrafl{40}{\rsmatze{v'}{-z}}\; Z'\; )
\]
are distinguished triangles, where $2\b\De_1\llaa{s}2\b\De_1$ is the morphism of periodic posets determined by $0s = 1^{-1}$, $1s = 0$, $0^{+1}s = 1$ and 
$1^{+1}s = 0^{+1}$.
\end{Lemma}

{\it Proof.} This follows by Lemma \ref{LemV3}.(1,\,2).\qed

\subsection{$n$-triangles and strictly exact functors}

Let $\Cl\lraa{F}\Cl'$ be a strictly exact functor between Heller triangulated categories $(\Cl,\TTT,\tht)$ and $(\Cl',\TTT',\tht')$. Let $n\ge 0$.

\begin{Lemma}
\label{LemExFunTri}
Given an $n$-triangle $X\in\Ob\,\Cl^{+,\,\tht=1}(\b\De_n^\#)$, the diagram $X (F^+(\b\De_n^\#))$, obtained by pointwise application of $F$
to $X$, is an $n$-triangle, i.e.\ an object of $\Cl'^{+,\,\tht'=1}(\b\De_n^\#)\ru{5}$.
\end{Lemma}

{\it Proof.} Using $F\TTT' = \TTT F$ as well as $[X]^{+1} = [X^{+1}]$, we obtain
\[
\ba{rclcl}
[X F^+(\b\De_n^\#)]^{+1} 
& = & X (F^+(\b\De_n^\#))(\Cl^+(\TTT_{\! n})) 
& = & X (F^+(\TTT_{\! n})) \\
& = & X (\Cl^+(\TTT_{\! n})) (F^+(\b\De_n^\#)) 
& = & [X]^{+1} (F^+(\b\De_n^\#)) \\ 
& = & [X^{+1}] (F^+(\b\De_n^\#)) 
& = & X (\TTT^+(\b\De_n^\#)) (F^+(\b\De_n^\#))\\ 
& = & X ((\TTT F)^+(\b\De_n^\#)) 
& = & X ((F\TTT')^+(\b\De_n^\#)) \\ 
& = & X (F^+(\b\De_n^\#))(\TTT'^+(\b\De_n^\#)) 
& = & [(X (F^+(\b\De_n^\#)))^{+1}] \; . \\ 
\ea
\]
Moreover, 
\[
X (F^+(\b\De_n^\#)) \tht'_n 
\= X \tht_n (F^+(\b\De_n^\#)) 
\= 1_{[X]^{+1}} (F^+(\b\De_n^\#)) 
\= 1_{[X F^+(\b\De_n^\#)]^{+1}}\; .  
\]

\qed

\subsection{A remark on spectral sequences}

\bq
 {\sc Verdier} calls certain pretriangles {\it objets spectraux} (spectral objects); cf.\ \bfcite{Ve67}{Sec.\ II.4}. We shall explain the connection to spectral sequences in our language.
\eq

Consider the linearly ordered set $\Z_\infty := \{-\infty\}\disj\Z\disj\{+\infty\}$. Let $\b\Z_\infty^{\#\#}$ be the subposet of $\b\Z_\infty^\#(\De_1)$ 
consisting of those $\de/\be\bby\ga/\al$ for which
\[
\de^{-1}\;\le\;\al\;\le\;\be\;\le\;\ga\;\le\;\de\;\le\;\al^{+1}\; , 
\]
where $\al,\,\be,\,\ga,\,\de\,\in\,\b\Z_\infty$. A {\it spectral object,} in a slightly different sense from \bfcite{Ve67}{II.4.1.2}, 
is an object of $\Cl^+(\b\Z_\infty^\#)$. The {\it spectral sequence functor}
\[
\barcl
\Cl^+(\b\Z_\infty^\#) & \lraa{\EE} & \h\Cl(\b\Z_\infty^{\#\#}) \\
X                     & \lramaps   & X\EE \; ,\\
\ea
\]
is defined by
\[
X\EE(\de/\be\bby\ga/\al) \; := \; \Img(X_{\ga/\al}\lra X_{\de/\be})
\]
for $\de/\be\bby\ga/\al\in\b\Z_\infty^{\#\#}$, equipped with the induced morphisms.

\begin{Lemma}
\label{LemUsual}
Given $\al,\,\be,\,\ga,\,\de,\,\eps\,\in\,\b\Z_\infty$ such that 
\[
\eps^{-1}\;\le\;\al\;\le\;\be\;\le\;\ga\;\le\;\de\;\le\;\eps\;\le\;\al^{+1}\; , 
\]
and given $X\in\Ob\,\Cl^+(\b\Z_\infty^\#)$, the morphisms appearing in $X\EE$ form a short exact sequence
\[
X\EE(\eps/\be\bby\ga/\al) \;\lramono\; X\EE(\eps/\be\bby\de/\al) \;\lraepi\; X\EE(\eps/\ga\bby\de/\al) \; .
\]
\end{Lemma}

{\it Proof.} This follows by Lemma \ref{Lem3_2_5}, applied to the diagram 
\[
(\, X_{\ga/\al},\, X_{\de/\al},\, X_{\eps/\al},\,X_{\ga/\be},\, X_{\de/\be},\, X_{\eps/\be},\, \ub{X_{\ga/\ga}}_{= 0},\, X_{\de/\ga},\, X_{\eps/\ga}\,)\;.
\]

\vspace*{-5mm}
\qed

\bq
 Note that we may apply a shift $\be/\al\lramaps\al^{+1}/\be$ to the indices, i.e.\ an outer shift to $X$, before applying Lemma \ref{LemUsual},
 to get another short exact sequence. 

 The usual exact sequences of spectral sequence terms can be derived from Lemma \ref{LemUsual}. Cf.\ \mb{\bfcite{Ve67}{II.4.2.6},} \bfcite{De94}{App.}.
\eq

%% file: h3frobhot.tex
\section{The stable category of a Frobenius category is \mb{Heller triangulated}}

Let $\Fl = (\Fl,\TTT,\III,\io,\PPP,\pi)$ be a functorial Frobenius category; cf.\ Definition \ref{DefFunFrob}.(3). Let $\Bl\tm\Fl$ denote the full subcategory of objects in the image of $\III$, 
coinciding with the full subcategory of the objects in the image of $\PPP$; then $\Bl$ is a sufficiently big full subcategory of bijectives in $\Fl$.

\bq
We shall prove in Theorem \ref{Th1Epi4} below that the classical stable category $\ulFl$ carries a Heller triangulation.
\eq

\subsection{Definition of $\Fl^\Box(\b\De_n^\#)$, modelling \ \ul{$\mb{\ul{$\Fl\!$}}\,^+(\b\De_n^\#)\!$}}

\bq
We shall model, in the sense of Proposition \ref{Prop1Epi3} below, the category $\ulFl^+(\b\De_n^\#)$ by a category $\Fl^\Box(\b\De_n^\#)$. Morally, we 
represent weak squares ($+$) in $\ulFl$ by pure squares ($\Box$) in $\Fl$. To do so, we have to represent the zeroes on the boundaries by bijective objects. 
\eq

Let $n\ge 0$. Concerning the notion of a pure square, see \S\ref{SubsubEL}. Let $\Fl^\Box(\b\De_n^\#) \tm \Fl(\b\De_n^\#)$ be the full subcategory defined by
\[
\Ob\Fl^\Box(\b\De_n^\#) \; :=\; \left\{ \;\; X\in \Ob\Fl(\b\De_n^\#)\;\;\; :\; 
\mb{
\begin{tabular}{rl}
1) & $X_{\al/\al}$ and $X_{\al^{+1}/\al}$ are in $\Ob\Bl$ for all $\al\in \b\De_n$ \\
2) & For all $\de^{-1}\le\al\le\be\le\ga\le\de\le\al^{+1}$ in $\b\De_n\,$, \\
   & the quadrangle \\
   & $\xymatrix{
      X_{\ga/\be}\ar[r]^x                       & X_{\de/\be} \\
      X_{\ga/\al}\ar[r]_x\ar[u]^x\ar@{}[ur]|{\Box} & X_{\de/\al}\ar[u]_x \\      
      }$ \\
   & is a pure square. \\
\end{tabular}
}\right\} \;\; .
\]
Given $n,\,m\,\ge\, 0$, a morphism $\b\De_n\llaa{p}\b\De_m$ induces a morphism $\Fl^\Box(p^\#)$, usually, and by abuse of notation, denoted by $p^\#$.

Given an exact functor $\Fl\lraa{F}\w\Fl$ between functorial Frobenius categories that sends bijectives to bijectives, we obtain an induced functor 
$\Fl^\Box(\b\De_n^\#)\mra{F^\Box(\b\De_n^\#)}\w\Fl^\Box(\b\De_n^\#)\ru{4.5}$ by pointwise application of $F$.

Denote by 
\[
\barcl
\Fl^\Box(\b\De_n^\#) & \lraa{M}   & \ulk{\ulFl^+(\b\De_n^\#)} \\
\Fl^\Box(\b\De_n^\#) & \lraa{M'}  & \ulFl^+(\b\De_n^\#) \\
\ulFl^+(\b\De_n^\#)  & \lraa{M''} & \ulk{\ulFl^+(\b\De_n^\#)} \\
\Fl                  & \lraa{N}  & \ulFl \\
\ulFl(\dDe_n)        & \lraa{N'} & \ulk{\ulFl(\dDe_n)} \\
\ea
\]
the respective residue class functors, welldefined by Lemma \ref{LemExToWEx}. In particular, $M = M'M''$.

\subsection{Folding for $\Fl^\Box(\b\De_n^\#)$}

\bq
We model, in the sense of Remark \ref{RemFold}, the folding operation $\ffk$ introduced in \S\ref{SecFold}. 
\eq

Suppose given $n\ge 0$. Let the periodic functor
\[
\barcl
\Fl^\Box((2\b \De_n)^\#) & \lraa{\w\ffk_n} & \Fl^\Box(\ol{\rh\disj \De_n}^{\,\#}) \\
X                        & \lramaps        & X\w\ffk_n \\
\ea
\]
be determined by the following data. Writing $Y := X\w\ffk_n$, we let
\[
\ba{lcl}
(Y_{\al/\rh} \lraa{y} Y_{\be/\rh}) 
& := & \left(X_{\al^{+1}/\al} \;\lraa{x}\; X_{\be^{+1}/\be}\right) \vspace*{2mm}\\
Y_{\rh/\rh}
& := & 0 \\
Y_{\rh^{+1}/\rh}
& := & 0\vspace*{2mm} \\
(Y_{\be/\rh} \lraa{y} Y_{\be/\al}) 
& := & \left(X_{\be^{+1}/\be} \;\lraa{\smatez{x}{x}}\; X_{\be^{+1}/\al^{+1}}\ds X_{\al^{+2}/\be}\right) \vspace*{4mm}\\ 
(Y_{\be/\al} \lraa{y} Y_{\de/\ga}) 
& := & \left(X_{\be^{+1}/\al^{+1}}\ds X_{\al^{+2}/\be} \;\lrafl{40}{\smatzz{x}{0}{0}{x}}\; X_{\de^{+1}/\ga^{+1}}\ds X_{\ga^{+2}/\de}\right) 
                                                                                                        \vspace*{4mm}\\ 
(Y_{\de/\ga} \lraa{y} Y_{\rh^{+1}/\ga}) 
& := & \left(X_{\de^{+1}/\ga^{+1}}\ds X_{\ga^{+2}/\de} \;\lrafl{40}{\rsmatze{x}{-x}}\; X_{\ga^{+2}/\ga^{+1}}\right) \\
\ea
\]
for $\al,\,\be,\,\ga,\,\de\,\in\, \De_n$ with $\al\le\be$, with $\ga\le\de$ and with $\be/\al \le \de/\ga$. 
The remaining morphisms are given by composition. 

We claim that $X\w\ffk_n$ is an object of $\Fl^\Box(\ol{\rh\disj \De_n}^{\,\#})$. 

In fact, by Lemma \ref{Lem2outof3}, we are reduced to considering the quadrangles of $Y$ inside $\De_n^{\trud}$, i.e.\ the quadrangles 
\begin{itemize}
\item[(i)] on $(\ga/\al,\,\de/\al,\,\ga/\be,\,\de/\be)$ for $\al,\,\be,\,\ga,\,\de\,\in\, \De_n$
with $\al\le\be\le\ga\le\de$;
\item[(ii)] on $(\ga/\rh,\,\de/\rh,\,\ga/\be,\,\de/\be)$ for $\be,\,\ga,\,\de\,\in\, \De_n$ with $\be\le\ga\le\de$; 
\item[(iii)] on $(\ga/\al,\,\rh^{+1}/\al,\,\ga/\be,\,\rh^{+1}/\be)$ for $\al,\,\be,\,\ga\,\in\, \De_n$ with $\al\le\be\le\ga$; 
\item[(iv)] and on $(\be/\rh,\,\rh^{+1}/\rh,\,\be/\al,\,\rh^{+1}/\al)$ for $\al,\,\be\,\in\, \De_n$ with $\al\le\be$.
\end{itemize}

Another application of loc.\ cit.\ reduces case (i) to case (ii) (or (iii)). Still another application of loc.\ cit.\ reduces the cases (ii) and (iii) to case (iv).
Now the quadrangle in case (iv) is in fact a pure square, as follows from $X\in\Ob\Fl^\Box(\b\De_n^\#)$ and the definition of a pure square via its pure short exact diagonal sequence.

The construction of $Y$ is functorial in $X$.

\begin{Remark}
\label{RemFold}\rm
We have $\w\ffk_n M' = M'\ffk_n$, and thus $\w\ffk_n M = M\ul{\ffk}_n$ for $n\ge 0$.
\[
\xymatrix{
\Fl^\Box((2\b\De_n)^\#) \ar[r]^{\w\ffk_n}\ar[d]_{M'} & \Fl^\Box(\ol{\rh\disj\De_n}^\#)\ar[d]^{M'} \\
\ulFl^+((2\b\De_n)^\#) \ar[r]^{\ffk_n}\ar[d]_{M''}         & \ulFl^+(\ol{\rh\disj\De_n}^\#)\ar[d]^{M''} \\
\ulk{\ulFl^+((2\b\De_n)^\#)} \ar[r]^{\ul{\ffk}_n}          & \ulk{\ulFl^+(\ol{\rh\disj\De_n}^\#)}          \\
}
\]
\end{Remark}

\bq
\begin{Example}
\label{ExFold}\rm
Let $n = 2$. Note that \mb{$2\b\De_2\iso\b\De_5$}. Let $X\in\Ob\Fl^\Box((2\b\De_2)^\#)$, depicted as follows.
\[
\hspace*{-20mm}
\xymatrix@C3mm{
& & & & & & & & & & \\
&  &  &  & X_{0^{+1}/0^{+1}}\ar[r]^x  & X_{1^{+1}/0^{+1}}\ar[u]\ar[r]^x\ar@{}[ur]|{\Box} & X_{2^{+1}/0^{+1}}\ar[u]\ar[r]^x\ar@{}[ur]|{\Box} & 
  X_{0^{+2}/0^{+1}}\ar[u]\ar[r]^x\ar@{}[ur]|{\Box} & X_{1^{+2}/0^{+1}}\ar[u]\ar[r]^x\ar@{}[ur]|{\Box} 
  & X_{2^{+2}/0^{+1}}\ar[u]\ar[r]^x\ar@{}[ur]|{\Box} & X_{0^{+3}/0^{+1}}\ar[u] \\
&  &  & X_{2/2}\ar[r]^x & X_{0^{+1}/2}\ar[u]^x\ar[r]^x\ar@{}[ur]|{\Box} & X_{1^{+1}/2}\ar[u]^x\ar[r]^x\ar@{}[ur]|{\Box} & 
  \fbox{$X_{2^{+1}/2}$}\ar[u]^x\ar[r]^x\ar@{}[ur]|{\Box} & X_{0^{+2}/2}\ar[u]^x\ar[r]^x\ar@{}[ur]|{\Box} 
  & X_{1^{+2}/2}\ar[u]^x\ar[r]^x\ar@{}[ur]|{\Box} & X_{2^{+2}/2}\ar[u]^x \\
&  & X_{1/1}\ar[r]^x & X_{2/1}\ar[u]^x\ar[r]^x\ar@{}[ur]|{\Box} & X_{0^{+1}/1}\ar[u]^x\ar[r]^x\ar@{}[ur]|{\Box} &
  \fbox{$X_{1^{+1}/1}$}\ar[u]^x\ar[r]^x\ar@{}[ur]|{\Box} & X_{2^{+1}/1}\ar[u]^x\ar[r]^x\ar@{}[ur]|{\Box} 
  & X_{0^{+2}/1}\ar[u]^x\ar[r]^x\ar@{}[ur]|{\Box} & X_{1^{+2}/1}\ar[u]^x \\
& X_{0/0}\ar[r]^x & X_{1/0}\ar[u]^x\ar[r]^x\ar@{}[ur]|{\Box} & X_{2/0}\ar[u]^x\ar[r]^x\ar@{}[ur]|{\Box} & 
  \fbox{$X_{0^{+1}/0}$}\ar[u]^x\ar[r]^x\ar@{}[ur]|{\Box} & X_{1^{+1}/0}\ar[u]^x\ar[r]^x\ar@{}[ur]|{\Box} 
  & X_{2^{+1}/0}\ar[u]^x\ar[r]^x\ar@{}[ur]|{\Box} & X_{0^{+2}/0}\ar[u]^x \\
X_{2^{-1}/2^{-1}}\ar[r]^x & X_{0/2^{-1}}\ar[u]^x\ar[r]^x\ar@{}[ur]|{\Box} & X_{1/2^{-1}}\ar[u]^x\ar[r]^x\ar@{}[ur]|{\Box} 
  & X_{2/2^{-1}}\ar[u]^x\ar[r]^x\ar@{}[ur]|{\Box} & X_{0^{+1}/2^{-1}}\ar[u]^x\ar[r]^x\ar@{}[ur]|{\Box} 
  & X_{1^{+1}/2^{-1}}\ar[u]^x\ar[r]^x\ar@{}[ur]|{\Box} & X_{2^{+1}/2^{-1}}\ar[u]^x \\
\ar[u]\ar@{}[ur]|{\Box}  & \ar[u]\ar@{}[ur]|{\Box} & \ar[u]\ar@{}[ur]|{\Box} & \ar[u]\ar@{}[ur]|{\Box} & \ar[u]\ar@{}[ur]|{\Box} & \ar[u] \\
}
\]
Note that the objects on the boundary of the diagram,
\[
\ba{c}
\dots\;,\;\; X_{2^{-1}/2^{-1}}\;,\;\; X_{0/0}\;,\;\; X_{1/1}\;,\;\; X_{2/2}\;,\;\;  X_{0^{+1}/0^{+1}}\;,\;\;\dots \\
\dots\;,\;\; X_{2^{+1}/2^{-1}}\;,\;\; X_{0^{+2}/0}\;,\;\; X_{1^{+2}/1}\;,\;\; X_{2^{+2}/2}\;,\;\;  X_{0^{+3}/0^{+1}}\;,\;\; \dots \\
\ea
\]
are all supposed to be in $\Ob\Bl$.

Note that $\rh\disj\De_2\iso\De_3$. Folding turns $X$ into $X\w\ffk_2\in\Ob\Fl^\Box(\ol{\rh\disj\De_2}^{\,\#})$, depicted as follows. 
\[
\hspace*{-12mm}
\xymatrix@C5mm{
& & & & & 0\ar[r] & \\
& & & & X_{2^{+1}/2^{+1}}\dk X_{2^{+2}/2}\ar[r]^(0.6){\rsmatze{x}{-x}} &  X_{2^{+2}/2^{+1}}\ar[u]\ar[r]\ar@{}[ur]|{\Box} & \\
& & &  X_{1^{+1}/1^{+1}}\dk X_{1^{+2}/1}\ar[r]^{\smatzz{x}{0}{0}{x}} 
        &  X_{2^{+1}/1^{+1}}\dk X_{1^{+2}/2}\ar[r]^(0.6){\rsmatze{x}{-x}}\ar[u]^{\smatzz{x}{0}{0}{x}}\ar@{}[ur]|{\Box} 
        & X_{1^{+2}/1^{+1}}\ar[u]^x\ar[r]\ar@{}[ur]|{\Box} & \\
& & X_{0^{+1}/0^{+1}}\dk X_{0^{+2}/0}\ar[r]^{\smatzz{x}{0}{0}{x}} 
        & X_{1^{+1}/0^{+1}}\dk X_{0^{+2}/1}\ar[r]^{\smatzz{x}{0}{0}{x}}\ar[u]^{\smatzz{x}{0}{0}{x}}\ar@{}[ur]|{\Box}
        & X_{2^{+1}/0^{+1}}\dk X_{0^{+2}/2}\ar[r]^(0.6){\rsmatze{x}{-x}}\ar[u]^{\smatzz{x}{0}{0}{x}}\ar@{}[ur]|{\Box} 
        & X_{0^{+2}/0^{+1}}\ar[u]^x\ar[r]\ar@{}[ur]|{\Box} &  \\ 
& 0\ar[r] & \fbox{$X_{0^{+1}/0}$}\ar[r]^x\ar[u]^{\smatez{x}{x}}\ar@{}[ur]|{\Box} & \fbox{$X_{1^{+1}/1}$}\ar[r]^x\ar[u]^{\smatez{x}{x}}\ar@{}[ur]|{\Box}
        & \fbox{$X_{2^{+1}/2}$}\ar[r]\ar[u]^{\smatez{x}{x}}\ar@{}[ur]|{\Box} & 0\ar[u] \\
X_{2^{-1}/2^{-1}}\dk X_{2/2^{-2}}\ar[r]^(0.65){\rsmatze{x}{-x}} 
        & X_{2/2^{-1}}\ar[r]^(0.35){\smatez{x}{x}}\ar[u]\ar@{}[ur]|{\Box} 
        & X_{2/0}\dk X_{0^{+1}/2^{-1}}\ar[r]^{\smatzz{x}{0}{0}{x}}\ar[u]^{\rsmatze{x}{-x}}\ar@{}[ur]|{\Box} 
        & X_{2/1}\dk X_{1^{+1}/2^{-1}}\ar[r]^{\smatzz{x}{0}{0}{x}}\ar[u]^{\rsmatze{x}{-x}}\ar@{}[ur]|{\Box} 
        & X_{2/2}\dk X_{2^{+1}/2^{-1}}\ar[u]^{\rsmatze{x}{-x}}\ar[u] \\
\ar[u]\ar@{}[ur]|{\Box}  & \ar[u]\ar@{}[ur]|{\Box} & \ar[u]\ar@{}[ur]|{\Box} & \ar[u] &         \\
}
\]
\end{Example}
\eq

\subsection{Some $1$-epimorphic functors}

Let $n\ge 0$. Concerning $1$-epimorphy, cf.\ \S\ref{Sec1epi}.

\begin{Lemma}
\label{Lem1Epi1}
The restriction functor
\[
\barcl
\Fl^\Box(\b\De_n^\#) & \lrafl{35}{(-)|_{\dDe_n}} & \Fl(\dDe_n) \\
X                    & \lramapsfl{0}{}           & X|_{\dDe_n}
\ea
\]
is $1$-epimorphic.
\end{Lemma}

{\it Proof.} We claim that the functor $(-)|_{\dDe_n}$ satisfies the requirements (i,\,ii) of Corollary \ref{CorE3_5}, which then implies that it is 
$1$-epimorphic.

Suppose given $Y\in\Ob\Fl(\dDe_n)$. We construct an object $\w Y$ of $\Fl^\Box(\b\De_n^\#)$ such that $\w Y|_{\dDe_n} = Y$ by the following procedure. 

Write $\b\De_n^{\tru,\,\cdot} := \b\De_n^\tru\ohne\{0/0\}$ and $\b\De_n^{\trd,\,\cdot} := \b\De_n^\trd\ohne\{0^{+1}/0\}$; cf.\ \S\ref{SecSomPos}.

On $\b\De_n^{\tru,\,\cdot}$, we proceed by induction to construct a diagram for which, moreover, the morphisms $\w Y_{\ga/\al} \lra \w Y_{\ga/\be}$ are
purely monomorphic for all $\al,\,\be,\,\ga\,\in\,\b\De_n$ with $0 \le \al \le \be \le \ga \le \al^{+1}$, and, moreover, for which 
$\w Y_{\al^{+1}/\al} = 0$ for all $0\le\al$.

First of all, let $\w Y|_{\dDe_n} := Y$.

Assume given $\ell\ge 0$ such that $\w Y_{\ka^\tru(\ell')}$, together with all diagram morphisms pointing to position $\ka^\tru(\ell')$,
is already constructed for all $\ell' < \ell$, but such that $\w Y_{\ka^\tru(\ell)}$ is not yet constructed; cf.\ \S\ref{ParKappa}. 

If $\ka^\tru(\ell)$ is of the form $\al/\al$ for some $\al\in \b\De_n$ with $0 < \al$, then choose a pure monomorphism $\w Y_{\al/(\al-1)}\lramono \w Y_{\al/\al}$
into an object $\w Y_{\al/\al}$ of $\Bl$ 

\bq
 We do not necessarily choose $\w Y_{\al/(\al-1)}\,\io$ here.
\eq

If $\ka^\tru(\ell) = \al^{+1}/\al$ for some $\al\in \b\De_n$ with $0 \le \al$, then let $\w Y_{\al^{+1}/\al} := 0$. 

If $\ka^\tru(\ell)$ is of the form $\be/\al$ for some $\al,\,\be\,\in\, \b\De_n$ with $0 < \al < \be < \al^{+1}$, then we let 
\[
(\; \w Y_{(\be - 1)/(\al - 1)}\,,\; \w Y_{(\be - 1)/\al}\,,\; \w Y_{\be/(\al - 1)}\,,\; \w Y_{\be/\al}\;  )
\]
be a pushout. Recall that by induction assumption, $\w Y_{(\be - 1)/(\al - 1)} \lramono \w Y_{(\be - 1)/\al}$ is purely monomorphic.
So $\w Y_{\be/(\al - 1)} \lramono \w Y_{\be/\al}$ is purely monomorphic as well.

On $\b\De_n^{\trd,\,\cdot}$, we proceed dually, and finally glue along $\dDe_n$ to obtain the sought $\w Y$.

Ad (i). The restriction map $\liu{\Fl(\b\De_n^\#)}{(\w Y_1,\w Y_2)} \;\lrafl{35}{(-)|_{\dDe_n}} \liu{\Fl(\dDe_n)}{(Y_1,Y_2)}$ is surjective for 
$Y_1,\,Y_2\,\in\,\Ob\Fl(\dDe_n)$, as we see by induction, using bijectivity to prolong morphisms and universal properties of occurring pushouts and pullbacks.

Ad (ii). Suppose given $X\in\Ob\Fl^\Box(\b\De_n^\#)$. Let $X'' := (X|_{\dDe_n})\w\,\in\Ob\Fl^\Box(\b\De_n^\#)\ru{-2}$. 
Let $X'\in\Ob\Fl^\Box(\b\De_n^\#)$ be defined by $X'|_{\b\De_n^{\tru,\,\cdot}} = X''|_{\b\De_n^{\tru,\,\cdot}}$ and 
by $X'|_{\b\De_n^{\trd,\,\cdot}} = X|_{\b\De_n^{\trd,\,\cdot}}$.

There is a morphism $X'\lra X$ that restricts to the identity of $X|_{\b\De_n^{\trd,\,\cdot}}$ on $\b\De_n^{\trd,\,\cdot}$, and hence to the identity of 
$X|_{\dDe_n}$ on $\dDe_n$.

There is a morphism $X'\lra X''$ that restricts to the identity of $X''|_{\b\De_n^{\tru,\,\cdot}}$ on $\b\De_n^{\tru,\,\cdot}$, and hence to the identity of 
$X|_{\dDe_n}$ on $\dDe_n$.

Now suppose given $X_1,\, X_2\,\in\,\Ob\Fl^\Box(\b\De_n^\#)$ such that $X_1|_{\dDe_n} = X_2|_{\dDe_n}$. Then there is a sequence of morphisms
\[
X_1 \;\;\lla\;\; X'_1\;\;\lra\;\; X''_1 = X''_2 \;\;\lla X_2'\;\;\lra X_2
\]
each of which restricts to the identity of $X_1|_{\dDe_n} = X_2|_{\dDe_n}$ on $\dDe_n$, as required.\qed

\begin{Lemma}
\label{Lem1Epi2}
The functors
\[
\barcl
\Fl(\dDe_n) & \lrafl{25}{N(\dDe_n)}    & \ulFl(\dDe_n) \\
\Fl(\dDe_n) & \lrafl{25}{N(\dDe_n) N'} & \ulk{\ulFl(\dDe_n)} \\
\ea
\]
are $1$-epimorphic.
\end{Lemma}

{\it Proof.} Since $N'$ is full and dense, it is $1$-epimorphic by Corollary \ref{CorE4}. Therefore, it suffices to show that $N(\dDe_n)$ is
$1$-epimorphic. 

We will apply Lemma \ref{LemE3}. Choosing representatives of the occurring morphisms $Z_i\lra Z_{i + 1}$ in an object $Z$ of $\ulFl(\dDe_n)$, where
$i\in [1,n-1]$, we see that $N(\dDe_n)$ is dense.

To fulfill condition (C) of loc.\ cit., we will show that given $X,\, Y\,\in\,\Ob\Fl(\dDe_n)$ and a morphism $(X)(N(\dDe_n))\lraa{f} (Y)(N(\dDe_n))$, there are morphisms
$X'\lraa{h} X$ and $X'\lraa{f'} Y$ in $\Fl(\dDe_n)$ such that $(h) (N(\dDe_n))$ is an isomorphism and such that $(h) (N(\dDe_n)) f = (f')(N(\dDe_n))$.

We proceed by induction on $k\in [1,n]$. Suppose given a diagram
\[
\xymatrix{
\w X_1 \ar[r]^{\w x}\ar[d]_{\w f_1} & \w X_2 \ar[r]^{\w x}\ar[d]_{\w f_2} & \w X_3 \ar[r]^{\w x}\ar[d]_{\w f_3} & \cdots \ar[r]^{\w x} & 
                                                                                    \w X_{n-1}\ar[r]^{\w x}\ar[d]_{\w f_{n-1}} & \w X_n\ar[d]_{\w f_n} \\
Y_1 \ar[r]^y                        & Y_2 \ar[r]^y                        & Y_3 \ar[r]^y                        & \cdots \ar[r]^y      & 
                                                                                      Y_{n-1}\ar[r]^y                          & Y_n \\
}
\]
in $\Fl$ such that $\w x \w f_{i+1} = \w f_i y$ for $i\in [1,k-1]$, and such that $\w x \w f_{i+1} \con_\Bl \w f_i y$ for $i\in [k,n]$, and a morphism $\w X\lraa{\w h} X$ in $\Fl(\dDe_n)$
such that $(\w h) (N(\dDe_n))$ is an isomorphism and such that $(\w h) (N(\dDe_n)) f$ is the morphism in $\ulFl(\dDe_n)$ represented by $\w f$. 

If $k < n$, we shall construct a morphism $\w X'\lraa{\w h'} \w X$ in $\Fl(\dDe_n)$ with each $\w h'_i N$ being an isomorphism, and a diagram 
\[
\xymatrix{
\w X'_1 \ar[r]^{\w x'}\ar[d]_{\w f'_1}& \w X'_2 \ar[r]^{\w x'}\ar[d]_{\w f'_2}& \w X'_3 \ar[r]^{\w x'}\ar[d]_{\w f'_3}& \cdots \ar[r]^{\w x'} & 
                                                                                \w X'_{n-1}\ar[r]^{\w x}\ar[d]_{\w f'_{n-1}}& \w X'_n\ar[d]_{\w f'_n}\\
Y_1 \ar[r]^y                          & Y_2 \ar[r]^y                          & Y_3 \ar[r]^y                          & \cdots \ar[r]^y       & 
                                                                                      Y_{n-1}\ar[r]^y                           & Y_n \\
}
\]
in $\Fl$ such that $\w x' \w f'_{i+1} = \w f'_i y$ for $i\in [1,k]$, such that $\w x' \w f'_{i+1} - \w f'_i y\con_\Bl 0$ for $i\in [k+1,n]$, 
and such that $\w h'_i \w f_i - \w f'_i \con_\Bl 0$ for all $i\in [1,n]$. For then we obtain a commutative diagram in $\ulFl(\dDe_n)$
\[
\xymatrix{
\w X'\ar[d]_{\w h'}\ar[ddr]^{\w f'}              & \\
\w X\ar[d]_{\w h}\ar[dr]_*+<0mm,-1mm>{\scm \w f} & \\
X\ar[r]_f                                        & Y \; ,\!\! \\
}
\]
in which we denoted morphisms by their representatives.

Let $\w X_k\lramonoa{j} B$ be a pure monomorphism to an object $B$ in $\Bl$, and let $\w x \w f_{k+1} - \w f_k y = jg$. Let
\[
\w X' \;\; :=\;\; (\;\w X_1 \lraa{\w x}\cdots\lraa{\w x} \w X_k\lrafl{25}{\smatez{\w x}{j}} \w X_{k+1}\ds B \lrafl{35}{\smatze{\w x}{0}} \w X_{k+2} 
          \lraa{\w x}\cdots\lraa{\w x} \w X_n\;)\; ,
\]
and let 
\[
\w f'_i \; :=\; 
\left\{
\ba{ll}
\w f_i                  & \mb{for $i\in [1,n]\ohne\{ k+1\}$} \vspace*{2mm}\\
\smatze{\w f_{k+1}}{-g} & \mb{for $i = k+1$} \\ 
\ea
\right. , \;\;\;
\w h'_i \; :=\; 
\left\{
\ba{ll}
1_{\w X_i}                         & \mb{for $i\in [1,n]\ohne\{ k+1\}$} \vspace*{2mm}\\
\smatze{1_{\w X_{k+1}}\ru{-3}}{0}  & \mb{for $i = k+1$} \\ 
\ea
\right.  .
\]
\qed

\begin{Proposition}
\label{Prop1Epi3}
The residue class functor $\Fl^\Box(\b\De_n^\#)\lraa{M}\ulk{\ulFl^+(\b\De_n^\#)}$ is $1$-epimorphic.
\end{Proposition}

{\it Proof.} Consider the commutative quadrangle
\[
\xymatrix@C=14mm{
\Fl^\Box(\b\De_n^\#)\ar[r]^M\ar[d]_{(-)|_{\dDe_n}} & \ulk{\ulFl^+(\b\De_n^\#)}\ar[d]^{(-)|_{\dDe_n}} \\
\Fl(\dDe_n)\ar[r]^{N(\dDe_n)N'}                    & \ulk{\ulFl(\dDe_n)}\; .\!\! \\
}
\]
Therein, the functor $\Fl^\Box(\b\De_n^\#)\;\lrafl{35}{(-)|_{\dDe_n}}\;\Fl(\dDe_n)$ is $1$-epimorphic by Lemma \ref{Lem1Epi1}.
The functor $\Fl(\dDe_n)\;\mrafl{25}{N(\dDe_n)\,N'}\;\ulk{\ulFl(\dDe_n)}$ is $1$-epimorphic by Lemma \ref{Lem1Epi2}.
The functor $\ul{\ulFl^+(\b\De_n^\#)}\;\lrafl{35}{(-)|_{\dDe_n}}\;\ul{\ulFl(\dDe_n)}$\vspace{1mm}\ is an equivalence by Proposition \ref{PropR5}. 
Hence by Remark \ref{LemE2}, the functor $\Fl^\Box(\b\De_n^\#)\lraa{M}\ulk{\ulFl^+(\b\De_n^\#)}$ is $1$-epimorphic.\qed

\bq
 We do not claim that the residue class functor $\Fl^\Box(\b\De_n^\#)\;\lrafl{25}{N(\b\De_n^\#)}\;\ulFl^+(\b\De_n^\#)$ is $1$-epimorphic.
\eq

\subsection{Construction of $\tht$}

Let $n\ge 0$. In the notation of Lemma \ref{LemComparison}, we let $C := \b\De_n^\#$; the role of the category called $\El$ there is played by $\Fl$ here;
we let $\Gl := \Fl^\Box(\b\De_n^\#)$; and finally, we let $\Hl := \ulFl^+(\b\De_n^\#)$. Note that $\ulFl^+(\b\De_n^\#)$ is a characteristic 
subcategory of $\ulFl(\b\De_n^\#)$.

The tuples
\[
\barcl
I(n)  & =  & \left(\left(I_{X,\,\be/\al}\right)_{\!\be/\al\in \b\De_n^\#}\right)_{\! X\in\Ob\Fl^\Box(\b\De_n^\#)} \vspace*{2mm}\\
      & := & \left(\left(X_{\be/\al}\mramono{\smatez{x}{x}} X_{\be/\be}\ds X_{\al^{+1}/\al}
         \mraepifl{40}{\rsmatze{x}{-x}} X_{\al^{+1}/\be}\right)_{\!\! \be/\al\in \b\De_n^\#}\right)_{\!\! X\in\Ob\Fl^\Box(\b\De_n^\#)}\vss\\
J(n)  & =  & \left(\left(J_{X,\,\be/\al}\right)_{\!\be/\al\in \b\De_n^\#}\right)_{\! X\in\Ob\Fl^\Box(\b\De_n^\#)}  \vspace*{2mm}\\
      & := & \left(\left(X_{\be/\al}\mramono{X_{\be/\al}\,\io} X_{\be/\al}\III = X_{\be/\al}^{+1}\PPP
          \mraepi{X_{\be/\al}^{+1}\,\pi} X_{\be/\al}^{+1}
          \right)_{\!\! \be/\al\in \b\De_n^\#}\right)_{\!\! X\in\Ob\Fl^\Box(\b\De_n^\#)}\\
\ea
\]
are $\b\De_n^\#$-resolving systems, inducing an isomorphism $\TTT_{\! I(n)}\mraisoa{\al_{I(n),\,J(n)}} \TTT_{\! J(n)}$ by Lemma \ref{LemComparison}.(2). 
Recall that $\ulFl^+(\b\De_n^\#)\lraa{M''} \ulk{\ulFl^+(\b\De_n^\#)}$ denotes the residue class functor. We have
\[
\ba{lclcl}
\TTT_{\! I(n)} M'' & = & M\,\ulk{\ulFl^+(\TTT_{\! n})} & = & M [-]^{+1}      \\
\TTT_{\! J(n)} M'' & = & M\,\ulk{\TTT^+(\b\De_n)}      & = & M [-^{+1}] \; . \\
\ea
\]

Since $M$ is $1$-epimorphic by Proposition \ref{Prop1Epi3}, we obtain
\[
\xymatrix{
\Fl^\Box(\b\De_n^\#)\ar@/^2.5mm/[rrr]^{\TTT_{\! I(n)}}="ti" \ar@/_2.5mm/[rrr]_{\TTT_{\! J(n)}}="tj"\ar[d]_M   & & & \ulFl^+(\b\De_n^\#)\ar[d]^{M''} \\
\ulk{\ulFl^+(\b\De_n^\#)} \ar@/^2.5mm/[rrr]^{[-]^{+1}}="out" \ar@/_2.5mm/[rrr]_{[-^{+1}]}="inn"         & & & \ulk{\ulFl^+(\b\De_n^\#)} \; ,\!\!
\ar@2 "ti"+<0mm,-3mm>;"tj"+<0mm,3mm>^{\;\al_{I(n),\,J(n)}}
\ar@2 "out"+<0mm,-3.1mm>;"inn"+<0mm,3.1mm>^{\;\tht_n}
}
\]
where $\tht_n$ is characterised by this commutative diagram, i.e.\ by
\[
\al_{I(n),\,J(n)} \st M'' \= M\st \tht_n\; .
\]
Since $\al_{I(n),\,J(n)}$ is an isomorphism, so is $\tht_n$. Varying $n$, this defines $\tht = (\tht_n)_{n\ge 0}$.

\begin{Theorem}
\label{Th1Epi4}
The tuple $\,\tht = (\tht_n)_{n\ge 0}\,$ is a Heller triangulation on $\ulFl$.
\end{Theorem}

{\it Proof.} According to Definition~\ref{Def4}.(i), we have to show that the following conditions $(\ast)$ and $(\ast\ast)$ hold. 
\begin{itemize}
\item[$(\ast)$] For $m,\, n\,\ge\, 0$, for a morphism $\b\De_n\llaa{p}\b\De_m$ and for an object $Y\in\Ob\ulk{\ulFl^+(\b\De_n^\#)}$, we have
\[
(Y\ul{p}^\#)\tht_m \= (Y\tht_n)\ul{p}^\#\; .
\]
\item[$(\ast\ast)$] For $n\ge 0$ and for an object $Y\in\Ob\ulk{\ulFl^+((2\b\De_n)^\#)}$, we have
\[
(Y\ul{\ffk}_n)\tht_{n+1} \= (Y\tht_{2n+1})\ul{\ffk}_n\; .
\]
\end{itemize}

Ad $(\ast)$. Recall that $\ul{p}^\#$ stands for $\ulk{\ulFl^+(p^\#)}$, and that $p^\#$ stands for $\Fl^\Box(p^\#)$.
By Proposition \ref{Prop1Epi3}, we may assume $Y = XM$ for some $X\in\Ob\Fl^\Box(\b\De_n^\#)$. Then
\[
\ba{lclcl}
(XM\ul{p}^\#)\tht_m & = & (Xp^\# M)\tht_m                  & = & (Xp^\#\al_{I(m),\,J(m)})M'' \\
(XM\tht_n)\ul{p}^\# & = & (X\al_{I(n),\,J(n)}M'')\ul{p}^\# & = & (X\al_{I(n),\,J(n)}p^\#)M''\; , \\
\ea
\]
so that it suffices to show that $Xp^\#\al_{I(m),\,J(m)} = X\al_{I(n),\,J(n)}p^\#$. 

Starting with $X p^\#\,$, the object $X p^\# \TTT_{\! I(m)}$ is calculated by means of $(I_{Xp^\#,\,\be/\al})_{\be/\al\in\b\De_m^\#}\,$; whereas $X \TTT_{\! I(n)}$ is 
calculated by means of $(I_{X,\,\de/\ga})_{\de/\ga\in\b\De_n^\#}\,$, so that $X \TTT_{\! I(n)} p^\#$ can be regarded as being calculated by means of $(I_{X,\,\be p/\al p})_{\be/\al\in\b\De_m^\#}$. But
\[
I_{X,\,\be p/\al p} 
\= \left(X_{\be p/\al p}\lramonoa{\smatez{x}{x}} X_{\be p/\be p}\ds X_{(\al p)^{+1}/\al p} \lraepifl{40}{\rsmatze{x}{-x}} X_{(\al p)^{+1}/\be p}\right) 
\= I_{X p^\#,\,\be/\al} \; , 
\]
whence $X p^\# \TTT_{\! I(m)} = X\TTT_{\! I(n)}p^\#$. 

Next, starting with $X p^\#\,$, the object $X p^\# \TTT_{\! J(m)}$ is calculated by means of $(J_{Xp^\#,\,\be/\al})_{\be/\al\in\b\De_m^\#}\,$; whereas 
$X \TTT_{\! J(n)} p^\#$ can be regarded as being calculated by means of $(J_{X,\,\be p/\al p})_{\be/\al\in\b\De_m^\#}$. But
\[
J_{X,\,\be p/\al p} 
\=  \left(X_{\be p/\al p}\mramono{X_{\be p/\al p}\io} X_{\be p/\al p}\III = X_{\be p/\al p}^{+1}\PPP
      \mraepi{X_{\be p/\al p}^{+1}\pi} X_{\be p/\al p}^{+1}\right) 
\=  J_{X p^\#,\,\be/\al} \; , 
\]
whence $X p^\# \TTT_{\! J(m)} = X \TTT_{\! J(n)} p^\#$. 

Now 
\[
X p^\# \TTT_{\! I(m)}\;\mraisoa{X p^\# \al_{I(m),\,J(m)}}\; X p^\# \TTT_{\! J(m)}
\]
is induced by $(I_{Xp^\#,\,\be/\al})_{\be/\al\in\b\De_m^\#}$ and by $(J_{Xp^\#,\,\be/\al})_{\be/\al\in\b\De_m^\#}\,$; whereas 
\[
X \TTT_{\! I(n)}p^\#\;\mraisoa{X \al_{I(n),\,J(n)}p^\#}\; X \TTT_{\! J(n)} p^\#
\]
can be regarded as being induced by $(I_{X,\,\be p/\al p})_{\be/\al\in\b\De_m^\#}$\vspace{1mm}\ and by $(J_{X,\,\be p/\al p})_{\be/\al\in\b\De_m^\#}$.
We have just seen, however, that these pairs of tuples coincide. 

Ad $(\ast\ast)$. By Proposition \ref{Prop1Epi3}, we may assume $Y = XM$ for some $X\in\Ob\Fl^\Box((2\b\De_n)^\#)$. 
By Remark \ref{RemFold}, we have
\[
\ba{lclcl}
(XM\ul{\ffk}_n)\tht_{n+1}  & = & (X\w\ffk_n M)\tht_{n+1}                & = & (X\w\ffk_n \al_{I(n+1),\,J(n+1)})M'' \\
(XM\tht_{2n+1})\ul{\ffk}_n & = & (X\al_{I(2n+1),\,J(2n+1)}M'')\ul{\ffk}_n & = & (X\al_{I(2n+1),\,J(2n+1)}\ffk_n)M''\; , \\
\ea
\]
so that it suffices to show that $X\w\ffk_n \al_{I(n+1),\,J(n+1)} = X\al_{I(2n+1),\,J(2n+1)}\ffk_n\,$. 

Starting with $X\w\ffk_n\,$, the object $X\w\ffk_n \TTT_{\! I(n+1)}$ is calculated by means of 
$(I_{X\w\ffk_n,\,\be/\al})_{\be/\al\,\in\,\ol{\rh\disj\De_n}^\#}\,$; whereas $X \TTT_{\! I(2n+1)}\ffk_n$ can be regarded as being calculated by means of 
the tuple of pure short exact sequences consisting of
\[
\left\{
\ba{ll}
0                                                               & \mb{at $(\rh/\rh)^{+z}$, $z\in\Z$} \\
I_{X,\,(\al^{+1}/\al)^{+z}}                                     & \mb{at $(\al/\rh)^{+z}$, $\al\in\De_n\,$, $z\in\Z$} \\
I_{X,\,(\be^{+1}/\al^{+1})^{+z}}\ds I_{X,\,(\al^{+2}/\be)^{+z}} & \mb{at $(\be/\al)^{+z}$, $\al,\,\be\,\in\,\De_n\,$, $\al\le\be$, $z\in\Z$} \; .\\
\ea
\right.
\]
We have $I_{X\w\ffk_n,\,(\rh/\rh)^{+z}} = 0$ for $z\in\Z$. For $\al\in\De_n\,$, we have
\[
I_{X\w\ffk_n,\,\al/\rh} \= \left(X_{\al^{+1}/\al} \lramonofl{30}{\smatez{x}{x}} X_{\al^{+1}/\al^{+1}}\ds X_{\al^{+2}/\al} 
\lraepifl{40}{\rsmatze{x}{-x}} X_{\al^{+2}/\al^{+1}}  \right) \= I_{X,\,\al^{+1}/\al} \; ,\\
\]
and accordingly at $(\al/\rh)^{+z}$ for $z\in\Z$. Moreover, for $\al,\,\be\,\in\,\De_n$ with $\al\le\be$, we have
\[
\ba{l}
I_{X\w\ffk_n,\,\be/\al} \= \vs\\
\hspace*{-5mm}\left(X_{\be^{+1}/\al^{+1}}\ds X_{\al^{+2}/\be} 
\mramonofl{40}{\rsmatzv{x}{\; 0}{\; x}{0}{0}{x}{0}{-x}}
X_{\be^{+1}/\be^{+1}}\ds X_{\be^{+2}/\be}\ds X_{\al^{+3}/\al^{+1}}\ds X_{\al^{+2}/\al^{+2}}
\;\;\lraepifl{70}{\rsmatvz{x}{0}{0}{-x}{-x}{0}{0}{-x}}\;\;
X_{\al^{+3}/\be^{+1}}\ds X_{\be^{+2}/\al^{+2}}\right)
\ea
\]
and
\[
\ba{l}
I_{X,\,\be^{+1}/\al^{+1}}\ds I_{X,\,\al^{+2}/\be} \= \vs\\
\hspace*{-5mm}\left(X_{\be^{+1}/\al^{+1}}\ds X_{\al^{+2}/\be} 
\mramonofl{40}{\rsmatzv{x}{x}{0}{0}{0}{0}{x}{x}}
X_{\be^{+1}/\be^{+1}}\ds X_{\al^{+3}/\al^{+1}}\ds X_{\al^{+2}/\al^{+2}}\ds X_{\be^{+2}/\be}
\;\;\lraepifl{70}{\rsmatvz{x}{0}{-x}{0}{0}{x}{0}{-x}}\;\;
X_{\al^{+3}/\be^{+1}}\ds X_{\be^{+2}/\al^{+2}}\right)
\ea
\]
Accordingly at $(\be/\al)^{+z}$ for $z\in\Z$.

Since there is an isomorphism from $I_{X\w\ffk_n,\,\al/\rh}$ to $I_{X,\,\be^{+1}/\al^{+1}}\ds I_{X,\,\al^{+2}/\be}$ that has identities on the 
first and on the third terms of the short exact sequences, completed by
\[
X_{\be^{+1}/\be^{+1}}\ds X_{\be^{+2}/\be}\ds X_{\al^{+3}/\al^{+1}}\ds X_{\al^{+2}/\al^{+2}}
\mraisofl{60}{
\enger{
\left(
\ba{rrrr}
\scm  1 &\scm\;  0 &\scm  0 &\scm\; 0 \\
\scm  0 &\scm    0 &\scm  0 &\scm   1 \\
\scm  0 &\scm    1 &\scm  0 &\scm   0 \\
\scm  0 &\scm    0 &\scm -1 &\scm   0 \\
\ea
\right)}
}
X_{\be^{+1}/\be^{+1}}\ds X_{\al^{+3}/\al^{+1}}\ds X_{\al^{+2}/\al^{+2}}\ds X_{\be^{+2}/\be} \ru{11}
\]
on the second terms, the characterisation in Lemma \ref{LemComparison}.(1) shows that we end up altogether with 
$X\w\ffk_n \TTT_{\! I(n+1)} = X\TTT_{\! I(2n+1)}\ffk_n$.

Starting with $X\w\ffk_n\,$, the object $X\w\ffk_n \TTT_{\! J(n+1)}$ is calculated by means of 
$(J_{X\w\ffk_n,\,\be/\al})_{\be/\al\,\in\,\ol{\rh\disj\De_n}^\#}\,$; whereas $X \TTT_{\! J(2n+1)}\ffk_n$ can be regarded as being calculated by means of 
the tuple of pure short exact sequences consisting of
\[
\left\{
\ba{ll}
0                                                               & \mb{at $(\rh/\rh)^{+z}$, $z\in\Z$} \\
J_{X,\,(\al^{+1}/\al)^{+z}}                                     & \mb{at $(\al/\rh)^{+z}$, $\al\in\De_n\,$, $z\in\Z$} \\
J_{X,\,(\be^{+1}/\al^{+1})^{+z}}\ds J_{X,\,(\al^{+2}/\be)^{+z}} & \mb{at $(\be/\al)^{+z}$, $\al,\,\be\,\in\,\De_n\,$, $\al\le\be$, $z\in\Z$} \; .\\
\ea
\right.
\]
We have $J_{X\w\ffk_n,\,(\rh/\rh)^{+z}} = 0$ for $z\in\Z$. For $\al\in\De_n\,$, we have
\[
J_{X\w\ffk_n,\,\al/\rh} \= \left(X_{\al^{+1}/\al}\mramono{X_{\al^{+1}/\al}\;\io} X_{\al^{+1}/\al}\III = X_{\al^{+1}/\al}^{+1}\PPP
          \mraepi{X_{\al^{+1}/\al}^{+1}\,\pi} X_{\al^{+1}/\al}^{+1} \right) \= J_{X,\,\al^{+1}/\al} \; ,\\
\]
and accordingly at $(\al/\rh)^{+z}$ for $z\in\Z$. 

Moreover, for $\al,\,\be\,\in\,\De_n$ with $\al\le\be$, we have
\[
\barcl
J_{X\w\ffk_n,\,\be/\al} 
& = & \Bigg(X_{\be^{+1}/\al^{+1}}\ds X_{\al^{+2}/\be} \Icm
      \mramonofl{35}{(X_{\be^{+1}/\al^{+1}}\ds X_{\al^{+2}/\be})\,\io} \vs\\
&   & (X_{\be^{+1}/\al^{+1}}\ds X_{\al^{+2}/\be})\III = (X_{\be^{+1}/\al^{+1}}\ds X_{\al^{+2}/\be})^{+1}\PPP \Icm
      \mraepifl{40}{(X_{\be^{+1}/\al^{+1}}\ds X_{\al^{+2}/\be})^{+1}\,\pi} \vs\\
&   & X_{\al^{+3}/\be^{+1}}\ds X_{\be^{+2}/\al^{+2}}\Bigg) \;\; \= \;\; J_{X,\,\be^{+1}/\al^{+1}}\ds J_{X,\al^{+2}/\be}\; ,\\
\ea
\]
and accordingly at $(\be/\al)^{+z}$ for $z\in\Z$.

Hence altogether, we conclude that $X\w\ffk_n \TTT_{\! J(n+1)} = X\TTT_{\! J(2n+1)}\ffk_n$.

Now
\[
X\w\ffk_n \TTT_{\! I(n+1)} \;\;\;\;\mraisofl{35}{X\w\ffk_n \al_{I(n+1),\,J(n+1)}}\;\;\;\; X\w\ffk_n \TTT_{\! J(n+1)}
\]
is induced by $(I_{X\w\ffk_n,\,\be/\al})_{\be/\al\in\b\De_{n+1}^\#}$ and by $(J_{X\w\ffk_n,\,\be/\al})_{\be/\al\in\b\De_{n+1}^\#}\,$; whereas 
\[
X \TTT_{\! I(2n+1)}\ffk_n \;\;\;\;\mraisofl{35}{X\al_{I(2n+1),\,J(2n+1)}\ffk_n}\;\;\;\; X \TTT_{\! J(2n+1)}\ffk_n 
\]
can be regarded as being induced by the tuple consisting of
\[
\left\{
\ba{ll}
0                                                               & \mb{at $(\rh/\rh)^{+z}$, $z\in\Z$} \\
I_{X,\,(\al^{+1}/\al)^{+z}}                                     & \mb{at $(\al/\rh)^{+z}$, $\al\in\De_n\,$, $z\in\Z$} \\
I_{X,\,(\be^{+1}/\al^{+1})^{+z}}\ds I_{X,\,(\al^{+2}/\be)^{+z}} & \mb{at $(\be/\al)^{+z}$, $\al,\,\be\,\in\,\De_n\,$, $\al\le\be$, $z\in\Z$} \; .\\
\ea
\right.
\]
and by $(J_{X\w\ffk_n,\,\be/\al})_{\be/\al\in\b\De_{n+1}^\#}$.

Since the respective former tuples are isomorphic by a tuple of isomorphisms that has identities on the first and on the third term, 
and since the respective latter tuples are equal, the characterisation in Lemma \ref{LemComparison}.(3) shows that
in fact $X\w\ffk_n \al_{I(n+1),\,J(n+1)} = X\al_{I(2n+1),\,J(2n+1)}\ffk_n.\;$\qed

\begin{Corollary}
\label{Cor1Epi5}
Let $\El$ be a Frobenius category. There exists a Heller triangulation on $(\uulEl,\TTT)$.
\end{Corollary}

Concerning the stable category $\uulEl$, cf.\ Definition \ref{DefHot}.

{\it Proof.} Let $\Bl\tm\El$ be the full subcategory of bijectives. The category $\Bl^\ac$ is functorially Frobenius by Example \ref{ExFunFrob}. Hence $\uulEl = \ulk{\Bl^\ac}$, equipped with 
the complex shift $\TTT$, carries a Heller triangulation by virtue of Theorem \ref{Th1Epi4}.\qed

\subsection{Exact functors induce strictly exact functors}

\begin{Proposition}
\label{LemExFun}
Suppose given an exact functor 
\[
\Fl\;\;\lraa{F}\;\; \w\Fl
\]
between functorial Frobenius categories $\Fl = (\Fl,\TTT,\III,\io,\PPP,\pi)$ and $\w\Fl = (\w\Fl,\w\TTT,\w\III,\w\io,\w\PPP,\w\pi)$ that satisfies 
\[
\ba{lcr}
F\w\TTT   & = & \TTT\! F \\
F\,\w\III & = & \III   F \\
F\w\PPP   & = & \PPP\! F \; .\!\!\! \\
\ea
\]
Then the induced functor
\[
\ulFl\;\lraa{\ulF}\;\ulk{\w\Fl}
\]
is strictly exact with respect to the Heller triangulations introduced in Theorem {\rm\ref{Th1Epi4}.}
\end{Proposition}

{\it Proof.} Condition (1) of Definition \ref{Def4}.(iii) is satisfied. Condition (2) of loc.\ cit.\ holds since each morphism has a weak kernel that
is sent to a weak kernel of its image; and dually. In fact, given a morphism represented by $X\lraa{f} Y$, the residue class of the kernel of 
$X\ds Y\PPP\lraepifl{40}{\smatze{f\ru{-1.5}}{Y\pi}} Y\ru{8}$ in $\Fl$, composed with $X\ds Y\PPP\lraepifl{40}{\smatze{1}{0}} X\ru{7.5}$, is a weak kernel of the residue class of $X\lraa{f} Y$ by 
Lemma \ref{LemExToWEx} and Remark \ref{RemBal}. Since pure short exact sequences and bijectives are preserved by $F$, this weak kernel is preserved by $\ulF$.

Consider condition (3) of loc.\ cit. Let $\tht$ resp.\ $\w\tht$ be the Heller triangulation on $\ulFl$ resp.\ on $\ulk{\w\Fl}$ characterised as in Theorem \ref{Th1Epi4} by
\[
\barcl
\al_{I(n),\,J(n)} M''          & = & M\tht_n \\
\al_{\w I(n),\,\w J(n)} \w M'' & = & \w M\w\tht_n \; ,\\
\ea
\]
where $\w M$, $\w M'$, $\w M''$, $\w I(n)$ resp.\ $\w J(n)$ is defined over 
$\w\Fl$ as $M$, $M'$, $M''$, $I(n)$ resp.\ $J(n)$ is over $\Fl$. To prove (3), i.e.\ to show that for $n\ge 0$ and $Y\in\Ob\ulk{\ulFl^+(\b\De_n^\#)}$, we have
\[
(Y \tht_n) \ulk{\ulF^+(\b\De_n^\#)} \= (Y \ulk{\ulF^+(\b\De_n^\#)}) \w\tht_n\; ,
\]
we may assume by Proposition \ref{Prop1Epi3} that $Y = XM$ for some $X\in\Ob\Fl^\Box(\b\De_n^\#)$. Since 
\[
\ba{rclcl}
(XM \tht_n) \ulk{\ulF^+(\b\De_n^\#)}   & = & (X\al_{I(n),\,J(n)} M'') \ulk{\ulF^+(\b\De_n^\#)} & = & (X\al_{I(n),\,J(n)} \ulF^+(\b\De_n^\#)) \w M'' \vspace*{1mm}\\
(XM \ulk{\ulF^+(\b\De_n^\#)}) \w\tht_n & = & (X F^\Box(\b\De_n^\#) \w M) \w\tht_n 
                                                                                   & = & (X F^\Box(\b\De_n^\#)\al_{\w I(n),\,\w J(n)})\w M''\;,\\
\ea
\]
it suffices to show that $X\al_{I(n),\,J(n)} \ulF^+(\b\De_n^\#) = X F^\Box(\b\De_n^\#)\al_{\w I(n),\,\w J(n)}$.

Starting with $X F^\Box(\b\De_n^\#)\,$, the object $X F^\Box(\b\De_n^\#) \TTT_{\! \w I(n)}$ is calculated by means of \linebreak
$(\w I_{X F^\Box(\b\De_n^\#),\,\be/\al})_{\be/\al\in\b\De_n^\#}\,$; whereas 
$X \TTT_{\!\w I(n)} \ulF^+(\b\De_n^\#)$ can be regarded as being calculated by means of $(I_{X,\,\be/\al} F)_{\be/\al\in\b\De_n^\#}\ru{4}\ru{-0.5}$,
where $I_{X,\,\be/\al} F$ is defined by an application of $F$ to all three terms and both morphisms of the pure short exact sequence 
$I_{X,\,\be/\al}$. Since $F$ is additive, we get 
\[
(\w I_{X F^\Box(\b\De_n^\#),\,\be/\al})_{\be/\al\in\b\De_n^\#} \= (I_{X,\,\be/\al} F)_{\be/\al\in\b\De_n^\#} \; ,
\]
whence $X F^\Box(\b\De_n^\#) \TTT_{\!\w I(n)} = X \TTT_{\! I(n)} \ulF^+(\b\De_n^\#)$.

Starting with $X F^\Box(\b\De_n^\#)\,$, the object $X F^\Box(\b\De_n^\#) \TTT_{\!\w J(n)}$ is calculated by means of \linebreak
$(\w J_{X F^\Box(\b\De_n^\#),\,\be/\al})_{\be/\al\in\b\De_n^\#}\,$; whereas 
$X \TTT_{\! J(n)} \ulF^+(\b\De_n^\#)$ can be regarded as being calculated by means of $(J_{X,\,\be/\al} F)_{\be/\al\in\b\De_n^\#}\ru{4}$, where
$J_{X,\,\be/\al} F$ is obtained by entrywise application of $F$. Since $F$ commutes with $\PPP$ and $\w\PPP$, and with $\III$ and $\w\III$, 
we get
\[
(\w J_{X F^\Box(\b\De_n^\#),\,\be/\al})_{\be/\al\in\b\De_n^\#} \= (J_{X,\,\be/\al} F)_{\be/\al\in\b\De_n^\#} \; ,
\]
whence $X F^\Box(\b\De_n^\#) \TTT_{\!\w J(n)} = X \TTT_{\! J(n)} \ulF^+(\b\De_n^\#)$.

Moreover, since the defining pairs of tuples coincide, we finally get 
$X F^\Box(\b\De_n^\#)\al_{\w I(n),\,\w J(n)} = X\al_{I(n),\,J(n)} \ulF^+(\b\De_n^\#)$.\qed

Suppose given an exact functor 
\[
\El\;\;\lraa{E}\;\; \w\El
\]
between Frobenius categories $\El$ and $\w\El$ that sends all bijective objects in $\El$ to bijective objects in $\w\El$.
Let $\Bl\tm\El$ resp.\ $\w\Bl\tm\w\El$ be the respective subcategories of bijectives. We obtain an induced functor $\Bl^\ac\lraa{E^\ac}\w\Bl^\ac$, inducing in turn a functor 
\[
\uulE := \ulk{E^\ac}\;\; :\;\;\uulEl = \ulk{\Bl^\ac}\;\;\lra\;\;\ulk{\w\Bl^\ac} = \w{\uulEl\ru{3}} 
\]
modulo split acyclic complexes; cf.\ Example \ref{ExFunFrob}.(2). 

\begin{Corollary}
\label{CorExFun}
The induced functor
\[
\uulEl\;\lrafl{30}{\uulE}\;\,\w{\!\uulEl\ru{3}}
\]
is strictly exact with respect to the Heller triangulations on $\uulEl$ and on $\,\w{\!\uulEl\ru{3}}$ introduced in Theorem {\rm\ref{Th1Epi4}} via the functorial Frobenius categories $\Bl^\ac$ and 
$\w\Bl^\ac$.
\end{Corollary}

{\it Proof.} We may apply Proposition \ref{LemExFun} to $(\Fl\lraa{F}\w\Fl) := (\Bl^\ac\lraa{E^\ac}\w\Bl^\ac)$.\qed

%% file: h3qcyc.tex
\section{Some quasicyclic categories}
\label{SecDefInf}

\bq
In the definition of a Heller triangulated category, the categories $\ulk{\Cl^+(\b\Delta_n^\#)}$ occur. Replacing this classical stable category by its stable counterpart, these turn out to be
Heller triangulated themselves. So we can iterate. Cf.\ \mb{\bfcite{Be00}{Prop.\ 8.4}.}
\eq

Let $\Cl$ be a weakly abelian category. Let $n\ge 0$.

\subsection{The category $\Cl^+(\b\De_n^\#)$ is Frobenius}

\subsubsection{The category $\Al^0(\b\De_n^\#)$ is Frobenius}
\label{SecAl0Frob}

\bq
 We proceed in a slightly more general manner than necessary. We generalise the fact that the category of complexes $\Al^0(\b\De_2^\#)$ over an additive category $\Al$ is
 a Frobenius category, to a category $\Al^0(\b\De_n^\#)$ for $n\ge 0$; cf.\ Lemma \ref{LemA0Frob} below. Then we will specialise 
 to our weakly abelian category $\Cl$ and pass to the full subcategory $\Cl^+(\b\De_n^\#) \tm \Cl^0(\b\De_n^\#)$; cf.\ Proposition \ref{PropFC1} below.
\eq

\paragraph{Notation}
\indent

Let $\Al$ be an additive category. Let $\Al^0(\b\De_n^\#)$ be the full subcategory of $\Al(\b\De_n^\#)$ defined by
\[
\Ob\Al^0(\b\De_n^\#) \; :=\; \big\{ X\in\Ob\Al(\b\De_n^\#)\; :\; \mb{$X_{\al/\al} = 0$ and $X_{\al^{+1}/\al} = 0$ for all $\al\in\b\De_n$}\big\}\; .
\]
A sequence $X'\lraa{i} X\lraa{p} X''$ in $\Al^0(\b\De_n^\#)$ is called {\it pointwise split short exact} if the sequence 
$X'_\xi\lraa{i_\xi} X_\xi\lraa{p_\xi} X''_\xi$ is split short exact for all $\xi\in\b\De_n^\#$. A morphism is called {\it pointwise split monomorphic}
(resp.\ {\it epimorphic}\,) if it appears as a kernel (resp.\ cokernel) in a pointwise split short exact sequence. 

The category $\Al^0(\b\De_n^\#)$ carries an {\it outer shift} functor $X\lramaps [X]^{+1}$, where $[X]^{+1}_{\be/\al} := X_{(\be/\al)^{+1}} = X_{\al^{+1}/\be}$ for $\be/\al\in\b\De_n^\#$.

Recall that $\Al$, together with the set of split short exact sequences, is an exact category; cf.\ Example \ref{Lem5_1}.
So the additive category $\Al^0(\b\De_n^\#)$, equipped with the set of pointwise split short exact sequences as pure short exact sequences, is an exact 
category; cf.\ Example \ref{Rem5_2}. 

Given $\be/\al,\, \de/\ga\,\in\,\b\De_n^\#$, we write $\be/\al\lessdot\de/\ga$ if $\al < \ga$ and $\be < \de$.

Given $A\in\Ob\Al$ and $\be/\al\in\b\De_n^\#$, we denote by $A_{]\al,\be]}$ the object in $\Al^0(\b\De_n^\#)$ consisting of identical 
morphisms wherever possible and having
\[
(A_{]\al,\be]})_{\de/\ga} \; :=\; 
\left\{
\ba{ll}
A & \text{if $\al/\be^{-1} \lessdot \de/\ga \le \be/\al$} \\
0 & \text{else} \\
\ea
\right.
\]
for $\de/\ga\in\b\De_n^\#$. Such an object is called an {\it extended interval.} 

\bq
Intuitively, it is a rectangle with upper right corner at $\be/\al$, and as large as possible in $\Al^0(\b\De_n^\#)$.
\eq

Let $\Al^{+,\,\spl}(\b\De_n^\#)$ be the full subcategory of $\Al^0(\b\De_n^\#)$ consisting of objects isomorphic to summands of objects of the form
\[
\Ds_{\be/\al\;\in\;\b\De_n^\#} (A_{\be/\al})_{]\al,\be]} \; ,
\]
where $A_{\be/\al}\in\Ob\Al$ for $\be/\al\in\b\De_n^\#$. This direct sum exists since it is a finite direct sum at each $\de/\ga\in\b\De_n^\#$. Concerning the 
notation $\Al^{+,\,\spl}(\b\De_n^\#)$, cf.\ also Remark \ref{RemA0Frob} below.

\paragraph{The periodic case}
\indent

Let $\Al'$ be an additive category, equipped with a {\it graduation shift}
automorphism $X\lramaps X[+1]$. We write $X\lramaps X[m]$ for its $m$-th iteration, where $m\in\Z$.

By entrywise application, there is also a graduation shift on $\Al'^0(\b\De_n^\#)$, likewise denoted by $X\lramaps X[+1]$.

As in \S\ref{SecNotQuite}, we define the subcategory $\Al'^{0,\per}(\b\De_n^\#) \tm \Al'^0(\b\De_n^\#)$ 
to consist of the morphisms $X\lraa{f} Y$ in $\Al'^0(\b\De_n^\#)$ that satisfy 
\[
(X[+1]\lraa{f[+1]} Y[+1]) \= ([X]^{+1}\lraa{[f]^{+1}} [Y]^{+1})
\]
So the subcategory $\Al'^{0,\per}(\b\De_n^\#) \tm \Al'^0(\b\De_n^\#)$ is not full in general.

Given $A\in\Ob\Al'$ and $0\le i\le j\le n$, we denote by $A_{]i,j]}$ the object in $\Al'^0(\b\De_n^\#)$ consisting only of 
zero and identical morphisms and having
\[
(A_{]i,j]})_{\de/\ga} \; :=\; 
\left\{
\ba{ll}
A[m] & \text{if $(i/j^{-1})^{+m} \lessdot \de/\ga \le (j/i)^{+m}$ for some $m\in\Z$} \\
0    & \text{else} \\
\ea
\right.
\]
for $\de/\ga\in\b\De_n^\#$. 

\bq
Intuitively, it is a rectangle with upper right corner at $j/i$, and as large as possible in $\Al'^{0,\,\per}(\b\De_n^\#)$,
repeated $\Z$-periodically up to according graduation shift. 
\eq

Let $\Al'^{+,\,\spl,\,\per}(\b\De_n^\#)$ be the full subcategory of $\Al'^{0,\,\per}(\b\De_n^\#)$ consisting of objects isomorphic to summands of objects 
of the form
\[
\Ds_{0\le i\le j\le n} (A_{j,i})_{]i,j]} \; ,
\]
where $A_{j,i}\in\Ob\Al'$ for $0\le i\le j\le n$. Such an object is called a {\it periodic extended interval.}

\begin{Lemma}
\label{LemA0FrobPer}
The category $\Al'^{0,\,\per}(\b\De_n^\#)$, equipped with the pointwise split short exact sequences, is a Frobenius category,
having $\Al'^{+,\,\spl,\,\per}(\b\De_n^\#)$ as its subcategory of bijectives.
\end{Lemma}

{\it Proof.} By duality, it suffices to show that the following assertions (1,\,2) hold.
\begin{itemize}
\item[(1)] The object $A_{]i,j]}$ is injective in $\Al'^{0,\,\per}(\b\De_n^\#)$ for any $A\in\Ob\Al'$ and any 
           \mb{$0\le i\le j\le n$.}
\item[(2)] For each object of $\Al'^{0,\,\per}(\b\De_n^\#)$, there exists a pure monomorphism into an object of $\Al'^{+,\,\spl,\,\per}(\b\De_n^\#)$.
\end{itemize}

Ad (1). Note that we have an adjunction isomorphism
\[
\barcl
\liu{\Al'^{0,\,\per}(\b\De_n^\#)}(X,A_{]i,j]}) & \lraiso  & \liu{\Al'}{(X_{j/i}\,,A)} \\
f                                              & \lramaps & f_{j/i} \; ,\\
\ea
\]
where $X\in\Ob\Al'^{0,\,\per}(\b\De_n^\#)$. Suppose given a pure monomorphism $A_{]i,j]}\lramono X$ for some $A\in\Ob\Al'$. Let $(A\lramono X_{j/i}\lraepi A) = 1_A$. Let $X\lra A_{]i,j]}$
correspond to $X_{j/i}\lraepi A$. The composition $(A_{]i,j]}\lramono X\lra A_{]i,j]})$ restricts to $1_A$ at $j/i$, hence equals $1_{A_{]i,j]}}$.

Ad (2). Suppose given $X\in\Ob\Al'^{0,\,\per}(\b\De_n^\#)$. Given $0\le i \le j \le n$, we let 
\[
X\;\lrafl{30}{X s_{j/i}}\; (X_{j/i})_{]i,j]}
\]
be the morphism corresponding to $1_{X_{j/i}}$ by adjunction, which is natural in $X$. Collecting these morphisms yields a morphism
\[
X\lraa{X s} \Ds_{0\le i \le j \le n} (X_{j/i})_{]i,j]}\; ,
\]
which is pointwise split monomorphic since at $j/i$, its component $X_{j/i}\lra X_{j/i}$ is an identity.
\qed

\paragraph{The general case}
\indent

\begin{Lemma}
\label{LemA0Frob}
The category $\Al^0(\b\De_n^\#)$, equipped with the pointwise split short exact sequences, is a Frobenius category,
having $\Al^{+,\,\spl}(\b\De_n^\#)$ as its subcategory of bijectives.
\end{Lemma}

{\it Proof.} To prove that $\Al^0(\b\De_n^\#)$ is a Frobenius category, we more precisely claim that $\Al^{+,\,\spl}(\b\De_n^\#)$ is a sufficiently 
big category of bijective objects in the exact category $\Al^0(\b\De_n^\#)$.

Abbreviate $\Al^\Z := \Al(\dot\Z)$, where $\dot\Z$ denotes the discrete category with $\Ob\dot\Z = \Z$ and only identical morphisms.
The category $\Al^\Z$ carries the graduation shift automorphism 
\[
\barcl
\Al^\Z         & \lraiso  & \Al^\Z \\
(X\lraa{f} Y)  & \lramaps & (X[+1]\lrafl{27}{f[+1]} Y[+1]) \; :=\; (X_{i+1}\lrafl{27}{f_{i+1}} Y_{i+1})_{i\in\Z} \; . \\
\ea
\]

We have an isomorphism of categories
\[
\barcl
\Al^+(\b\De_n^\#)                         & \lraisoa{\Phi}  & (\Al^{\Z})^{+,\,\per}(\b\De_n^\#) \\
X                                         & \lramaps        & \left((X_{(\be/\al)^{+i}})_{i\in\Z}\right)_{\be/\al\in\b\De_n^\#} \\
((Y_{\be/\al})_0)_{\be/\al\in\b\De_n^\#}  & \llamaps        & Y \; .\\
\ea
\]
Both categories are exact when equipped with pointwise split short exact sequences, and $\Phi$ and $\Phi^{-1}$ are exact functors.
We have $\Al^{+,\,\spl}(\b\De_n^\#)\Phi = (\Al^{\Z})^{+,\,\spl,\,\per}(\b\De_n^\#)$. 

Putting $\Al' := \Al^\Z$, the result follows by Lemma \ref{LemA0FrobPer}.\qed

If $\Al = \Cl$ is a weakly abelian category, we have two definitions of $\Cl^{+,\,\spl}(\b\De_n^\#)$.

The first one, given in \S\ref{SecDefUl}, defines this category as a full subcategory of $\Cl^+(\b\De_n^\#)$ containing those diagrams in which
all morphisms are split in $\h\Cl$.

The second one, just given, defines this category as a full subcategory of $\Cl^0(\b\De_n^\#)$ containing, up to isomorphism, summands of direct sums of
extended intervals. 

\begin{Remark}
\label{RemA0Frob}
If $\Al = \Cl$ is a weakly abelian category, then the two aforementioned definitions of $\Cl^{+,\,\spl}(\b\De_n^\#)$ coincide. 
\end{Remark}

{\it Proof.} First, we notice that an extended interval lies in $\Cl^+(\b\De_n^\#)$, and that all its diagram morphisms are split in $\h\Cl$. 

It remains to be shown that an object in $\Cl^+(\b\De_n^\#)$ all of whose diagram morphisms are split in $\h\Cl$, is, up to isomorphism, a summand of a direct sum
of extended intervals.

Passing to $(\Cl^\Z)^{+,\,\per}(\b\De_n^\#)$, we have to show that an object $X\in\Ob(\Cl^\Z)^{+,\,\per}(\b\De_n^\#)$ all of whose diagram morphisms are split in $\h\Cl^\Z$, 
is, up to isomorphism, a summand of a direct sum of periodic extended intervals.

By Lemma \ref{LemSpl}, applied to the abelian Frobenius category $\h\Cl^\Z$, the object $X|_{\dDe_n}\in\Ob\,\Cl^\Z(\dDe_n)$ is isomorphic to a summand of a finite direct 
sum of intervals. Hence, by Lemma \ref{LemNQE}, the object $X$ is isomorphic to a summand of a finite direct sum of images of intervals under $\b S$, i.e.\ of periodic extended
intervals, as was to be shown.
\qed

\subsubsection{The subcategory $\Cl^+(\b\De_n^\#)\tm\Cl^0(\b\De_n^\#)$}
\label{SecCl+Frob}

Recall that $\Cl$ is a weakly abelian category.

\begin{Lemma}
\label{LemFC0_5}
Suppose given a pure short exact sequence
\[
X'\;\lramono\; X\;\lraepi\; X''
\]
in $\Cl^0(\b\De_n^\#)$. If two out of the three objects $X'$, $X$ and $X''$ are in $\Cl^+(\b\De_n^\#)$, so is the third.
\end{Lemma}

{\it Proof.} For an object $X\in\Ob\,\Cl^0(\b\De_n)$ to lie in $\Ob\,\Cl^+(\b\De_n)$, it suffices to know that the complex
\[
X(\al,\be,\ga) \;:=\;  \big(\cdots\;\lra\; X_{\be/\ga^{-1}}\;\lra\; X_{\be/\al}\;\lra\; X_{\ga/\al}\;\lra\; X_{\ga/\be}\;\lra\; X_{\al^{+1}/\be}\;\lra\;\cdots\big)
\]
is acyclic in $\h\Cl$ for all $\al,\,\be,\,\ga\,\in\,\b\De_n$ with $\al\le\be\le\ga\le\al^{+1}$; which is true, as we take from Lemma \ref{Lem2}; cf.\ Remark \ref{RemBal}.

Now the long exact homology sequence, applied in $\h\Cl$ to the short exact sequence $X'(\al,\be,\ga)\lramono X(\al,\be,\ga)\lraepi  X''(\al,\be,\ga)$ of complexes, shows that 
if two of these complexes are acyclic, so is the third.
\qed

\begin{Proposition}
\label{PropFC1}\Absit
\begin{itemize}
\item[{\rm (1)}] The category $\Cl^+(\b\De_n^\#)$, equipped with the pointwise split short exact sequences, is a Frobenius 
category, having $\Cl^{+,\,\spl}(\b\De_n^\#)$ as its subcategory of bijectives.

Hence its stable category $\uul{8}{\Cl^+(\b\De_n^\#)}$\vspace{1mm}\ is equivalent to its classical stable category
$\ulk{\Cl^+(\b\De_n^\#)} = \Cl^+(\b\De_n^\#)/\Cl^{+,\,\spl}(\b\De_n^\#)$\vspace{1mm}. So both $\uul{8}{\Cl^+(\b\De_n^\#)}$ and $\ulk{\Cl^+(\b\De_n^\#)}$ are weakly abelian.

\item[{\rm (2)}] Suppose $\Cl$ to be equipped with an automorphism $X\lramaps X^{+1}$.
The category $\Cl^{+,\,\per}(\b\De_n^\#)$ (cf.\ \S\ref{SecNotQuite}), equipped with the pointwise split short exact sequences, 
is an additively functorial Frobenius category, having 
$\Cl^{+,\,\spl,\,\per}(\b\De_n^\#)$ as its subcategory of bijectives.
\end{itemize}
\end{Proposition}

We remark that $\ulk{\Cl^+(\b\De_n^\#)}$ appears in Definition {\rm\ref{Def4}}.

{\it Proof.} Ad (1). To prove that $\Cl^+(\b\De_n^\#)$ is an exact category, it remains to be shown, in view of Lemma \ref{LemA0Frob} and of \S\ref{SecGQL}, that a pure short exact sequence 
in $\Cl^0(\b\De_n^\#)$ that has the first and the third term in $\Ob\,\Cl^+(\b\De_n^\#)$, has the second term in $\Ob\,\Cl^+(\b\De_n^\#)$, too. This follows by Lemma \ref{LemFC0_5}.

To prove that $\Cl^+(\b\De_n^\#)$ is Frobenius, we may use that $\Cl^0(\b\De_n^\#)$ is Frobenius, with the bijective objects already lying in $\Cl^+(\b\De_n^\#)$, thus being a fortiori
bijective with respect to $\Cl^+(\b\De_n^\#)$. By duality, it remains to be shown that the kernel of a pointwise split epimorphism of a bijective object to a given object in $\Cl^+(\b\De_n^\#)$ 
is again in $\Cl^+(\b\De_n^\#)$, thus showing that this epimorphism is actually pure in $\Cl^+(\b\De_n^\#)$. This follows by Lemma \ref{LemFC0_5}.

Ad (2). In view of Lemma \ref{LemA0FrobPer}, this follows as (1).
\qed

\pagebreak

\bq
 I do not know whether $\Cl^{+,\,\tht=1}(\b\De_n^\#)$ is Frobenius. It seems doubtful, since this question is \mb{reminiscent} of the example of 
 {\sc A.\ Neeman} that shows that the mapping cone of a morphism of distinguished triangles in the sense of Verdier need not be distinguished \bfcite{Ne91}{p.\ 234}.
\eq

\subsubsection{Two examples}

\bq
Suppose $\Cl$ to be equipped with an automorphism $X\lramaps X^{+1}$.

The category $\Cl^{+,\,\per}(\b\De_n^\#)$ being a Frobenius category by Proposition \ref{PropFC1}.(2), its classical stable category
$\ul{\Cl^{+,\,\per}(\b\De_n^\#)}\ru{4}\ru{-1}$ carries a Heller operator, defined on $X\in\Ob\,\Cl^{+,\,\per}(\b\De_n^\#)$ as the kernel of $B\lraepi X$, where 
$B\in\Ob\,\Cl^{+,\,\spl,\,\per}(\b\De_n^\#)$. As examples, we calculate the Heller operator for $n\in\{ 2,3\}$ on periodic $n$-pretriangles.

Suppose $n = 2$. Let $X\in\Ob\,\Cl^{+,\,\per}(\b\De_2^\#)$ be a periodic $2$-pretriangle. We obtain
\begin{center}
\begin{picture}(1500,1150)
\put(   0,   0){$X_{1/0}$}
\put( 100,  10){\vector(1,0){280}}
\put( 220,  25){$\scm x$}
\put( 400,   0){$X_{2/0}$}
\put( 480,  40){\vector(2,1){200}}
\put( 560, 105){$\scm x$}
\put( 700, 150){$X_{2/1}$}
\put( 800, 160){\vector(1,0){280}}
\put( 920, 175){$\scm x$}
\put(1100, 150){$X_{1/0}^{+1}$}
\put(1180, 190){\vector(2,1){200}}
\put(1260, 255){$\scm x$}
\put(1400, 300){$X_{2/0}^{+1}$}

\put(  25, 380){\vector(0,-1){330}}
\put( -45, 200){$\smateckze{x}{1}$}
\put( 425, 380){\vector(0,-1){330}}
\put( 355, 200){$\smateckze{x}{1}$}
\put( 725, 530){\vector(0,-1){330}}
\put( 655, 350){$\smateckze{x}{1}$}
\put(1125, 530){\vector(0,-1){330}}
\put(1055, 350){$\smateckze{x}{1}$}
\put(1425, 680){\vector(0,-1){330}}
\put(1355, 500){$\smateckze{x}{1}$}

\put( -64, 400){$X_{2/1}^{-1}\dk X_{1/0}$}
\put( 140, 410){\vector(1,0){165}}
\put( 180, 435){$\scm\smateckzz{0}{0}{1}{0}$}
\put( 336, 400){$X_{1/0}\dk X_{2/0}$}
\put( 460, 450){\vector(2,1){160}}
\put( 500, 525){$\scm\smateckzz{0}{0}{1}{0}$}
\put( 636, 550){$X_{2/0}\dk X_{2/1}$}
\put( 840, 560){\vector(1,0){165}}
\put( 880, 585){$\scm\smateckzz{0}{0}{1}{0}$}
\put(1036, 550){$X_{2/1}\dk X_{1/0}^{+1}$}
\put(1160, 600){\vector(2,1){160}}
\put(1200, 675){$\scm\smateckzz{0}{0}{1}{0}$}
\put(1336, 700){$X_{1/0}^{+1}\dk X_{2/0}^{+1}$}

\put(  25, 780){\vector(0,-1){330}}
\put( -80, 600){$\smateckez{1}{-x}$}
\put( 425, 780){\vector(0,-1){330}}
\put( 320, 600){$\smateckez{1}{-x}$}
\put( 725, 930){\vector(0,-1){330}}
\put( 620, 750){$\smateckez{1}{-x}$}
\put(1125, 930){\vector(0,-1){330}}
\put(1020, 750){$\smateckez{1}{-x}$}
\put(1425,1080){\vector(0,-1){330}}
\put(1320, 900){$\smateckez{1}{-x}$}

\put(   0, 800){$X_{2/1}^{-1}$}
\put( 100, 810){\vector(1,0){280}}
\put( 220, 825){$\scm -x$}
\put( 400, 800){$X_{1/0}$}
\put( 480, 840){\vector(2,1){200}}
\put( 550, 905){$\scm -x$}
\put( 700, 950){$X_{2/0}$}
\put( 800, 960){\vector(1,0){280}}
\put( 920, 975){$\scm -x$}
\put(1100, 950){$X_{2/1}$}
\put(1180, 990){\vector(2,1){200}}
\put(1250,1055){$\scm -x$}
\put(1400,1100){$X_{1/0}^{+1}$}
\end{picture}
\end{center}
In particular, if $X$ is a $2$-triangle, i.e.\ an object of $\Cl^{+,\,\tht=1}(\b\De_2^\#)$, then this Heller shift of $X$ is also a $2$-triangle; cf.\ Lemma \ref{LemV4}.

\pagebreak

\begin{landscape} 
\hspace*{-10cm} Suppose $n = 3$. Let $X\in\Ob\,\Cl^{+,\,\per}(\b\De_3^\#)$ be a periodic $3$-pretriangle. We obtain
\begin{center}
\begin{picture}(2300,1350)(800,0)
\put(   0,   0){$X_{1/0}$}
\put( 100,  10){\vector(1,0){280}}
\put( 220,  25){$\scm x$}
\put( 400,   0){$X_{2/0}$}
\put( 500,  10){\vector(1,0){280}}
\put( 620,  25){$\scm x$}
\put( 800,   0){$X_{3/0}$}

\put( 480,  40){\vector(2,1){200}}
\put( 560, 105){$\scm x$}
\put( 880,  40){\vector(2,1){200}}
\put( 960, 105){$\scm x$}

\put( 700, 150){$X_{2/1}$}
\put( 800, 160){\line(1,0){20}}
 \put( 840, 160){\vector(1,0){240}}
\put( 920, 175){$\scm x$}
\put(1100, 150){$X_{3/1}$}
\put(1200, 160){\vector(1,0){280}}
\put(1320, 175){$\scm x$}
\put(1500, 150){$X_{1/0}^{+1}$}

\put(1180, 190){\vector(2,1){200}}
\put(1260, 255){$\scm x$}
\put(1580, 190){\vector(2,1){200}}
\put(1660, 255){$\scm x$}

\put(1400, 300){$X_{3/2}$}
\put(1500, 310){\line(1,0){20}}
 \put(1540, 310){\vector(1,0){240}}
\put(1620, 325){$\scm x$}
\put(1800, 300){$X_{2/0}^{+1}$}

\put(1880, 340){\vector(2,1){200}}
\put(1960, 405){$\scm x$}

\put(2100, 450){$X_{3/0}^{+1}$}


\put(  30, 380){\vector(0,-1){330}}
\put( -30, 200){$\smateckde{x}{x}{1}$}
\put( 430, 380){\vector(0,-1){330}}
\put( 370, 200){$\smateckve{x}{x}{x}{1}$}
\put( 830, 380){\vector(0,-1){330}}
\put( 770, 230){$\smateckde{x}{x}{1}$}
\put( 730, 530){\line(0,-1){100}}
 \put( 730, 395){\vector(0,-1){195}}
\put( 670, 300){$\smateckde{x}{x}{1}$}
\put(1130, 530){\vector(0,-1){330}}
\put(1070, 350){$\smateckve{x}{x}{x}{1}$}
\put(1530, 530){\vector(0,-1){330}}
\put(1470, 380){$\smateckde{x}{x}{1}$}
\put(1430, 680){\vector(0,-1){330}}
\put(1370, 470){$\smateckde{x}{x}{1}$}
\put(1830, 680){\vector(0,-1){330}}
\put(1770, 500){$\smateckve{x}{x}{x}{1}$}
\put(2130, 830){\vector(0,-1){330}}
\put(2070, 650){$\smateckde{x}{x}{1}$}


\put( -74, 405){$\scm X_{2/1}^{-1}\dk X_{3/1}^{-1}\dk X_{1/0}$}
\put( 140, 410){\vector(1,0){130}}
\put( 150, 440){$\scm\enger{\left[
                 \ba{rrrr}
                 \sscm 0 &\sscm 0 &\sscm 0 &\sscm 0 \\
                 \sscm 1 &\sscm 0 &\sscm 0 &\sscm 0 \\
                 \sscm 0 &\sscm 0 &\sscm 1 &\sscm 0 \\
                 \ea
                 \right]}$}
\put( 283, 405){$\scm X_{3/1}^{-1}\dk X_{3/2}^{-1}\dk X_{1/0}\dk X_{2/0}$}
\put( 580, 410){\vector(1,0){130}}
\put( 590, 455){$\scm\enger{\left[
                 \ba{rrr}
                 \sscm 0 &\sscm 0 &\sscm 0 \\
                 \sscm 0 &\sscm 0 &\sscm 0 \\
                 \sscm 1 &\sscm 0 &\sscm 0 \\
                 \sscm 0 &\sscm 1 &\sscm 0 \\
                 \ea
                 \right]}$}
\put( 720, 405){$\scm X_{1/0}\dk X_{2/0}\dk X_{3/0}$}

\put( 460, 450){\vector(2,1){160}}
\put( 480, 545){$\scm\enger{\left[
                 \ba{rrr}
                 \sscm 0 &\sscm 0 &\sscm 0 \\
                 \sscm 1 &\sscm 0 &\sscm 0 \\
                 \sscm 0 &\sscm 0 &\sscm 0 \\
                 \sscm 0 &\sscm 1 &\sscm 0 \\
                 \ea
                 \right]}$}
\put( 860, 450){\vector(2,1){160}}
\put( 960, 460){$\scm\enger{\left[
                 \ba{rrrr}
                 \sscm 0 &\sscm 0 &\sscm 0 &\sscm 0 \\
                 \sscm 1 &\sscm 0 &\sscm 0 &\sscm 0 \\
                 \sscm 0 &\sscm 1 &\sscm 0 &\sscm 0 \\
                 \ea
                 \right]}$}

\put( 616, 555){$\scm X_{3/2}^{-1}\dk X_{2/0}\dk X_{2/1}$}
\put( 840, 560){\vector(1,0){130}}
\put( 850, 590){$\scm\enger{\left[
                 \ba{rrrr}
                 \sscm 0 &\sscm 0 &\sscm 0 &\sscm 0 \\
                 \sscm 1 &\sscm 0 &\sscm 0 &\sscm 0 \\
                 \sscm 0 &\sscm 0 &\sscm 1 &\sscm 0 \\
                 \ea
                 \right]}$}
\put( 988, 555){$\scm X_{2/0}\dk X_{3/0}\dk X_{2/1}\dk X_{3/1}$}
\put(1280, 560){\vector(1,0){140}}
\put(1290, 605){$\scm\enger{\left[
                 \ba{rrr}
                 \sscm 0 &\sscm 0 &\sscm 0 \\
                 \sscm 0 &\sscm 0 &\sscm 0 \\
                 \sscm 1 &\sscm 0 &\sscm 0 \\
                 \sscm 0 &\sscm 1 &\sscm 0 \\
                 \ea
                 \right]}$}
\put(1440, 555){$\scm X_{2/1}\dk X_{3/1}\dk X_{1/0}^{+1}$}

\put(1160, 600){\vector(2,1){160}}
\put(1170, 690){$\scm\enger{\left[
                 \ba{rrr}
                 \sscm 0 &\sscm 0 &\sscm 0 \\
                 \sscm 1 &\sscm 0 &\sscm 0 \\
                 \sscm 0 &\sscm 0 &\sscm 0 \\
                 \sscm 0 &\sscm 1 &\sscm 0 \\
                 \ea
                 \right]}$}
\put(1560, 600){\vector(2,1){160}}
\put(1660, 610){$\scm\enger{\left[
                 \ba{rrrr}
                 \sscm 0 &\sscm 0 &\sscm 0 &\sscm 0 \\
                 \sscm 1 &\sscm 0 &\sscm 0 &\sscm 0 \\
                 \sscm 0 &\sscm 0 &\sscm 1 &\sscm 0 \\
                 \ea
                 \right]}$}

\put(1315, 705){$\scm X_{3/0}\dk X_{3/1}\dk X_{3/2}$}
\put(1540, 710){\vector(1,0){130}}
\put(1545, 740){$\scm\enger{\left[
                 \ba{rrrr}
                 \sscm 0 &\sscm 0 &\sscm 0 &\sscm 0 \\
                 \sscm 1 &\sscm 0 &\sscm 0 &\sscm 0 \\
                 \sscm 0 &\sscm 1 &\sscm 0 &\sscm 0 \\
                 \ea
                 \right]}$}
\put(1688, 705){$\scm X_{3/1}\dk X_{3/2}\dk X_{1/0}^{+1}\dk X_{2/0}^{+1}$}

\put(1860, 750){\vector(2,1){160}}
\put(1870, 840){$\scm\enger{\left[
                 \ba{rrrr}
                 \sscm 0 &\sscm 0 &\sscm 0 \\
                 \sscm 0 &\sscm 0 &\sscm 0 \\
                 \sscm 1 &\sscm 0 &\sscm 0 \\
                 \sscm 0 &\sscm 1 &\sscm 0 \\
                 \ea
                 \right]}$}

\put(2025, 855){$\scm X_{1/0}^{+1}\dk X_{2/0}^{+1}\dk X_{3/0}^{+1}$}


\put(  30, 780){\vector(0,-1){330}}
\put( -80, 600){$\scm\enger{\left[
                 \ba{rrr}
                 \sscm 1 &\sscm 0 &\sscm -x \\
                 \sscm 0 &\sscm 1 &\sscm -x \\
                 \ea
                 \right]}$}
\put( 430, 780){\vector(0,-1){330}}
\put( 300, 600){$\scm\enger{\left[
                 \ba{rrrr}
                 \sscm 1 & \sscm 0 &\sscm 0 &\sscm -x \\
                 \sscm 0 & \sscm 1 &\sscm 0 &\sscm -x \\
                 \sscm 0 & \sscm 0 &\sscm 1 &\sscm -x \\
                 \ea
                 \right]}$}
\put( 830, 780){\vector(0,-1){330}}
\put( 835, 680){$\scm\enger{\left[
                 \ba{rrr}
                 \sscm 1 &\sscm 0 &\sscm -x \\
                 \sscm 0 &\sscm 1 &\sscm -x \\
                 \ea
                 \right]}$}
\put( 730, 930){\line(0,-1){110}}
 \put( 730, 800){\vector(0,-1){200}}
\put( 620, 730){$\scm\enger{\left[
                 \ba{rrr}
                 \sscm 1 &\sscm 0 &\sscm -x \\
                 \sscm 0 &\sscm 1 &\sscm -x \\
                 \ea
                 \right]}$}
\put(1130, 930){\vector(0,-1){330}}
\put(1000, 750){$\scm\enger{\left[
                 \ba{rrrr}
                 \sscm 1 & \sscm 0 &\sscm 0 &\sscm -x \\
                 \sscm 0 & \sscm 1 &\sscm 0 &\sscm -x \\
                 \sscm 0 & \sscm 0 &\sscm 1 &\sscm -x \\
                 \ea
                 \right]}$}
\put(1530, 930){\vector(0,-1){330}}
\put(1535, 830){$\scm\enger{\left[
                 \ba{rrr}
                 \sscm 1 &\sscm 0 &\sscm -x \\
                 \sscm 0 &\sscm 1 &\sscm -x \\
                 \ea
                 \right]}$}
\put(1430,1080){\line(0,-1){110}}
 \put(1430, 950){\vector(0,-1){200}}
\put(1320, 870){$\scm\enger{\left[
                 \ba{rrr}
                 \sscm 1 &\sscm 0 &\sscm -x \\
                 \sscm 0 &\sscm 1 &\sscm -x \\
                 \ea
                 \right]}$}
\put(1830,1080){\vector(0,-1){330}}
\put(1700, 900){$\scm\enger{\left[
                 \ba{rrrr}
                 \sscm 1 & \sscm 0 &\sscm 0 &\sscm -x \\
                 \sscm 0 & \sscm 1 &\sscm 0 &\sscm -x \\
                 \sscm 0 & \sscm 0 &\sscm 1 &\sscm -x \\
                 \ea
                 \right]}$}
\put(2130,1230){\vector(0,-1){330}}
\put(2020,1050){$\scm\enger{\left[
                 \ba{rrr}
                 \sscm 1 &\sscm 0 &\sscm -x \\
                 \sscm 0 &\sscm 1 &\sscm -x \\
                 \ea
                 \right]}$}


\put( -37, 805){$\scm X_{2/1}^{-1}\dk X_{3/1}^{-1}$}
\put( 110, 810){\vector(1,0){200}}
\put( 160, 835){$\scm\enger{\left[
                 \ba{rrr}
                 \sscm 0 &\sscm 0 &\sscm -x \\
                 \sscm 1 &\sscm 0 &\sscm -x \\
                 \ea
                 \right]}$}
\put( 330, 805){$\scm X_{3/1}^{-1}\dk X_{3/2}^{-1} \dk X_{1/0}$}
\put( 550, 810){\vector(1,0){200}}
\put( 620, 845){$\scm\enger{\left[
                 \ba{rr}
                 \sscm 0 &\sscm -x \\
                 \sscm 0 &\sscm -x \\
                 \sscm 1 &\sscm -x \\
                 \ea
                 \right]}$}
\put( 762, 805){$\scm X_{1/0}\dk X_{2/0}$}

\put( 520, 850){\vector(2,1){160}}
\put( 490, 905){$\scm\enger{\left[
                 \ba{rr}
                 \sscm 0 &\sscm -x \\
                 \sscm 1 &\sscm -x \\
                 \sscm 0 &\sscm -x \\
                 \ea
                 \right]}$}
\put( 900, 840){\vector(2,1){180}}
\put( 860, 890){$\scm\enger{\left[
                 \ba{rrr}
                 \sscm 0 &\sscm -x &\sscm 0 \\
                 \sscm 1 &\sscm -x &\sscm 0 \\
                 \ea
                 \right]}$}

\put( 662, 955){$\scm X_{3/2}^{-1}\dk X_{2/0}$}
\put( 810, 960){\vector(1,0){200}}
\put( 860, 985){$\scm\enger{\left[
                 \ba{rrr}
                 \sscm 0 &\sscm 0 &\sscm -x \\
                 \sscm 1 &\sscm 0 &\sscm -x \\
                 \ea
                 \right]}$}
\put(1030, 955){$\scm X_{2/0}\dk X_{3/0}\dk X_{2/1}$}
\put(1250, 960){\vector(1,0){200}}
\put(1310, 995){$\scm\enger{\left[
                 \ba{rr}
                 \sscm 0 &\sscm -x \\
                 \sscm 0 &\sscm -x \\
                 \sscm 1 &\sscm -x \\
                 \ea
                 \right]}$}
\put(1462, 955){$\scm X_{2/1}\dk X_{3/1}$}

\put(1200, 990){\vector(2,1){180}}
\put(1180,1050){$\scm\enger{\left[
                 \ba{rr}
                 \sscm 0 &\sscm -x \\
                 \sscm 1 &\sscm -x \\
                 \sscm 0 &\sscm -x \\
                 \ea
                 \right]}$}
\put(1600, 990){\vector(2,1){180}}
\put(1580,1050){$\scm\enger{\left[
                 \ba{rrr}
                 \sscm 0 &\sscm 0 &\sscm -x \\
                 \sscm 1 &\sscm 0 &\sscm -x \\
                 \ea
                 \right]}$}

\put(1362,1105){$\scm X_{3/0}\dk X_{3/1}$}
\put(1510,1110){\vector(1,0){200}}
\put(1550,1135){$\scm\enger{\left[
                 \ba{rrr}
                 \sscm 0 &\sscm -x &\sscm 0 \\
                 \sscm 1 &\sscm -x &\sscm 0 \\
                 \ea
                 \right]}$}
\put(1730,1105){$\scm X_{3/1}\dk X_{3/2}\dk X_{1/0}^{+1}$}

\put(1900,1140){\vector(2,1){180}}
\put(1880,1200){$\scm\enger{\left[
                 \ba{rr}
                 \sscm 0 &\sscm -x \\
                 \sscm 0 &\sscm -x \\
                 \sscm 1 &\sscm -x \\
                 \ea
                 \right]}$}

\put(2062,1255){$\scm X_{1/0}^{+1}\dk X_{2/0}^{+1}$}
\end{picture}
\end{center}
\hspace*{-10cm} If $X$ is a $3$-triangle, i.e.\ an object $\Cl^{+,\,\tht=1}(\b\De_3^\#)$, I do not know whether this Heller shift of $X$ is again a $3$-triangle.
\end{landscape}
\eq

\subsection{A quasicyclic category}
\label{SecQC}

The category of {\it quasicyclic categories} is defined to be the category of contravariant functors from $\b\Deltab^{\!\circ}\ru{4.5}$ to the ($1$-)category $(\Cat)$ of 
categories. Recall that we have a functor $\Deltab\lra\b\Deltab$, $\De_n\lramaps\b\De_n$ that allows to restrict a quasicyclic category to its
{\it underlying simplicial category.} 

Given a category $\Ul$, we denote by $\text{Iso\,\,}\Ul\tm\Ul$ its subcategory consisting of isomorphisms. Given a functor $\Ul\lraa{U}\Ul'$, we denote by 
$\text{Iso\,} F:\text{Iso\,\,}\Ul\lra \text{Iso\,\,}\Ul'$ the induced functor.

We define
\[
\barcl
\b\Deltab^{\!\circ}          & \lrafl{30}{\qcyc_\bt\Cl} & (\Cat) \\
(\b\De_n\llaa{p}\b\De_m)     & \lramaps                 & \Big(\text{Iso\,}\Cl^+(\b\De_n^\#)\mrafl{25}{\text{Iso\,}\Cl^+(p^\#)} \text{Iso\,}\Cl^+(\b\De_m^\#)\Big)\; . \\
\ea
\]
More intuitively written, $\qcyc_\bt\Cl := \text{Iso\,}\Cl^+(\b\De_\bt^\#)$.
Note that $\qcyc_0\Cl$ consists only of zero-objects. 

A strictly exact functor $\Cl\lraa{F}\Cl'$ induces a functor $\Cl^+(\b\De_n^\#)\mrafl{30}{F^+(\b\De_n^\#)}\Cl'^+(\b\De_n^\#)$ for $n\ge 0$, and thus a morphism
\[
\qcyc_\bt\Cl\;\mrafl{30}{\qcyc_\bt F}\; \qcyc_\bt\Cl'
\]
of quasicyclic categories.

As variants, we mention 
\[
\barcl
\qcyc^{\per}_\bt\Cl     & := & \text{Iso\,}\Cl^{+,\,\per}(\b\De_\bt^\#) \\
\qcyc^{\tht = 1}_\bt\Cl & := & \text{Iso\,}\Cl^{+,\,\tht = 1}(\b\De_\bt^\#)\; . \\
\ea
\]

\subsection{A biquasicyclic category}
\label{SecBQC}

\bq
 As an attempt in the direction described in \bfcite{Wa85}{p.\ 330}, we define a first step of an ``iteration'' of the construction $\Cl\ramaps\qcyc_\bt\Cl$.
\eq

By Proposition \ref{PropFC1}, we may form the category $\uul{8}{\Cl^+(\b\De_n^\#)}^+(\b\De_m^\#)$. Note that a morphism $\b\De_m\llaa{f}\b\De_{m'}$ of 
periodic linearly ordered sets induces a functor $\uul{8}{\Cl^+(\b\De_n^\#)}^+(f^\#)$ in the second variable.

By Lemma \ref{LemExFun}, a morphism $\b\De_n\llaa{g}\b\De_{n'}$ of periodic
linearly ordered sets induces a strictly exact functor $\uul{8}{\Cl^+(g^\#)}$, and so a functor $\uul{8}{\Cl^+(g^\#)}^+(\b\De_m^\#)$ in the first variable for $m\ge 0$.

The functors induced by $f$ and by $g$ commute.

We may define 
\[
\qcyc_{\bt\bt}\Cl \; :=\; \text{Iso}\left(\uul{8}{\Cl^+(\b\De_\bt^\#)}^+(\b\De_\bt^\#)\right)\; , 
\]
which yields a {\it biquasicyclic category,} i.e.\ a functor from $\b\Deltab^{\!\circ}\ti\b\Deltab^{\!\circ}$ to the ($1$-)category $(\Cat)$ of categories.

By Lemma \ref{LemExFun}, a strictly exact functor $\Cl\lraa{F}\Cl'$ induces a functor 
$\uul{8}{\Cl^+(\b\De_n^\#)}\mrafl{40}{\uul{6}{F^+(\b\De_n^\#)}}\uul{8}{\Cl'^+(\b\De_n^\#)}$ for $n\ge 0$, and thus a morphism
\[
\qcyc_{\bt\bt}\Cl\;\mrafl{25}{\qcyc_{\bt\bt} F}\; \qcyc_{\bt\bt} \Cl'
\]
of quasicyclic categories.

As variants, we mention
\[
\barcl
\qcyc^{\per}_{\bt\bt}\Cl     & := & \text{Iso}\left(\uul{8}{\Cl^+(\b\De_\bt^\#)}^{+,\,\per}(\b\De_\bt^\#)\right) \vspace*{1mm}\\
\qcyc^{\tht = 1}_{\bt\bt}\Cl & := & \text{Iso}\left(\uul{8}{\Cl^+(\b\De_\bt^\#)}^{+,\,\tht = 1}(\b\De_\bt^\#)\right)\; . \\
\ea
\]
Cf.\ Remark \ref{RemV1+2}.

\bq
 This procedure can be iterated to obtain triquasicyclic categories etc.
\eq

%% file: h3general.tex
\appendix

\begin{footnotesize}
\section{Some general lemmata}

\sbq
This appendix is a tool kit consisting of folklore lemmata (with proof) and known results (mainly without proof). We do not claim originality.
\seq

\subsection{An additive lemma}

Let $\Al$ and $\Bl$ be additive categories, and let $\Al\lraa{F}\Bl$ be a full and dense additive functor. Let $\Nl\tm\Bl$ be a full
additive subcategory. Let $\Ml\tm\Al$ be the full subcategory determined by
\[
\Ob\Ml\; :=\; \{ A\in\Ob\Al\; :\; \mb{$AF$ is isomorphic to an object of $\Nl$}\} \; .
\]

\begin{Lemma}
\label{Lem5}
Suppose that for each morphism $A\lraa{a_0} A'$ in $\Al$ such that $a_0 F = 0$, there exists a factorisation
\[
(A\lraa{a_0} A') \= (A\lraa{a_0'} M_0\lraa{a_0''} A')
\]
with $M_0\in\Ob\Ml$. Then the induced functor
\[
\barcl
\Al/\Ml        & \lraa{\ulF} & \Bl/\Nl \\
(A\lraa{a} A') & \lramaps    & (AF\lraa{aF}A'F) \\
\ea
\]
is an equivalence.
\end{Lemma}

{\it Proof.} We have to show that $\ulF$ is faithful. Suppose given $A\lraa{a} A'$ in $\Al$ such that 
\[
(AF\lraa{aF} A'F) \= (AF\lraa{b'} N \lraa{b''} A'F)\; ,
\]
where $N\in\Ob\Nl$. Since $F$ is dense, we may assume $N = MF$ for some $M\in\Ob\Ml$. Since $F$ is full, there exist $a'$ and $a''$ in $\Al$ with $a'F = b'$ and 
$a''F = b''$. Then
\[
(A\lraa{a} A') \= (A\lraa{a'a''} A') + (A\lraa{a_0} A')
\]
with $a_0 F = 0$. Since $a'a''$ factors over $M\in\Ob\Ml$, and since $a_0$ factors over an object of $\Ml$ by assumption on $F$, we conclude that
$a$ factors over an object of $\Ml$.\qed

\subsection{Exact categories}
\label{SecExCat}

\subsubsection{Definition}
\label{SecDefEx}

\sbq
 The concept of exact categories is due to {\sc Quillen} \bfcit{Qu73}, who uses a different, but equivalent set of axioms. In \mb{\bfcite{Ke90}{App.\ A}}, {\sc Keller} has cut down redundancies 
 in this set of axioms. We give still another equivalent reformulation.
\seq

An {\it additive category} $\Al$ is a category with zero object, binary products and binary coproducts such that the natural map from the coproduct to 
the product is an isomorphism; which allows to define a commutative and associative addition $(+)$ on $\liu{\Al}{(X,Y)}$, where $X,\, Y\,\in\,\Ob\Al$; and such that
there exists an endomorphism $-1_X$ for each $X\in\Ob\Al$ that is characterised by $1_X + (-1_X) = 0_X$.

A sequence $X\lraa{f} Y\lraa{g} Z$ in $\Al$ is called {\it short exact} if $f$ is a kernel of $g$ and $g$ is a cokernel of $f$.

A short exact sequence isomorphic to a short exact sequence of the form
\[
X\;\lraa{\smatez{1}{0}}\; X\ds Y \;\lrafl{35}{\smatze{0}{1}}\; Y\; , 
\]
where $X,\, Y\,\in\,\Ob\Al$, is called {\it split short exact.} A morphism appearing as a kernel in a split short exact sequence is {\it split monomorphic}, a morphism appearing as a cokernel in 
a split short exact sequence is called {\it split epimorphic}. A split short exact sequence is isomorphic to a sequence of the form just displayed by an isomorphism having an identity on the
first and on the third term.

An {\it exact category} $(\El,\Sl)$ consists of an additive category $\El$ and an isomorphism closed set $\Sl$ of short exact sequences in $\El$, called {\it pure short exact sequences} 
(\footnote{\scr This notion is borrowed from the particular cases of pure short exact sequences of lattices over orders and of $\ts$-pure short exact sequences of modules. 
Other frequently used names are {\it admissible short exact sequence}, consisting of an {\it admissible monomorphism} and an {\it admissible epimorphism}; and {\it conflation,} consisting 
of an {\it inflation} and a {\it deflation.}}), such that the following axioms \hyperlink{Ex.1}{(Ex 1, 2, 3, 1$^\circ$, 2$^\circ$, 3$^\circ$)} are satisfied.

A monomorphism fitting into a pure short exact sequence is called a {\it pure monomorphism}, denoted by $\lramono$; an epimorphism 
fitting into a pure short exact sequence is called a {\it pure epimorphism}, denoted by $\lraepi$. A morphism which can be written 
as a composition of a pure epimorphism followed by a pure monomorphism is called {\it pure}.\hypertarget{Ex.1}{}%

\begin{tabular}{lp{12cm}}
{\rm (Ex 1)}         & Split monomorphisms are pure monomorphisms. \hypertarget{Ex.1o}{}\\
{\rm (Ex 1$^\circ$)} & Split epimorphisms are pure epimorphisms. \hypertarget{Ex.2}{}\\
{\rm (Ex 2)}         & The composite of two pure monomorphisms is purely monomorphic. \hypertarget{Ex.2o}{}\\
{\rm (Ex 2$^\circ$)} & The composite of two pure epimorphisms is purely epimorphic. \hypertarget{Ex.3}{}\\
{\rm (Ex 3)} 
& Given a commutative diagram
  $$
  \xymatrix{
                                & Y\ar~+{|*\dir{|}}[dr] & \\
  X\ar[ur] \ar~+{|*\dir{*}}[rr] &                       & Z\; , \!\! \\
  }
  $$
  we may insert it into a commutative diagram
  $$
  \xymatrix{
  A\ar~+{|*\dir{*}}[dr] \ar~+{|*\dir{*}}[rr] &                                           & B \\
                                             & Y\ar~+{|*\dir{|}}[dr]\ar~+{|*\dir{|}}[ur] & \\
  X\ar~+{|*\dir{*}}[ur] \ar~+{|*\dir{*}}[rr] &                                           & Z \\
  }
  $$
  with $(X,Y,B)$ and $(A,Y,Z)$ pure short exact sequences. \hypertarget{Ex.3o}{}\\
\end{tabular}

\begin{tabular}{lp{12cm}}
{\rm (Ex 3$^\circ$)} & 
  Given a commutative diagram
  $$
  \xymatrix{
                                             & Y\ar[dr] & \\
  X\ar~+{|*\dir{*}}[ur] \ar~+{|*\dir{|}}[rr] &                       & Z\; , \!\! \\
  }
  $$
  we may insert it into a commutative diagram
  $$
  \xymatrix{
  A\ar~+{|*\dir{*}}[dr] \ar~+{|*\dir{|}}[rr] &                                           & B \\
                                             & Y\ar~+{|*\dir{|}}[dr]\ar~+{|*\dir{|}}[ur] & \\
  X\ar~+{|*\dir{*}}[ur] \ar~+{|*\dir{|}}[rr] &                                           & Z \\
  }
  $$
  with $(X,Y,B)$ and $(A,Y,Z)$ pure short exact sequences. \\
\end{tabular}

An {\it exact functor} from an exact category $(\El,\Sl)$ to an exact category $(\Fl,\Tl)$ is given by an additive functor $\El\lraa{F}\Fl$ such that
$\Sl F\tm\Tl$, where, by abuse of notation, $F$ also denotes the functor induced by $F$ on diagrams of shape $\bullet\lra\bullet\lra\bullet$.

Frequently, the exact category $(\El,\Sl)$ is simply referred to by $\El$.

\begin{Example}
\label{Ex5_0_5}\rm\Absit
\begin{itemize}
\item[(1)] An abelian category, equipped with the set of all short exact sequences as pure short exact sequences, is an exact category.
\item[(2)] If $\El$ is an exact category, so is $\El^\circ$, equipped with the pure short exact sequences of $\El$ considered as short exact sequences in $\El^\circ$, with the roles of kernel and 
cokernel interchanged.
\end{itemize}
\end{Example}

\begin{Example}
\label{Lem5_1}\rm
An additive category $\Al$, equipped with the set of split short exact sequences as pure short exact sequences, is an exact category.

In fact, \hyperlink{Ex.1}{(Ex 1, 2)} are fulfilled, and it remains to prove \hyperlink{Ex.3}{(Ex 3)}; then the dual axioms ensue by duality. Given
$$
\xymatrix{
                                              & X\dk Y\dk Z\ar~+{|*\dir{|}}[dr]^{\smatdz{1}{0}{0}{1}{0}{0}} & \\
X\ar[ur] \ar~+{|*\dir{*}}[rr]^{\smatez{1}{0}} &                                                             & X\dk Y\; , \!\! \\
}
$$
we get
$$
\xymatrix@!@R=-1mm{
Z\ar~+{|*\dir{*}}[dr]_{\smated{0}{0}{1}} \ar~+{|*\dir{*}}[rr]^{\smatez{0}{1}} &                                                                                      & Y\dk Z \\
                                                       & X\dk Y\dk Z\ar~+{|*\dir{|}}[dr]^(0.7){\smatdz{1}{0}{0}{1}{0}{0}}\ar~+{|*\dir{|}}[ur]_(0.7){\rsmatdz{0}{-a}{1}{0}{0}{1}} & \\
X\ar~+{|*\dir{*}}[ur]^{\smated{1}{0}{a}}\ar~+{|*\dir{*}}[rr]^{\smatez{1}{0}}  &                                                                                      & X\dk Y \; , \!\!\\
}
$$
where $X\lraa{a} Z$ is the third component of the given morphism $X\lra X\ds Y\ds Z$.
\end{Example}

\begin{Example}
\label{Rem5_2}\rm
Suppose given an exact category $\El$ and a category $D$. Let a short exact sequence $(X,Y,Z)$ in $\El(D)$ be pure if the sequence 
$(X_d,Y_d,Z_d)$ is a pure short exact sequence in $\El$ for all $d\in\Ob D$. Then $\El(D)$ is an exact category.
\end{Example}

\subsubsection{Embedding exact categories}
\label{SecGQL}

By a theorem 
\begin{itemize}
\item stated by {\sc Quillen} \bfcite{Qu73}{p.\ 100}, 
\item proven by {\sc Laumon} \bfcite{La83}{Th. 1.0.3}, 
\item re-proven by {\sc Keller} \bfcite{Ke90}{Prop.\ A.2},
\item where {\sc Quillen} resp.\ {\sc Keller} refer to \bfcit{Ga62} for a similar resp.\ an auxiliary technique, 
\end{itemize}
for any exact category $\El$, there exists an abelian category $\w\El$ containing $\El$ as a full subcategory closed under extensions, 
the pure short exact sequences in $\El$ being the short exact sequences in $\w\El$ with all three objects in $\Ob\El$. 

Conversely, suppose given an exact category $\El$ and a full subcategory $\El'\tm\El$ such that whenever $(X,Y,Z)$ is a pure short exact sequence in $\El$ 
with $X,\, Z\,\in\,\Ob\El'$, then also $Y\in\Ob\El'$. Then the subcategory $\El'$, equipped with the pure short exact sequences in $\El$ with
all three terms in $\Ob\El'$ as pure short exact sequences in $\El'$, is an exact category. 

\subsubsection{Frobenius categories\,: definitions}
\label{SecFrobCat}

\begin{Definition}
\label{DefFunFrob}\rm\Absit
\begin{itemize}
\item[(1)] A {\it bijective object} in an exact category $\El$ is an object $B$ for which $\liu{\El}{(B,-)}$ and $\liu{\El}{(-,B)}$ are exact functors from $\El$ resp.\ from $\El^\circ$ to $\Z\Modl$. 
\item[(2)] A {\it Frobenius category} is an exact category for which each object $X$ allows for a diagram $B\lraepi X\lramono B'$ with $B$ and $B'$ bijective.
\item[(3)] Suppose given an exact category $\Fl$ carrying a shift automorphism $\TTT: X\lramaps X\TTT = X^{+1}$ and two additive endofunctors $\III$ and $\PPP$ together with natural transformations
$1_\Fl\lraa{\io}\III$ and $\PPP\lraa{\pi} 1_\Fl$ such that $\TTT\PPP = \III$ and such that 
\[
X\lraa{X\io} X\III = X^{+1}\PPP\lraa{X^{+1}\pi} X^{+1}
\]
is a pure short exact sequence with bijective middle term for all $X\in\Ob\Cl$. Then $(\Fl,\TTT,\III,\io,\PPP,\pi)$ is called a {\it functorial Frobenius category.} Often we write just 
$\Fl$ for $(\Fl,\TTT,\III,\io,\PPP,\pi)$.
\end{itemize}
\end{Definition}

\begin{Example}
\label{ExFunFrob}\rm\Absit
\begin{itemize}
\item[(1)] Let $\Al$ be an additive category. Let $\dot\Z$ denote the {\sf discrete} category with $\Ob\dot\Z = \Z$ and only identical morphisms. The category $\Al(\dot\Z)$ carries the 
shift functor $X^\bt\lramaps X^{\bt+1}$, where $(X^{\bt+1})^i = X^{i+1}$.
An object in the category $\CC(\Al)$ of complexes with entries in $\Al$ is written $(X^\bt,d^\bt)$, where $X$ is an object of $\Al(\dot\Z)$ and where $X^\bt\lraa{d^\bt} X^{\bt+1}$ with 
$d^\bt d^{\bt+1} = 0$. The category $\CC(\Al)$, equipped with pointwise split short exact sequences, is an exact category; cf.\ Examples \ref{Lem5_1}, \ref{Rem5_2}.
Given a complex $(X^\bt,d^\bt)$, we let $(X^\bt,d^\bt)\TTT = (X^\bt,d^\bt)^{+1} := (X^{\bt+1},-d^{\bt+1})$ and
\[
\ba{rl}
   & \Big((X^\bt,d^\bt)\mra{(X^\bt,d^\bt)\io} (X^\bt,d^\bt)\III = (X^\bt,d^\bt)^{+1}\PPP\mra{(X^\bt,d^\bt)^{+1}\pi} (X^\bt,d^\bt)^{+1}\Big) \vspace*{2mm}\\
:= & \Big((X^\bt,d^\bt)\mrafl{25}{\smatez{1}{d^\bt}} \left(X^\bt\ds X^{\bt+1},\smatzz{0}{0}{1}{0}\right)\mrafl{35}{\smatze{-d^\bt}{1}} (X^{\bt+1},-d^{\bt+1})\Big)\; . \\
\ea
\] 
Then $(\CC(\Al),\TTT,\III,\io,\PPP,\pi)$ is a functorial Frobenius category.
\item[(2)] Suppose $\El$ to be a Frobenius category. Let $\Bl\tm\El$ be a {\it sufficiently big full subcategory of bijective objects,} i.e.\ each object of $\Bl$ is bijective in $\El$, and
each object $X$ of $\El$ admits $B\lraepi X\lramono B'$ with $B,\, B'\in\Ob\Bl$. In other words, each bijective object of $\El$ is isomorphic to a direct summand of an object of $\Bl$.

Let $\Bl^\ac\tm\CC(\Bl)$ denote the full subcategory of {\it purely acyclic} complexes, i.e.\ complexes $(X^\bt,d^\bt)$ such that all differentials $X^i\lraa{d} X^{i+1}$ are pure, factoring in $\El$ as 
$d = \b d \dot d$ with $\b d$ purely epi- and $\dot d$ purely monomorphic, and such that all resulting sequences $(\dot d,\b d)$ are purely short exact. For short, a complex is
purely acyclic if it decomposes into pure short exact sequences.

Then $\Bl^\ac$ is a functorial Frobenius category, equipped with the restricted functors and transformations of $\CC(\Bl)$ as defined in (1); cf.\ \bfcite{Ne90}{Lem.\ 1.1}. Let 
$\Bl^\spac\tm\Bl^\ac$ be the full subcategory of split acyclic complexes, i.e.\ of complexes isomorphic to a complex of the form $(T^\bt\ds T^{\bt+1},\smatzz{0}{0}{1}{0})$ for some 
\mb{$T^\bt\in\Ob\Bl(\dot\Z)$}. Then $\Bl^\spac$ is a sufficiently big full subcategory of bijective objects in $\Bl^\ac$.
\end{itemize}
\end{Example}

\begin{Definition}
\label{DefHot}\rm
Suppose given a Frobenius category $\El$, and a sufficiently big full subcategory $\Bl\tm\El$ of bijectives. Let
\[
\ba{rcll}
\ulEl  & := & \El/\Bl           & \mb{be the {\it classical stable category} of $\El$}\; ; \\
\uulEl & := & \Bl^\ac/\Bl^\spac & \mb{be the {\it stable category} of $\El$}\; . \\
\ea
\]
In other words, the stable category $\uulEl$ of $\El$ is defined to be the classical stable category $\ulk{\Bl^\ac}$ of $\Bl^\ac$. The shift functor induced by the automorphism $\TTT$ of $\Bl^\ac$
on $\uulEl$ is also denoted by $\TTT$.
\end{Definition}

\begin{Lemma}
\label{LemBothSideRes}
The functor
\[
\barcl
\Bl^\ac & \lraa{I} & \El \\
(X,d)   & \lramaps & \Img(X^0\lraa{d} X^1)
\ea
\]
induces an equivalence
\[
\uulEl \= \ulk{\Bl^\ac}\;\;\lraisoa{\ulk{I}}\;\ulEl\; .
\]
\end{Lemma}

Cf.\ \bfcite{Ke94}{Sec.\ 4.3}.

{\it Proof.} This is an application of Lemma \ref{Lem5}. \qed

We choose an inverse equivalence $R$ to $\ulk{I}$.
We have the residue class functor $\El\lraa{N}\ulEl$, and, by abuse of notation, a second residue class functor $(\El\lraa{N}\uulEl) := (\El\lraa{N}\ulEl\lraisoa{R}\uulEl)$.

\sbq
A morphism $X\lraa{f} Y$ is zero in $\ulEl$ if and only if for any monomorphism $X\lramonoa{i} X'$ and any epimorphism $Y'\lraepia{p} Y$, there is a factorisation
$f = if'p$. This defines $\ulEl$ without mentioning bijective objects in $\El$. So one might speculate whether the class of Frobenius categories within the class of exact categories could be extended 
still without losing essential properties of Frobenius categories.
\seq

\subsection{Kernel-cokernel-criteria}

Let $\Al$ be an abelian category. The {\it circumference lemma} states that given a commutative triangle in $\Al$, the induced sequence on kernels and cokernels, with zeroes attached to the ends, 
is long exact.

\begin{Definition}
\label{DefSquare}\rm
A {\it weak square} in $\Al$ is a commutative quadrangle $(A,B,C,D)$ in $\Al$ whose diagonal sequence 
$(A,B\ds C,D)$ is exact at $B\ds C$. It is denoted by a ``$\,+\,$''-sign in the commutative diagram,
\[
\xymatrix{
C \ar[r]                   & D \\
A \ar[r]\ar[u]\ar@{}[ur]|+ & B\ar[u] \; . \!\! \\
}
\]
A {\it pullback} is a weak square with first morphism in the diagonal sequence being monomorphic. It is denoted 
\[
\xymatrix{
C \ar[r]                              & D \\
A \ar[r]\ar[u]\ar@{}[ur]|(0.3)\hookpb & B\ar[u] \; . \!\! \\
}
\]
A {\it pushout} is a weak square with second morphism in the diagonal sequence being epimorphic. It is denoted 
\[
\xymatrix{
C \ar[r]                              & D \\
A \ar[r]\ar[u]\ar@{}[ur]|(0.7)\hookpo & B\ar[u] \; . \!\! \\
}
\]
A {\it square} is a commutative quadrangle that is a pullback and a pushout, i.e.\ that has a short exact diagonal sequence. It is denoted 
\[
\xymatrix{
C \ar[r]                      & D \\
A \ar[r]\ar[u]\ar@{}[ur]|\Box & B\ar[u] \; . \!\! \\
}
\]
\end{Definition}

\begin{Remark}
\label{LemKCcritSquares}
If a commutative quadrangle in $\Al$
\[
\xymatrix{
C \ar[r]^c          & D \\
A \ar[r]_b\ar[u]^a & B\ar[u]_d  \\
}
\]
is a square, then the induced morphism from the kernel of $A\lraa{a} C$ to the kernel of $B\lraa{d} D$ is an isomorphism and the induced morphism from 
the cokernel of $A\lraa{a} C$ to the cokernel of $B\lraa{d} D$ is an isomorphism.
\end{Remark}

{\it Proof.} 
If $(A,B,C,D)$ is a square, then the circumference lemma, applied to the commutative triangle
\[
\xymatrix{
C                                 &                               \\
                                  & B\ds C\ar[ul]_{\smatze{0}{1}} \\
A\ar[ur]_{\smatez{b}{a}}\ar[uu]^a & ,                             \\
}
\]
yields a long exact sequence
\[
0\;\lra\;\KK_a\;\lraa{j}\; B\;\lraa{-d} D\;\lraa{q}\;\CC_a\;\lra\; 0\; ,
\]
where $\KK_a\lraa{i} A$ is the kernel of $a$, and where $C\lraa{p}\CC_a$ is the cokernel of $a$. Since $ib = j$ and $cb = p$, the induced morphisms on the kernels and on the cokernels of $a$ and $d$ 
are isomorphisms.
\qed

\begin{Lemma}
\label{Rem0_5}
A commutative quadrangle in $\Al$
\[
\xymatrix{
C \ar[r]^c          & D \\
A \ar[r]_b\ar[u]^a & B\ar[u]_d  \\
}
\]
is a weak square if and only if the induced morphism $\KK_a\lra\KK_d$ from the kernel of $A\lraa{a} C$ to the kernel of $B\lraa{d} D$ is an epimorphism and the induced morphism 
$\CC_a\lra\CC_d$ from the cokernel of $A\lraa{a} C$ to the cokernel of $B\lraa{d} D$ is a monomorphism. 

It is a pullback if and only if $\KK_a\lraiso\KK_d$ and $\CC_a\lramono\CC_d$.

It is a pushout if and only if $\KK_a\lraepi\KK_d$ and $\CC_a\lraiso\CC_d$.

It is a square if and only if $\KK_a\lraiso\KK_d$ and $\CC_a\lraiso\CC_d$.
\end{Lemma}

{\it Proof.}  Let $A'$ be the pullback of $(C,B,D)$, and let $D'$ be the pushout of $(A',C,B)$. We obtain induced morphisms $A\lra A'$ and $D'\lramono D$. 
The circumference lemma, applied to $(B,D',D)$, shows $\CC_{B\to D'}\lramono\CC_{B\to D}$.

The quadrangle $(A,B,C,D)$ is a weak square if and only if $A\lraepi A'$; which in turn, by the circumference lemma applied to $(A,A',C)$, is equivalent to 
$\KK_{A\to C}\lraepi\KK_{A'\to C}$ and $\CC_{A\to C}\lraiso\CC_{A'\to C}$; which, by composition and by Remark \ref{LemKCcritSquares}, applied to the square $(A',B,C,D')$,
is equivalent to $\KK_{A\to C}\lraepi\KK_{B\to D}$ and $\CC_{A\to C}\lramono\CC_{B\to D}$.

The quadrangle $(A,B,C,D)$ is a pullback if and only if $A\lraiso A'$; which in turn, by the circumference lemma applied to $(A,A',C)$, is equivalent to 
$\KK_{A\to C}\lraiso\KK_{A'\to C}$ and $\CC_{A\to C}\lraiso\CC_{A'\to C}$; which, by composition and by Remark \ref{LemKCcritSquares}, applied to the square $(A',B,C,D')$,
is equivalent to $\KK_{A\to C}\lraiso\KK_{B\to D}$ and $\CC_{A\to C}\lramono\CC_{B\to D}$.

The quadrangle $(A,B,C,D)$ is a square if and only if $A\lraiso A'$ and $D'\lraiso D$; which in turn, by the circumference lemma applied to $(A,A',C)$, is equivalent to 
$\KK_{A\to C}\lraiso\KK_{A'\to C}$, $\CC_{A\to C}\lraiso\CC_{A'\to C}$ and $\CC_{B\to D'}\lraiso\CC_{B\to D}$; which, by composition and by Remark \ref{LemKCcritSquares}, applied to the square 
$(A',B,C,D')$, is equivalent to $\KK_{A\to C}\lraiso\KK_{B\to D}$ and $\CC_{A\to C}\lraiso\CC_{B\to D}$.
\qed

\subsection{An exact lemma}
\label{SubsubEL}

Let $\El$ be an exact category. A {\it pure square} in $\El$ is a commutative quadrangle $(A,B,C,D)$ in $\El$ that has a pure short exact diagonal sequence 
$(A,B\ds C,D)$. Just as a square in abelian categories, a pure square is denoted by a box ``$\,\Box\,$''.

\begin{Lemma}
\label{Lem2outof3}
Suppose given a composition
\[
\xymatrix{
X'\ar[r]      & Y'\ar[r]      & Z' \\
X\ar[r]\ar[u] & Y\ar[r]\ar[u] & Z\ar[u] \\
}
\]
of commutative quadrangles in $\El$. If two out of the three quadrangles $(X,Y,X',Y')$, $(Y,Z,Y',Z')$, 
\mb{$(X,Z,X',Z')$} are pure squares, so is the third. 
\end{Lemma}

{\it Proof.}
In an abelian category, this follows from Lemma \ref{Rem0_5}.

As explained in \S\ref{SecGQL}, we may embed $\El$ fully, faithfully and additively into an abelian category $\w\El$ such that the pure short 
exact sequences in $\El$ are precisely the short exact sequences in $\w\El$ with all three objects in $\Ob\El$. In particular, the pure squares in $\El$ are 
precisely the squares in $\w\El$ with all four objects in $\Ob\El$, and the assertion in $\El$ follows from the assertion in $\w\El$.
\qed

\subsection{Some abelian lemmata}
\label{SubsubTL}

Let $\Al$ be an abelian category. 

\begin{Lemma}
\label{LemDec}
Inserting images, a weak square $(A,B,C,D)$ in $\Al$ decomposes into 
\[
\xymatrix{
C\ar~+{|*\dir{|}}[r]                                              & \ar~+{|*\dir{*}}[r]                                              & D \\
\ar~+{|*\dir{|}}[r]\ar~+{|*\dir{*}}[u]\ar@{}[ur]|{\Box}           & \ar~+{|*\dir{*}}[r]\ar~+{|*\dir{*}}[u]\ar@{}[ur]|(0.42){\hookpb} & \ar~+{|*\dir{*}}[u] \\
A\ar~+{|*\dir{|}}[r]\ar~+{|*\dir{|}}[u]\ar@{}[ur]|(0.58){\hookpo} & \ar~+{|*\dir{*}}[r]\ar~+{|*\dir{|}}[u]\ar@{}[ur]|{\Box}          & B\ar~+{|*\dir{|}}[u]  \\
}
\]
\end{Lemma}

{\it Proof.} The assertion follows using the characterisation of weak squares, pullbacks and pushouts given in Lemma \ref{Rem0_5}.\qed

\begin{Lemma}
\label{Lem1}
If, in a commutative diagram
\[
\xymatrix{
X'\ar[r]                    & Y'\ar[r]                    & Z' \\
X\ar[u]\ar[r]\ar@{}[ur]|{+} & Y\ar[u]\ar[r]\ar@{}[ur]|{+} & Z\ar[u] \\
}
\]
in $\Al$, the quadrangles $(X,Y,X',Y')$ and $(Y,Z,Y',Z')$ are weak squares, then the composite quadrangle $(X,Z,X',Z')$ is also a weak 
square.
\end{Lemma}

{\it Proof.} The assertion follows using the characterisation of weak squares given in Lemma \ref{Rem0_5}.\qed

\begin{Lemma}
\label{Lem1_5}
If, in a commutative diagram
\[
\xymatrix{
X'\ar[r]                              & Y'\ar[r]      & Z' \\
X\ar[u]\ar[r]\ar@{}[ur]^(.7){\hookpo} & Y\ar[u]\ar[r] & Z\ar[u] \\
}
\]
in $\Al$, the left hand side quadrangle $(X,Y,X',Y')$ is a pushout, as indicated, and the outer quadrangle $(X,Z,X',Z')$ is a weak 
square, then the right hand side quadrangle $(Y,Z,Y',Z')$ is also a weak square.

If the left hand side quadrangle $(X,Y,X',Y')$ and the outer quadrangle $(X,Z,X',Z')$ are pushouts, 
then the right hand side quadrangle $(Y,Z,Y',Z')$ is also a pushout.
\end{Lemma}

{\it Proof.} This follows using Lemma \ref{Rem0_5}.\qed

\begin{Lemma}
\label{Lem1_7}
If, in a commutative quadrangle in $\Al$
\[
\xymatrix{
X'\ar[r]      & Y'\\
X\ar[u]\ar[r] & Y\ar[u]\; , \!\! \\
}
\]
the morphism $X\lra Y$ is an epimorphism and the morphism $X'\lra Y'$ is a monomorphism, then the quadrangle is a weak square.
\end{Lemma}

{\it Proof.} This follows using Lemma \ref{Rem0_5}, applied horizontally.\qed

\begin{Lemma}
\label{Lem2}
Given a commutative diagram
\[
\xymatrix{
0 \ar[r]      & Y'\ar[r]      & Z'\ar[r]      & W' \\
X\ar[u]\ar[r] & Y\ar[u]\ar[r] & Z\ar[u]\ar[r] & 0\ar[u] \\
}
\]
in $\Al$ such that $(X,Z,0,Z')$ and $(Y,0,Y',W')$ are weak squares, then $(Y,Z,Y',Z')$ is a weak square.
\end{Lemma}

{\it Proof.} This follows using Lemma \ref{Rem0_5}. \qed 

\begin{Lemma}
\label{Lem3}
Given a diagram
\[
\xymatrix{
0\ar[r]                             & Y''\ar[r]^{v''}                        & Z''           \\
X'\ar[u]\ar[r]^{u'}\ar@{}[ur]|{+}   & Y'\ar[r]^{v'}\ar[u]^{y'}\ar@{}[ur]|{+} & Z'\ar[u]^{z'} \\
X\ar[u]^{x}\ar[r]^{u}\ar@{}[ur]|{+} & Y\ar[r]\ar[u]^{y}\ar@{}[ur]|{+}        & 0\ar[u]       \\
}
\]
in $\Al$ consisting of weak squares, as indicated by $+$, the sequence
\[
X\;\lraa{xu'}\; Y'\;\lrafl{30}{\smatez{y'}{v'}}\; Y''\ds Z'\;\lrafl{45}{\rsmatze{v''}{-z'}}\; Z''\ru{8}
\]
is exact at $Y'$ and at $Y''\ds Z'$.
\end{Lemma}

{\it Proof.} At $Y'$, we reduce to the case $u$, $u'$, $x$ and $y$ monomorphic and $(X,Y,X',Y')$ being a pullback via Lemma \ref{LemDec}. Suppose given $T\lraa{t} Y'$ with $ty' = 0$ and $tv' = 0$.
First of all, there exist $T\lraa{a} X'$ and $T\lraa{b} Y$ such that $a u' = t = b y$. Thus there exists $T\lraa{c} X$ such that $c x = a$ and $c u = b$. In particular, $c x u' = a u' = t$. 
Hence a factorisation of $t$ over $xu'$ exists. Uniqueness follows by monomorphy of $xu'$.\qed  

\begin{Lemma}
\label{Lem3_1}
\indent
\begin{itemize}
\item[{\rm (1)}] Suppose given a weak square in $\Al$
\[
\xymatrix{
X'\ar[r]                     & Y' \\
X \ar[r]\ar[u]\ar@{}[ur]|{+} & Y\ar[u] \\
}
\]
with $X'$ bijective. If the images of $X\lra Y$, of $X\lra X'$ and of $Y\lra Y'$ are bijective, then the images of $X'\lra Y'$ and of $X\lra Y'$ are 
bijective, too.
\item[{\rm (2)}] Suppose given a weak square in $\Al$
\[
\xymatrix{
X'\ar[r]                     & Y' \\
X \ar[r]\ar[u]\ar@{}[ur]|{+} & Y\ar[u] \\
}
\]
with $Y$ bijective. If the images of $X'\lra Y'$, of $X\lra X'$ and of $Y\lra Y'$ are bijective, then the images of $X\lra Y$ and of $X\lra Y'$ are 
bijective, too.
\end{itemize}
\end{Lemma}

{\it Proof.} Ad (1). We decompose $(X,Y,X',Y')$ according to Lemma \ref{LemDec} and denote the image of $X\lra Y$ by $\Img_{X,Y}$, etc. 

The diagonal sequence of the square $(\Img_{X,Y},\; Y,\; \Img_{X,Y'},\; \Img_{Y,Y'})$ shows that $\Img_{X,Y'}$ is bijective. 

The diagonal sequence of the square $(\Img_{X,X'},\;\Img_{X,Y'},\; X',\;\Img_{X',Y'})$ shows that $\Img_{X',Y'}$ is bijective. \qed

\begin{Lemma}
\label{LemSplHelp}
Given a pullback
\[
\xymatrix{
X\ar[r]^f                                                                  & Y \\
X'\ru{3.5}\ar[r]_{f'}\ar@{ >->}~+{|*\dir{*}}[u]^x\ar@{}[ur]|(0.3){\hookpb} & Y'\ar[u]_y\; ,\!\!\! \\
}
\]
in $\Al$ with $Y'$ injective, the morphism $(X',Y')\lraa{(x,y)} (X,Y)$ is split monomorphic in $\Al(\De_1)$. More precisely, any
retraction for $x$ may be extended to a retraction for $(x,y)$.
\end{Lemma}

{\it Proof.}
Let $xx' = 1_{X'}$. We form the pushout.
\[
\xymatrix{
& & Y \\
X\ar[r]\ar[urr]^f                                                          & P\ar~+{|*\dir{*}}[ur] \\
X'\ru{3.5}\ar[r]_{f'}\ar@{ >->}~+{|*\dir{*}}[u]^x\ar@{}[ur]|(0.7){\hookpo} & Y'\ru{3.5}\ar@{ >->}~+{|*\dir{*}}[u]\ar[uur]_y \\
}
\]
There is an induced morphism $P\lra Y'$ such that $(X\lra P\lra Y') = (X\lraa{x'f'} Y')$ and such that 
$(Y'\lra P\lra Y') = (Y'\lraa{1_{Y'}} Y')$. Since $Y'$ is injective, we obtain a factorisation $(P\lra Y') = $ 
\mb{$(P\lramono Y\lra Y')$.}\qed

\begin{Lemma}
\label{Lem3_2}
Suppose given a morphism $X\lra Y$ of commutative quadrangles in $\Al$, i.e.\ a morphism in $\Al(\De_1\ti\De_1)$.
\begin{itemize}
\item[{\rm (1)}] If $X$ is a pushout and $Y$ is a weak square, then the cokernel of $X\lra Y$ is a weak square.
\item[{\rm (2)}] If $X$ is a weak square and $Y$ is a pullback, then the kernel of $X\lra Y$ is a weak square.
\end{itemize}
\end{Lemma}

{\it Proof.} Ad (1). A morphism of commutative quadrangles gives rise to a morphism of the diagonal sequences; namely from a sequence that is
exact in the middle and has an epimorphic second morphism, stemming from $X$, to a sequence that is exact in the middle, stemming from $Y$. In order to prove that the 
cokernel sequence is exact in the middle, we reduce by insertion of the image of the first morphism of the diagonal sequence and by an application of the circumference lemma to the case 
in which the sequence stemming from $Y$ has a monomorphic first morphism. Then the snake lemma yields the result.
\qed

\begin{Lemma}
\label{Lem3_2_5}
Suppose given a diagram 
\[
\xymatrix{
0  \ar[r]                   & Y''\ar[r]                  & Z'' \\
X' \ar[r]\ar[u]\ar@{}[ur]|+ & Y'\ar[r]\ar[u]\ar@{}[ur]|+ & Z'\ar[u]\\
X  \ar[r]\ar[u]\ar@{}[ur]|+ & Y \ar[r]\ar[u]\ar@{}[ur]|+ & Z \ar[u]\\
}
\]
in $\Al$, consisting of weak squares. The induced morphisms furnish a short exact sequence
\[
\Img(X\lra Z')\;\;\lra\;\; \Img(Y\lra Z')\;\;\lra\;\; \Img(Y\lra Z'')\; .
\]
\end{Lemma}

{\it Proof.} Abbreviate $\Img(X\lra Z')$ by $\Img_{X,Z'}$ etc. The morphism $\Img_{X,Z'}\lra \Img_{Y,Z'}$ is monomorphic by composition, and, dually, the morphism 
$\Img_{Y,Z'}\lra \Img_{Y,Z''}$ is epimorphic. Now since $\Img_{X,X'}\lra \Img_{X,Z'}$ is epimorphic and $\Img_{Y,Z''}\lra \Img_{Y'',Z''}$ is monomorphic, it suffices to show that 
\[
\Img_{X,X'}\;\lra\; \Img_{Y,Z'}\;\lra\; \Img_{Y'',Z''}
\]
is exact at $\Img_{Y,Z'}$. This follows from the diagram obtained by Lemma \ref{LemDec} 
\[
\xymatrix{
0 \ar[r]                                                      & Y''\ar~+{|(0.45)*\dir{|}}[r]                                             & \Img_{Y'',Z''} \\
X'\ar[r]\ar[u]\ar@{}[ur]|+                                    & Y'\ar~+{|(0.45)*\dir{|}}[r]\ar[u]\ar@{}[ur]|(0.56){\hookpo}               & \Img_{Y',Z'}\ar[u] \\
\Img_{X,X'}\ar[r]\ar~+{|*\dir{*}}[u]\ar@{}[ur]|(0.44){\hookpb} & \Img_{Y,Y'}\ar~+{|(0.45)*\dir{|}}[r]\ar~+{|*\dir{*}}[u]\ar@{}[ur]|{\Box} & \Img_{Y,Z'}\ar~+{|*\dir{*}}[u]\; , \!\!\\
}
\]
since by Lemma \ref{Lem1}, weak squares are stable under composition.
\qed

\subsection{On Frobenius categories}
\label{SubsubFL}

\subsubsection{Some Frobenius-abelian lemmata}

Suppose given an abelian Frobenius category $\Al$; cf.\ Definition \ref{DefFunFrob}. Let $\Bl$ be its full subcategory of bijective objects.
Recall that the classical stable category of $\Al$ is defined as $\ul{\Al} = \Al/\Bl$; cf.\ Definition \ref{DefHot}. A morphism in $\Al$ whose residue class in $\ul{\Al}$ is an isomorphism 
is called a {\it homotopism.} A morphism in $\Al$ whose residue class in $\ul{\Al}$ is a retraction is called a {\it retraction up to homotopy.}

\begin{Lemma}
\label{Lem3_3}
Given a retraction up to homotopy $X\lraa{f} Y$ and an epimorphism $Y'\lraepia{y} Y$ in $\Al$, in the pullback
\[
\xymatrix{
X' \ar~+{|*\dir{|}}[r]^*+[u]{\scm x}\ar[d]_{f'}\ar@{}[dr]|(0.3){\hookpbo} & X \ar[d]^f \\
Y' \ar~+{|*\dir{|}}[r]^*+[u]{\scm y}                                      & Y \; ,\!\!\!\\        
}
\]
the morphism $X'\lraa{f'} Y'$ is a retraction up to homotopy, too. More precisely, if $gf \con_\Bl 1_Y$, then we may find a morphism $g'$ with $g' f'\con_\Bl 1_{Y'}$ as a pullback of 
$g$ along $x$.
\end{Lemma}

{\it Proof.} Let $Y\lraa{g} X$ be such that $gf = 1_Y + h$, where 
\[
(Y\lraa{h} Y) \= (Y\lraa{h_1} B \lraa{h_2} Y)
\]
for some $B\in\Ob\Bl$ and some morphisms $h_1$ and $h_2$ in $\Al$. Let $B\lraa{h'_2} Y'$ be a morphism such that
\[
(B\lraa{h'_2} Y' \lraepia{y} Y) \= (B\lraa{h_2} Y)\; ,
\]
which exists since $B$ is projective and $y$ is epimorphic. The commutative quadrangle
\[
\xymatrix{
Y' \ar~+{|*\dir{|}}[r]^*+[u]{\scm y}\ar[d]_{1_{Y'} + y h_1 h'_2}            & Y\ar[d]^{1_Y + h} \\        
Y' \ar~+{|*\dir{|}}[r]^*+[u]{\scm y}                                        & Y \\        
}
\]
is a pullback since the induced morphism on the horizontal kernels is an identity; cf.\ Lemma \ref{Rem0_5}. So we may form the diagram
\[
\xymatrix{
Y' \ar~+{|*\dir{|}}[r]^*+[u]{\scm y}\ar[d]_{g'}\ar@{}[dr]|(0.3){\hookpbo} & Y \ar[d]^g \\        
X' \ar~+{|*\dir{|}}[r]^*+[u]{\scm x}\ar[d]_{f'}\ar@{}[dr]|(0.3){\hookpbo} & X \ar[d]^f \\
Y' \ar~+{|*\dir{|}}[r]^*+[u]{\scm y}                                      & Y \; ,\!\!\!\\        
}
\]
in which $g'$ with $g'x = yg$ and $g'f' = 1_{Y'} + y h_1 h'_2$ is induced by the universal property of the lower pullback $(X',X,Y',Y)$, and in which 
the resulting upper quadrangle $(Y',Y,X',X)$ is a pullback by Lemma \ref{Rem0_5}.\qed

\begin{Lemma}
\label{Lem3_4}
Given a homotopism $X\lraa{f} Y$ and an epimorphism $Y'\lraepia{y} Y$ in $\Al$, in the pullback
\[
\xymatrix{
X' \ar~+{|*\dir{|}}[r]^*+[u]{\scm x}\ar[d]_{f'}\ar@{}[dr]|(0.2){\hookpbo} & X \ar[d]^f \\
Y' \ar~+{|*\dir{|}}[r]^*+[u]{\scm y}                                      & Y \; ,\!\!\!\\        
}
\]
the morphism $X'\lraa{f'} Y'$ is a homotopism, too.
\end{Lemma}

{\it Proof.} Let $gf \con_\Bl 1_Y$ and $fg\con_\Bl 1_X$. By Lemma \ref{Lem3_3}, we may form the diagram
\[
\xymatrix{
Y' \ar~+{|*\dir{|}}[r]^*+[u]{\scm y}\ar[d]_{g'}\ar@{}[dr]|(0.3){\hookpbo} & Y \ar[d]^g \\        
X' \ar~+{|*\dir{|}}[r]^*+[u]{\scm x}\ar[d]_{f'}\ar@{}[dr]|(0.3){\hookpbo} & X \ar[d]^f \\
Y' \ar~+{|*\dir{|}}[r]^*+[u]{\scm y}                                      & Y \; ,\!\!\!\\        
}
\]
in which $g' f' \con_\Bl 1_{Y'}$.
Since $g$ is a retraction up to homotopy, so is $g'$ by Lemma \ref{Lem3_3}. Therefore $g'$ is a homotopism. Hence also $f'$ is a homotopism.\qed

\subsubsection{Decomposing split diagrams in intervals}
\label{SecDecSplit}

Let $\Al$ be an abelian Frobenius category, and let $\Bl$ be its full subcategory of bijective objects.
Suppose given $n\ge 1$. Write $\dDe_n := \De_n\ohne\{ 0\}$. An object $X$ in $\Al(\dDe_n)$ is called {\it split} if
$X_k\lraa{x} X_l$ is split for all $k,\, l\,\in\, [1,n]$ with $k \le l$. 

Given $C\in\Ob\Al$ and $k,\, l\,\in\, [1,n]$ with $k \le l$, we denote by $C_{[k,l]}$ the object of $\Al(\dDe_n)$ given by 
$(C_{[k,l]})_j = 0$ for $j\in [1,n]\ohne [k,l]$, by $(C_{[k,l]})_j = C$ for $j\in [k,l]$, and
by $\big((C_{[k,l]})_{j}\lraa{c} (C_{[k,l]})_{j'}\big) = (C\lraa{1_C} C)$ for $j,\, j'\,\in\, [k,l]$ with $j \le j'$.
An object in $\Al(\dDe_n)$ of the form $C_{[k,l]}$ for some $C\in\Ob\Al$ and some $k,\, l\,\in\, [1,n]$ with $k \le l$, is called an {\it interval.}

\begin{Lemma}
\label{LemSpl}
Any split object in $\Bl(\dDe_n)$ is isomorphic to a finite direct sum of intervals.
\end{Lemma}

{\it Proof.} We proceed by induction on $n$. Suppose given a split object $X$ in $\Bl(\dDe_n)$. Let $X' := X\ind_0$ be defined as a 
pointwise pullback at $n$, using $0\lraa{0} X_n$ (cf.\ \S\ref{SubsecPPP} below). We have $X'\in\Ob\Bl(\dDe_n)$ with $X'_n = 0$. Hence, by induction, $X'$ is isomorphic to 
a finite direct sum of intervals. There is a pure 
monomorphism $X'\lramono X$ whose cokernel is a diagram in $\Ob\Bl(\dDe_n)$ consisting of split monomorphisms; cf.\ Lemma \ref{Rem0_5}. Moreover, by an iterated application of
Lemma \ref{LemSplHelp}, starting at position $1$, this pure monomorphism $X'\lramono X$ is split as a morphism of $\Al(\dDe_n)$ (\footnote{\scr At this point, we use that
$\dDe_n$ is linearly ordered.}).  Thus $X$ is isomorphic to the direct sum of $X'$ and the cokernel of $X'\lramono X$, and it remains to be shown that 
this cokernel is isomorphic to a finite direct sum of intervals.

Therefore, we may assume that $X$ consists of split monomorphisms $\xymatrix{X_k\;\ar@{ >->}~+{|*\dir{*}}[r]^x & X_l}$ for \mb{$k,\, l\,\in\, [1,n]$.}
We have a monomorphism $(X_1)_{[1,n]}\lramonoa{i} X$. Choosing a retraction to $\xymatrix{X_1\;\ar@{ >->}~+{|*\dir{*}}[r]^x & X_n}$ and composing, we obtain a coretraction to $i$,
so that $X$ is isomorphic to the direct sum of the interval $(X_1)_{[1,n]}$ and the cokernel of $i$. Since the cokernel of $i$ has a zero term at position $1$, we are done by induction.\qed

\subsubsection{A Freyd category reminder}
\label{SecFreyd}

\sbq
 The construction of the Freyd category and its properties are due to {\sc Freyd} \bfcite{Fr65}{Th.\ 3.1}.
\seq

\begin{Definition}
\label{DefFreydic}\rm
Suppose given an additive category $\Cl$ and a morphism $X\lraa{f} Y$ in $\Cl$.

\begin{itemize}
\item[(1)] A morphism $K\lraa{i} X$ is a {\it weak  kernel} of $X\lraa{f} Y$ if the sequence of abelian groups
\[
(T,K)\;\lraa{(-)i}\; (T,X)\;\lraa{(-)f}\; (T,Y)
\]
is exact at $(T,X)$ for every $T\in\Ob\Cl$.

\item[(2)] A morphism $Y\lraa{p} C$ is a {\it weak  cokernel} of $X\lraa{f} Y$ if the sequence of abelian groups
\[
(X,T)\;\llaa{f(-)}\; (Y,T)\;\llaa{p(-)}\; (C,T)
\]
is exact at $(Y,T)$ for every $T\in\Ob\Cl$.

\item[(3)] The category $\Cl$ is called {\it weakly abelian} if every morphism has a weak kernel and a weak cokernel, and if every morphism is a weak kernel (of some morphism) and a 
weak cokernel (of some morphism).
\end{itemize}
\end{Definition}

Let $\Cl$ be a weakly abelian category. Let $\Cl^0(\De_1)$ be the full subcategory of $\Cl(\De_1)$ whose objects are zero morphisms. The {\it Freyd category} $\h\Cl$ of $\Cl$ is defined to be the 
quotient category
\[
\h\Cl\; :=\; \Cl(\De_1)/\Cl^0(\De_1)\; .
\]
We collect some elementary facts and constructions and mention some conventions.

\begin{itemize}
\item[(1)] The category $\h\Cl$ is abelian. The kernel and the cokernel of a morphism $X\lraa{f} Y$ represented by $(f',f'')$ are constructed as 
\[
\xymatrix{
K\ar[r]^i\ar[d]_{ix} & X'\ar[d]_x \ar[r]^{f'} & Y'\ar[d]^{y}\ar[r]^{1_{Y'}} & Y'\ar[d]^{yp} \\
X''\ar[r]^{1_{X''}}  & X''\ar[r]^{f''}        & Y''\ar[r]^p                 & C\; , \!\! \\
}
\]
where $i$ is a chosen weak kernel and $p$ a chosen weak cokernel of the diagonal morphism $f'y = xf''$. If $f'y = xf'' = 0$, we choose $X'\lraa{1_{X'}} X'$ as weak kernel and 
$Y''\lraa{1_{Y''}} Y''$ as weak cokernel. 

Choosing a kernel and a cokernel for each object in $\h\Cl(\De_1)$, we obtain a kernel and a cokernel functor $\h\Cl(\De_1)\lradouble\h\Cl$, as for any abelian category.

\item[(2)] We stipulate that the pullback resp.\ the pushout of an identity morphism along a morphism is chosen to be an identity morphism.
\item[(3)] We have a full and faithful functor $\Cl\lra\h\Cl$, $X\lramaps (X\lraa{1_X} X)$. Its image, identified with $\Cl$, consists of bijective objects.
\item[(4)] For each $X = (X'\lraa{x} X'')\in\Ob\h\Cl$, we may define objects and morphisms 
\[
X\PPP\;\lraepia{X\pi}\; X\;\lramonoa{X\io}\; X\III
\]
by
\[
\xymatrix{
X'\ar[r]^{1_{X'}}\ar[d]_{1_{X'}} & X'\ar[r]^x\ar[d]^x  & X''\ar[d]^{1_{X''}} \\
X'\ar[r]^x                       & X''\ar[r]^{1_{X''}} & X'' \; . \!\!\\
}
\]
As already mentioned in (3), the objects $X\PPP$ and $X\III$ are bijective, and thus $\h\Cl$ is Frobenius.

Sometimes, we write just $\io$ for $X\io$ and $\pi$ for $X\pi$. Note that $X\pi = 1_X$ and $X\io = 1_X$ if $X\in\Ob\Cl$.

\sbq
This construction $X\PPP\lraepia{X\pi} X\lramonoa{X\io} X\III$ is {\sf not} meant to be functorial in $(X'\lraa{x} X'')$, however.
\seq
\end{itemize}

\begin{Remark}
\label{RemBal}
Suppose given morphisms $X\lraa{f} Y\lraa{g} Z$ in $\,\Cl$. The following assertions are equivalent.
\begin{itemize}
\item[{\rm (i)}] The morphism $f$ is a weak kernel of $g$.
\item[{\rm (ii)}] The morphism $g$ is a weak cokernel of $f$.
\item[{\rm (iii)}] The sequence $(f,g)$ is exact at $Y$ when considered in $\h\Cl$.
\end{itemize}
\end{Remark}

{\it Proof.} Ad (i) $\imp$ (iii). Suppose that $f$ is a weak kernel of $g$. Let $K\lraa{i} Y$ be the kernel of $g$ in $\h\Cl$. Factor $f = f'i$. 
Since $f$ is a weak kernel of $g$ in $\Cl$, we may factor $(K\pi)i = uf$, whence $K\pi = uf'$. Hence $f'$ is epimorphic.

Ad (iii) $\imp$ (i). Suppose $(f,g)$ to be exact at $Y$. Let $T\lraa{t} Y$ in $\Cl$ be such that $tg = 0$. Then $t$ factors over the kernel of $g$, taken in $\h\Cl$, and therefore, by 
projectivity of $T$ in $\h\Cl$, also over $X$. \qed

\begin{Remark}
\label{RemSpl}
A morphism $X\lraa{f} Y$ in $\Cl$ is monomorphic if and only if it is a coretraction. Dually, it is epimorphic if and only if it is a retraction. 
\end{Remark}

{\it Proof.} Suppose $f$ to be monomorphic in $\Cl$. It suffices to show that $f$ is monomorphic in $\h\Cl$, for then $f$ is a coretraction since $X$ is injective in $\h\Cl$.
Let $K\lraa{i} X$ be the kernel of $f$ in $\h\Cl$. From $(K\pi)if = 0$, we conclude $(K\pi)i = 0$ since $f$ is monomorphic in $\Cl$, and thus $K\iso 0$
since $K\pi$ is epimorphic and $i$ is monomorphic in $\h\Cl$.\qed

\sbq
In particular, an abelian category is weakly abelian if and only if it is semisimple, i.e.\ if and only if every morphism in $\Al$ splits. Hence
the notion ``weakly abelian'' is slightly abusive.
\seq

Let $\El$ be a Frobenius category; cf.\ \S\ref{SecFrobCat}.

\begin{Lemma}
\label{LemExToWEx1}
Suppose given a pure short exact sequence $X'\lramonoa{i} X\lraepia{p} X''$ in $\El$. In $\ulEl$, the residue class $iN$ is a weak kernel of $pN$, and the residue class $pN$ is a weak cokernel
of $iN$.
\end{Lemma}

{\it Proof.} By duality, it suffices to show that $iN$ is a weak kernel of $pN$. So suppose given $T\lraa{t} X$ in $\El$ with $tp\con_\Bl 0$. We have to show that there exists a morphism 
$T\lraa{t'} X'$ such that $t'i \con_\Bl t$. Let $(T\lraa{t} X\lraa{p} X'') = (T\lraa{u} B\lraa{q} X'')$, where $B$ is bijective. Let $P\lraepia{\w p} B$ be the pullback of $p$ along $q$. We 
have a factorisation $(T\lraa{t} X) = (T\lraa{v} P\lraa{w} X)$. We have a factorisation $(X'\lramonoa{i} X) = (X'\lramonoa{\w i} P\lraa{w} X)$; moreover, $(\w i,\w p)$ is a pure short exact sequence, 
hence split by projectivity of $B$; cf.\ Lemma \ref{Rem0_5} and \S\ref{SecGQL}. Let $\w i r = 1$. Then $r\w i - 1\con_\Bl 0$, since it factors over $B$. We obtain 
$(vr)i = vr\w i w \con_\Bl vw = t$.\qed

\begin{Remark}
\label{RemFrobWeak}
The stable category $\uulEl$ and the classical stable category $\ulEl$ of the Frobenius category $\El$ are weakly abelian. 
The stable category $\uulEl$ carries an automorphism $\TTT$, induced by shifting an acyclic complex to the left by one position and negating the differentials.
\end{Remark}

Cf.\ Definition \ref{DefHot}.

{\it Proof.} By Lemma \ref{LemBothSideRes}, it remains to prove that $\ulEl$ is weakly abelian. Suppose given a morphism $X\lraa{f} Y$ in $\El$.  By duality, it suffices to show that the residue 
class of $X\lraa{f} Y$ in $\ulEl$ is a weak cokernel and has a weak kernel. Substituting isomorphically in $\ulEl$ by adding a bijective object to $X$, we may assume $f$ to be a pure epimorphism 
in $\El$. So we may complete to a pure short exact sequence and apply Lemma \ref{LemExToWEx1}. \qed

\begin{Lemma}
\label{LemExToWEx}
A pure short exact sequence $X'\lramonoa{i} X\lraepia{p} X''$ in $\El$ is mapped via the residue class functor $N$ to a sequence in $\ulEl$ that is exact at $X$ when considered in the Freyd category 
$\ulk{\h\El}$ of $\ulEl$. In particular, a pure square in $\El$ is mapped to a weak square in $\ulEl$.
\end{Lemma}

{\it Proof.} By Remark \ref{RemBal}, we may apply Lemma \ref{LemExToWEx1}.\qed

\subsubsection{Heller operators for diagrams}

\sbq
In Definition \ref{Def4}, the central role is attributed to the tuple $\tht = (\tht_n)_{n\ge 0}$ of isomorphisms. In the case of $\Cl$ being the stable category 
of a Frobenius category, such an isomorphism $\tht_n$ arises from different choices
of pure monomorphisms into bijective objects. To that end, we 
provide a comparison lemma, which suitably organises wellknown facts.
\seq

Let $C$ be a category.

Given a category $\Dl$ and a full subcategory $\Ul\tm\Dl(C)$, we say that $\Ul$ is {\it characteristic} in $\Dl(C)$ if the image of 
$\Ul$ under $A(C)$ is contained in $\Ul$ for any autoequivalence $\Dl\lraisoa{A}\Dl$, and if $\Ul$ is closed under isomorphy in $\Dl(C)$,
i.e.\ $X\iso X'$ in $\Dl(C)$ and $X'\in\Ob\,\Ul$ implies $X\in\Ob\,\Ul$.

Let $\El$ be a Frobenius category. Denote by $\ulEl$ its classical stable category, and denote by $\El\lraa{N}\ulEl$ the residue class
functor. Let $\Gl\tm\El(C)$ be a full additive subcategory. Let $\Hl\tm\ulEl(C)$ be a full additive characteristic subcategory such that $(\Gl)(N(C)) \tm \Hl$.
\[
\xymatrix{
\Gl\;\ar@{^{(}->}[r]\ar[d] & \El(C)\ar[d]^{N(C)}  \\
\Hl\;\ar@{^{(}->}[r]       & \ulEl(C) \\
}
\]
A {\it $C$-resolving system} $I$ consists of pure short exact sequences
\[
I \= \left(\left(X_c\lramonoa{i_{X,c}} I_{X,c}\lraepia{p_{X,c}} \w X_c\right)_{\!\! c\in\Ob C}\right)_{\!\! X\in\Ob\Gl}\; ,
\]
with bijective objects $I_{X,c}$ in $\El$ as middle terms.

\pagebreak

\begin{Lemma}
\label{LemComparison}\Absit
\begin{itemize}
\item[{\rm (1)}] Given a $C$-resolving system
\[
I \= \left(\left(X_c\lramonoa{i_{X,c}} I_{X,c}\lraepia{p_{X,c}} \w X_c\right)_{\!\! c\in\Ob C}\right)_{\!\! X\in\Ob\Gl}\; ,
\]
there exists a functor
\[
\Gl\;\lraa{\TTT_I}\;\Hl
\]
that is uniquely characterised by the following properties. 

On objects $X\in\Ob\Gl\tm\Ob\El(C)$, the image $X\TTT_I\in\Ob\Hl\tm\Ob\ulEl(C)$ is characterised as follows.

\begin{center}
\begin{tabular}{lp{13cm}}
$(\ast)$ &
For any $(c\lraa{\ga} d)\in C$, there exist
\begin{itemize}
\item[$\bullet$] a representative $(X \TTT_I)_\ga^\sim$ in $\El$ of the evaluation
$(X \TTT_I)_c\lrafl{30}{(X \TTT_I)_\ga} (X \TTT_I)_d$ in $\ulEl$ at $c\lraa{\ga} d$ of the diagram $X\TTT_I\in\Ob\Hl\tm\Ob\ulEl(C)$, and 
\item[$\bullet$] a morphism $I_{X,c}\lra I_{X,d}$ in $\El$ 
\end{itemize}
such that 
\[
\xymatrix{
X_c \ar~+{|(0.45)*\dir{*}}[r]^(0.4){i_{X,c}}\ar[d]_{X_\ga} & I_{X,c}\ar~+{|(0.4)*\dir{|}}[r]^(0.4){p_{X,c}}\ar[d] & (X\TTT_I)_c\ar[d]^{(X\TTT_I)_\ga^\sim} \\
X_d \ar~+{|(0.45)*\dir{*}}[r]^(0.4){i_{X,d}}               & I_{X,d}\ar~+{|(0.4)*\dir{|}}[r]^(0.4){p_{X,d}}       & (X\TTT_I)_d \\
}
\]
is a morphism of pure short exact sequences. \\
\end{tabular}
\end{center}

On morphisms $(X\lraa{f}Y)\in\Gl\tm\El(C)$, the image $(X\TTT_I\lrafl{30}{f\TTT_I} Y\TTT_I)\in\Hl\tm\ulEl(C)$ is characterised as follows.

\begin{center}
\begin{tabular}{lp{13cm}}
$(\ast\ast)$ &
For any $c\in\Ob C$, there exist
\begin{itemize}
\item[$\bullet$] a representative $(f \TTT_I)_c^\sim$ in $\El$ of the evaluation 
$(X \TTT_I)_c\lrafl{30}{(f \TTT_I)_c} (Y \TTT_I)_c$ in $\ulEl$ at $c$ of the diagram morphism $(X\TTT_I\lrafl{30}{f\TTT_I} Y\TTT_I)\in\Hl\tm\ulEl(C)$, and 
\item[$\bullet$] a morphism $I_{X,c}\lra I_{Y,c}$ in $\El$ 
\end{itemize}
such that 
\[
\xymatrix{
X_c \ar~+{|(0.45)*\dir{*}}[r]^(0.4){i_{X,c}}\ar[d]_{f_c} & I_{X,c}\ar~+{|(0.4)*\dir{|}}[r]^(0.4){p_{X,c}}\ar[d] & (X\TTT_I)_c\ar[d]^{(f\TTT_I)_c^\sim} \\
Y_c \ar~+{|(0.45)*\dir{*}}[r]^(0.4){i_{Y,c}}             & I_{Y,c}\ar~+{|(0.4)*\dir{|}}[r]^(0.4){p_{Y,c}}       & (Y\TTT_I)_c \\
}
\]
is a morphism of pure short exact sequences. \\
\end{tabular}
\end{center}

\item[{\rm (2)}] Given $C$-resolving systems
\[
\barcl
I  & = & \left(\left(X_c\lramonoa{i_{X,c}} I_{X,c}\lraepia{p_{X,c}} \w X_c\right)_{\!\! c\in\Ob C}\right)_{\!\! X\in\Ob\Gl}\; ,\vss\\
I' & = & \left(\left(X_c\lramonoa{i'_{X,c}} I'_{X,c}\lraepia{p'_{X,c}} \w X'_c\right)_{\!\! c\in\Ob C}\right)_{\!\! X\in\Ob\Gl}\; ,\\ 
\ea
\]
there exists an isomorphism
\[
\TTT_I\;\lraisofl{30}{\al_{I,I'}} \TTT_{I'}
\]
that is uniquely characterised by the following property.

\begin{center}
\begin{tabular}{lp{13cm}}
$(\ast\!\ast\!\ast)$ &
For any $X\in\Ob\Gl\tm\Ob\El(C)$ and for any $c\in\Ob C$, there exist
\begin{itemize}
\item[$\bullet$] a representative $(X\al_{I,I'})^\sim$ in $\El$ of the evaluation
$(X \TTT_I)_c\;\lrafl{30}{(X\al_{I,I'})_c}\; (X \TTT_{I'})_c$ in $\ulEl$ at $c$ of the evaluation $X\TTT_I\lrafl{30}{X\al_{I,I'}} X\TTT_{I'}$ in $\Hl\tm\ulEl(C)$ of 
$\al_{I,I'}$ at $X$, and 
\item[$\bullet$] a morphism $I_{X,c}\lra I'_{X,c}$ in $\El$
\end{itemize}
such that 
\[
\xymatrix{
X_c \ar~+{|(0.45)*\dir{*}}[r]^(0.4){i_{X,c}}\ar@{=}[d] & I_{X,c}\ar~+{|(0.4)*\dir{|}}[r]^(0.4){p_{X,c}}\ar[d] & (X\TTT_I)_c\ar[d]^{(X\al_{I,I'})^\sim} \\
X_c \ar~+{|(0.45)*\dir{*}}[r]^(0.4){i'_{X,c}}          & I'_{X,c}\ar~+{|(0.4)*\dir{|}}[r]^(0.4){p'_{X,c}}     & (X\TTT_{I'})_c \\
}
\]
is a morphism of pure short exact sequences. \\
\end{tabular}
\end{center}
\end{itemize}
\end{Lemma}

{\it Proof.} Let us first assume that $\Hl = \ulEl(C)$. Having proven all assertions in this case, it then finally will remain to be shown that given 
$\Hl\tm\ulEl(C)$ and a $C$-resolving system $I$, we have $X\TTT_I\in\Ob\Hl\tm\Ob\ulEl(C)$ for $X\in\Ob\Gl$.

We remark that starting from a morphism $U\lraa{u} U'$ in $\El$ and from chosen pure short exact sequences $(U,B,V)$ and $(U',B',V')$ with bijective
middle terms $B$ resp.\ $B'$, we may define a morphism $V\lraa{v} V'$ in $\ulEl$ by the existence
of a morphism
\[
\xymatrix{
U \ar~+{|*\dir{*}}[r]^{i\ru{-1}}\ar[d]_u  & B\ar~+{|*\dir{|}}[r]^{p\ru{-1}}\ar[d] & V\ar[d]^{v^\sim\ru{-1}} \\
U' \ar~+{|*\dir{*}}[r]^{i'\ru{-1}}        & B'\ar~+{|*\dir{|}}[r]^{p'\ru{-1}}     & V' \\
}
\]
of pure short exact sequences in $\El$, where $V\lraa{v}V'$ is the image in $\ulEl$ of the morphism $V\lraa{v^\sim} V'$ in $\El$.

Ad (1). Given $X\in\Ob\Gl$, we define $X\TTT_I\in\ulEl(C)$ at the morphism $c\lraa{\ga} d$ of $C$ by the diagram in $(\ast)$. The characterisation $(\ast)$ shows that $X\TTT_I$ is in fact
in $\Ob\ulEl(C)$.

Given a morphism $X\lraa{f} Y$ in $\Gl$, we define the morphism $X\TTT_I\lraa{f\TTT_I} Y\TTT_I$ in $\ulEl(C)$ at 
\mb{$c\in\Ob C$} by the diagram in $(\ast\ast)$. Combining $(\ast)$ and $(\ast\ast)$, we see that $f\TTT_I$ is in fact in $\ulEl(C)$. From $(\ast\ast)$ we conclude that $\TTT_I$ is indeed a functor.

Ad (2). Given $X\in\Ob\Gl$, we define $X\TTT_I\lrafl{30}{X\al_{I,I'}} X\TTT_{I'}$ at $c\in\Ob C$ by the diagram in $(\ast\!\ast\!\ast)$.

Combining $(\ast\!\ast\!\ast)$ and $(\ast)$, we see that $X\TTT_I\lrafl{30}{X\al_{I,I'}} X\TTT_{I'}$ is indeed in $\ulEl(C)$. Combining $(\ast\!\ast\!\ast)$ and $(\ast\ast)$, we see
that $\al_{I,I'}$ is indeed a transformation.

Suppose given resolving systems $I$, $I'$ and $I''$.
The characterisation of $\al_{I,I'}$ etc.\ implies that $\al_{I,I'}\al_{I',I''} = \al_{I,I''}$ and that $\al_{I,I} = 1_{\TTT_I}$. Hence in particular, 
$\al_{I,I'}\al_{I',I} = 1_{\TTT_I}$ and $\al_{I',I}\al_{I,I'} = 1_{\TTT_{I'}}$, and so $\al_{I,I'}$ is an isomorphism from $\TTT_I$ to $\TTT_{I'}$.

Consider the case $C = \De_0$, i.e.\ the terminal category, let $\Gl = \El(\De_0) = \El$ and let $\Hl = \ulEl(\De_0) = \ulEl$. 
For a $\De_0$-resolving system $J$, we obtain a functor $\El\lraa{\TTT_J}\ulEl$ that factors as
\[
(\El\lraa{\TTT_J}\ulEl) \= (\El\lraa{N} \ulEl \lraisoa{\b \TTT_J} \ulEl)\; .
\]
In fact, for a morphism $b$ that factors over a bijective object $B$, we can choose $0$ as a representative of $b \TTT_J$, inserting the pure
short exact sequence $(B,B,0)$. Moreover, $\b \TTT_J$ is an equivalence, for it is full; faithful, using the dual of the argument just given;
and dense, since given a morphism of short exact sequences in $\El$ with bijective middle terms and an identity on the kernels, the morphism
on the cokernels is a homotopism.

Now return to the general case $\Hl\tm\ulEl(C)$.
Let $J'$ be a $C$-resolving system consisting of pure short exact sequences with bijective middle term that already occur in the chosen $\De_0$-resolving system 
$J$. Then, for $X\in\Ob\Gl$, we have $X\TTT_{J'} = X(N(C))(\b \TTT_J(C))$. Since $X(N(C))\in\Ob\Hl$ by assumption, and since, moreover, $\Hl$ is assumed to be a characteristic subcategory 
of $\ulEl(C)$, we conclude that $X(N(C))(\b \TTT_J(C)) = X\TTT_{J'}$ is in $\Ob\Hl$. Finally, let $I$ be an arbitrary $C$-resolving system. We have 
$X\TTT_I\lraisofl{30}{\al_{I,J'}} X\TTT_{J'}$ in $\ulEl(C)$, and thus $X\TTT_{J'}\in\Ob\Hl$ implies $X\TTT_I\in\Ob\Hl$, since a characteristic 
subcategory of $\ulEl(C)$ is, by definition, closed under isomorphy.\qed

\subsection{Pointwise pullback and pushout}
\label{SubsecPPP}

Suppose given an abelian category $\Al$, a poset $E$ and an element $\eps\in E$. Let $E^\eps := E\disj\{ \eps'\}$ be the poset defined by requiring that $\eps\le\eps'$, that $\al\not\le\eps'$
whenever $\al\not\le\eps$ and that $\eps'\not\le\al$ for all $\al\in E$; and the remaining relations within $E\tm E^\eps$ inherited from $E$. We define the {\it pushout at $\eps$}
\[
\barcl
\Al(E^\eps)             & \lraa{(-)\sind^{(=)}}  & \Al(E) \\
X'                      & \lramaps               & X\ind^{x'} \; ,\\
\ea
\]
where $X := X'|_E$, and $(X'_\eps\lraa{x'} X'_{\eps'}) = (X'_\eps\lrafl{30}{X'_{\eps'/\eps}} X'_{\eps'})$; and a transformation
\[
X'|_E = X\;\mra{i\;\=\; i X'}\; X\ind^{x'}\; ,
\]
natural in $X'$, by the following construction. Abbreviating $X\ind^{x'}$ by $\w X$, we let
\[
\xymatrix{
\hspace*{-9.6mm} X'_{\eps'} = \w X_\eps\ar[r]^{\w X_{\al/\eps}}  & \w X_\al \\
X_\eps \ar[r]^{X_{\al/\eps}}\ar[u]^{x'}\ar@{}[ur]|(0.7)\hookpo & X_\al\ar[u]_{i_\al} \\
}
\]
for $\al\in E$ with $\eps\le\al$. If $\eps\not\le\al$, we let $\w X_\al = X_\al$ and $i_\al = 1_{X_\al}$.

Given $\al\le\be$ in $E$, we let 
\[
\ba{rcll}
(\w X_\al\lraa{\w X_{\be/\al}} \w X_\be) &    & \mb{be induced by pushout}                          & \mb{if $\eps\le\al\le\be$}\; , \\
(\w X_\al\lraa{\w X_{\be/\al}} \w X_\be) & := & (X_\al\lraa{X_{\be/\al}}X_\be\lraa{i_\be} \w X_\be) & \mb{if $\eps\not\le\al$, but $\eps\le\be$}\; , \\
(\w X_\al\lraa{\w X_{\be/\al}} \w X_\be) & := & (X_\al\lraa{X_{\be/\al}} X_\be)                     & \mb{if $\eps\not\le\be$}\; . \\ 
\ea
\]
The morphism $X\lraa{i} X\ind^{x'}$ is the solution to the following universal problem. Suppose given a morphism $X\lraa{f} Y$ in $\Al(E)$ such 
that at $\eps\in E$ we have a factorisation 
\[
(X_\eps\lraa{f_\eps} Y_\eps) \= (X_\eps\lraa{x'} X'_{\eps'}\lra Y_\eps)\; . 
\]
Then there is a unique morphism $X\ind^{x'}\lraa{g} Y$ such that 
\[
(X\lraa{f} Y) \= (X\lraa{i} X\ind^{x'}\lraa{g} Y)\; .
\]

Dually, let $E_\eps := E\disj\{ \eps'\}$ be the poset defined by requiring that $\eps\ge\eps'$, that $\al\not\ge\eps'$
whenever $\al\not\ge\eps$ and that $\eps'\not\ge\al$ for all $\al\in E$; and the remaining relations within $E\tm E_\eps$ inherited from $E$. We define the {\it pullback at $\eps$}
\[
\barcl
\Al(E_\eps) & \lraa{(-)\sind_{(=)}}  & \Al(E) \\
X'          & \lramaps               & X\ind_{x'} \; ,\\
\ea
\]
where $X := X'|_E$, and $(X'_{\eps'}\lraa{x'} X'_\eps) = (X'_{\eps'}\lrafl{30}{X'_{\eps/\eps'}} X'_\eps)$; and a transformation
\[
X'|_E = X\;\mla{p\;\=\; pX'}\; X\ind_{x'}\; ,
\]
natural in $X'$, being the solution to the universal problem dual to the one described above.

\pagebreak

\subsection{$1$-epimorphic functors}
\label{Sec1epi}

Let $\Cl\lraa{F}\Dl$ be a functor between categories $\Cl$ and $\Dl$.

\begin{Definition}
\label{DefE1}\rm
The functor $\Cl\lraa{F}\Dl$ is {\it $1$-epimorphic} if the induced functor ``restriction along $F\;$''
\[
\fbo\Cl,\El\fbc\;\;\llaa{F(-)}\;\;\fbo\Dl,\El\fbc
\]
is full and faithful for any category $\El$. In particular, given functors
$\;
\Dl\;\lradoublea{G}{H\,}\;\El\vspace*{1mm}
\;$
with $FG\iso FH$, we can conclude that $G\iso H$; whence the notion of $1$-epimorphy.
\end{Definition}

\begin{Remark}
\label{LemE2}
Suppose given a diagram of categories and functors
\[
\xymatrix{
\Cl\ar[r]^F\ar[d]_S^\wr & \Dl\ar[d]^T_\wr \\
\Cl'\ar[r]^{F'}         & \Dl' \\
}
\]
with equivalences $S$ and $T$, and with $FT\iso SF'$. Then $F$ is $1$-epimorphic if and only if $F'$ is $1$-epimorphic.
\end{Remark}

Let $\,C,\, C'\,\in\,\Ob\Cl\,$. An {\it $F$-epizigzag} (resp.\ an {\it $F$-monozigzag}) $C\auf{u}{\leadsto} C'$ is a finite sequence of morphisms
\[
C \= C_0 \;\lraa{u_0}\; Z_0 \;\llaa{u'_0}\; C_1 \;\lraa{u_1}\; Z_1 \;\llaa{u'_1}\; C_1 \;\lraa{u_2}\; 
\cdots 
\;\llaa{u'_{k-2}}\; C_{k-1} \;\lraa{u_{k-1}}\; Z_{k-1} \;\llaa{u'_{k-1}} C_k \= C' 
\]
in $\Cl$ of {\it length} $k\ge 0$ such that $u'_i F$ is an isomorphism for all $i\in [0,k]$, and such that
\[
uF \;\; := \;\; (u_0 F)(u'_0 F)^-(u_1 F)(u'_1 F)^- \cdots (u_{k-1} F)(u'_{k-1} F)^-\; :\; CF\;\lra\; C'F
\] 
is a retraction (resp.\ a coretraction) in $\Dl$.

\begin{Lemma}
\label{LemE3}
Suppose the functor
\[
\Cl\;\;\lraa{F}\;\;\Dl
\]
to be dense, and to satisfy the following condition {\rm (C).}

{\rm (C)} $\hspace*{10mm}\left\{
\mb{\rm\begin{tabular}{p{12cm}}
Given objects $C,\, C'\in\Ob\,\Cl$ and a morphism $CF\lraa{d} C'F$ in $\Dl$, there exists an $F$-epizigzag $C_\text{s}\auf{c_\text{s}}{\leadsto} C$,
an $F$-monozigzag $C'\auf{c'_\text{t}}{\leadsto} C'_\text{t}$ and a morphism $C_\text{s}\lraa{c} C'_\text{t}$ such that
$$
(C_\text{s}F\lraa{c_\text{s}F} CF\lraa{d} C'F\lraa{c'_\text{t}F} C'_\text{t}F) \; =\; (C_\text{s}F\lraa{cF} C'_\text{t}F)\; .
$$
\end{tabular}}\right.$

Then $F$ is $1$-epimorphic.

\rm
{\it Proof.} Since $F$ is dense, Remark \ref{LemE2} allows to assume that $F$ is surjective on objects, i.e.\ $(\Ob\Cl)F = \Ob\Dl$. 

Let us prove that $\El(\Cl)\;\llaa{F(-)}\;\El(\Dl)$ is faithful. Suppose given functors $\Cl\lraa{F}\Dl\lradoublea{G}{H\,}\El\ru{-1}$ and morphisms 
$G\lraa{\ga} H$ and $G\lraa{\ga'} H$ such that $F \ga = F \ga'$. Given $D\in\Ob\Dl$, we have to show that $D \ga = D \ga'$. Writing $D = CF$ for some 
$C\in\Ob\Cl$, this follows from $D\ga = CF\ga = CF\ga' = D\ga'$.

Let us prove that $\El(\Cl)\;\llaa{F(-)}\;\El(\Dl)$ is full. Suppose given functors $\Cl\lraa{F}\Dl\lradoublea{G}{H\,}\El$ and a morphism 
$FG \lraa{\de} FH$. Define $G\lraa{\h\de} H$ by $(CF)\h\de := C\de$. 

We have to prove that $D\h\de$ is a welldefined morphism for $D\in\Ob\Dl$. So suppose that $D = CF = C'F$. We have to show that $C\de = C'\de$. 
By assumption (C), applied to $d = 1_D = 1_{CF} = 1_{C'F}$, there exist an $F$-epizigzag $C_\text{s}\auf{c_\text{s}}{\leadsto} C$, an 
$F$-monozigzag $C'\auf{c'_\text{t}}{\leadsto} C'_\text{t}$ and a morphism $C_\text{s}\lraa{c} C'_\text{t}$ such that $(c_\text{s} F)(c'_\text{t} F) = cF$. We obtain
\[
\barcl
(c_\text{s}FG)(C\de)(c'_\text{t} FH)
& = & (C_\text{s} \de)(c_\text{s} FH)(c'_\text{t} FH) \\
& = & (C_\text{s} \de)(c FH) \\
& = & (cFG)(C'_\text{t} \de) \\
& = & (c_\text{s}FG)(c'_\text{t} FG)(C'_\text{t} \de) \\
& = & (c_\text{s}FG)(C'\de)(c'_\text{t} FH)\; , \\
\ea
\]
whence $C\de  = C'\de$ by epimorphy of $c_\text{s}FG$ and by monomorphy of $c'_\text{t}FH$.

We have to prove that $\h\de$ is natural. Suppose given $CF\lraa{d} C'F$ in $\Dl$ for some $C,\, C'\in\Ob\Cl$. We have to show that $(dG)((C'F)\h\de) = ((CF)\h\de)(dH)$,
i.e.\ that $(dG)(C'\de) = (C\de)(dH)$. By  assumption (C), there exist an $F$-epizigzag $C_\text{s}\auf{c_\text{s}}{\leadsto} C$, an $F$-monozigzag
$C'\auf{c'_\text{t}}{\leadsto} C'_\text{t}$ and a morphism $C_\text{s}\lraa{c} C'_\text{t}$ such that $(c_\text{s} F) d (c'_\text{t} F) = cF$. We obtain
\[
\barcl
(c_\text{s} F G)(d G)(C' \de)(c'_\text{t} F H) 
& = & (c_\text{s} F G)(d G)(c'_\text{t} F G)(C'_\text{t} \de) \\
& = & (c G F)(C'_\text{t} \de) \\
& = & (C_\text{s} \de)(c F H) \\
& = & (C_\text{s} \de)(c_\text{s} FH)(dH)(C'_\text{t} F H) \\
& = & (c_\text{s} F G)(C \de)(d H)(c'_\text{t} F H)\; , \\
\ea
\]
whence $(d G)(C'\de) = (C \de)(d H)$ by epimorphy of $c_\text{s} F G$ and by monomorphy of $c'_\text{t} F H$.\qed
\end{Lemma}

\begin{Corollary}
\label{CorE3_5}
If $\Cl\lraa{F}\Dl$ is a functor such that {\rm (i,\, ii)} hold, then $F$ is $1$-epimorphic.

\begin{itemize}
\item[{\rm (i)}] For all morphisms $D\lraa{d} D'$ in $\Dl$, there is a morphism $C\lraa{c} C'$ in $\Cl$ such that 
\[
(C\lraa{c} C')F \= (D\lraa{d} D')\; .
\]
\item[{\rm (ii)}] For any $C,\, C'\,\in\,\Ob\Cl$ such that $CF = C'F$, there exists a finite sequence of morphisms
\[
C \= C_0 \;\lraa{u_0}\; Z_0 \;\llaa{u'_0}\; C_1 \;\lraa{u_1}\; Z_1 \;\llaa{u'_1}\; C_2 \;\lraa{u_2}\; 
\cdots 
\;\llaa{u'_{k-2}}\; C_{k-1} \;\lraa{u_{k-1}}\; Z_{k-1} \;\llaa{u'_{k-1}} C_k \= C' 
\]
from $C$ to $C'$ such that $u_i F = u'_i F = 1_{CF} = 1_{C'F}$ for all $i\in [0,k]$.
\end{itemize}
\end{Corollary}

{\it Proof.} The functor $F$ is dense, even surjective on objects, because identities have inverse images under $F$. To fulfill condition (C) of 
Lemma \ref{LemE3}, given objects $C,\, C'\,\in\,\Ob\Cl$ and a morphism $CF\lraa{d} C'F$ in $\Dl$, we may take some morphism
$C_\text{s}\lraa{c} C'_\text{t}$ in $\Cl$ such that $(C_\text{s}\lraa{c} C'_\text{t})F = (CF\lraa{d} C'F)$, we may take for $c_\text{s}$ a sequence as 
given by assumption because of $C_\text{s} F = CF$, and we may take for $c_\text{t}$ a sequence as given by assumption because of $C'_\text{t} F = C'F$.\qed

\begin{Corollary}
\label{CorE4}
If $\Cl\lraa{F}\Dl$ is a full and dense functor, then $F$ is $1$-epimorphic. 
\end{Corollary}

{\it Proof.} In fact, in condition (C) of Lemma \ref{LemE3}, we may take an $F$-monozigzag and an $F$-epizigzag of length $0.\;$\qed

\end{footnotesize}

%% file: h3ref.tex
\parskip0.0ex
\begin{footnotesize}

\parskip1.2ex

\vspace*{5mm}

\begin{flushright}
Matthias K\"unzer\\
Lehrstuhl D f\"ur Mathematik\\
RWTH Aachen\\
Templergraben 64\\
D-52062 Aachen \\
kuenzer@math.rwth-aachen.de \\
www.math.rwth-aachen.de/$\sim$kuenzer\\
\end{flushright}
\end{footnotesize}